\documentclass[11pt, final]{article}
\usepackage[english]{babel}
\usepackage[latin1]{inputenc}
\usepackage{exscale}
\usepackage[centertags]{amsmath}
\usepackage{amsfonts,amstext,amssymb,amsthm, amsmath}
\usepackage[usenames,dvipsnames]{xcolor}
\usepackage[colorlinks=true,urlcolor=blue,
citecolor=ForestGreen,linkcolor=blue,linktocpage,pdfpagelabels,
bookmarksnumbered,bookmarksopen,backref]{hyperref}
\usepackage{comment}
\usepackage{newlfont}
\usepackage{enumerate}
\usepackage{graphicx}
\usepackage{color}
\usepackage{url}
\usepackage{geometry}
\usepackage{float}
\usepackage{multirow}
\usepackage{graphicx}
\usepackage{amssymb}
\usepackage{verbatim} 
\usepackage{color}
\usepackage{epstopdf}
\usepackage{epsfig}
\usepackage[all]{xy}
\usepackage{mathrsfs}
\usepackage{stmaryrd}
\usepackage{caption}
\usepackage{subcaption}
\usepackage{xfrac}
\usepackage{aligned-overset}

\makeatletter
\def\@tocline#1#2#3#4#5#6#7{\relax
  \ifnum #1>\c@tocdepth 
  \else
    \par \addpenalty\@secpenalty\addvspace{#2}%
    \begingroup \hyphenpenalty\@M
    \@ifempty{#4}{%
      \@tempdima\csname r@tocindent\number#1\endcsname\relax
    }{%
      \@tempdima#4\relax
    }%
    \parindent\z@ \leftskip#3\relax \advance\leftskip\@tempdima\relax
    \rightskip\@pnumwidth plus4em \parfillskip-\@pnumwidth
    #5\leavevmode\hskip-\@tempdima
      \ifcase #1
       \or\or \hskip 1em \or \hskip 2em \else \hskip 3em \fi%
      #6\nobreak\relax
    \dotfill\hbox to\@pnumwidth{\@tocpagenum{#7}}\par
    \nobreak
    \endgroup
  \fi}
\makeatother
\usepackage{hyperref}

\DeclareMathAlphabet{\mathcalligra}{T1}{calligra}{a}{c}
\newcommand{\Cc}{\kern-4pt\mathcalligra{C\kern0.2pt a\kern0.2pt c\kern0.2pt c}\kern1pt}

\DeclareMathOperator{\spt}{spt}

\numberwithin{equation}{section}
\def\ttimes{\times{\kern-10.5pt}\times}

\newcommand\res{\mathop{\hbox{\vrule height 7pt width .3pt depth 0pt
\vrule height .3pt width 5pt depth 0pt}}\nolimits}

\newcommand{\gr}[1]{\operatorname{graph}(#1)}

\newcommand{\bG}{\mathbf{G}}

\newcommand{\sW}{{\mathscr{W}}}
\newcommand{\sC}{{\mathscr{C}}}
\newcommand\sS{{\mathscr S}}

\newcommand{\G}{\mathcal{G}}

\newcommand{\bL}{{\mathbf{L}}}

\newcommand{\bE}{{\mathbf{E}}}

\newcommand{\bB}{{\mathbf{B}}}
\newcommand{\bC}{{\mathbf{C}}}

\newcommand{\be}{{\mathbf{e}}}

\newcommand\N{{\mathbb N}}

\newcommand\R{{\mathbb R}}

\newcommand{\eps}{{\varepsilon}}
\newcommand{\bA}{\mathbf{A}}

\def\Xint#1{\mathchoice
{\XXint\displaystyle\textstyle{#1}}%
{\XXint\textstyle\scriptstyle{#1}}%
{\XXint\scriptstyle\scriptscriptstyle{#1}}%
{\XXint\scriptscriptstyle\scriptscriptstyle{#1}}%
\!\int}
\def\XXint#1#2#3{{\setbox0=\hbox{$#1{#2#3}{\int}$ }
\vcenter{\hbox{$#2#3$ }}\kern-.6\wd0}}
\def\mint{\Xint-}

\newcommand{\bp}{{\mathbf{p}}}

\newcommand{\bh}{{\mathbf{h}}}
\newcommand{\bM}{{\mathbf{M}}}

\newcommand{\diam}{{\operatorname{diam}}}

\def\Is#1{{\mathcal{A}}_{#1} (\R^{n})}

\def\a#1{\left\llbracket{#1}\right\rrbracket}
\newcommand{\abs}[1]{\left|#1\right|}

\newcommand{\D}{\textup{Dir}}

\newcommand{\breg}{\mathrm{Reg}_b}
\newcommand{\bsing}{\mathrm{Sing}_b}

\DeclareMathOperator{\dist}{dist}
\DeclareMathOperator{\spine}{spine}
\DeclareMathOperator{\Lip}{\operatorname{Lip}}
\DeclareMathOperator{\interior}{int}
\def\ttimes{\times{\kern-11.5pt}\times}
\def\restrict#1{\raise-.8ex\hbox{\ensuremath|}_{#1}}

\newcommand{\sigmaexp}{\sigma}

\newcommand{\infimo}[2]{\inf\left\{ #1 : #2 \right\}}

\newcommand{\euc}[1]{\mathbb{R}^{#1}}

\newcommand{\ball}[2]{\mathbf{B}\!\left(#1, #2\right)}

\newcommand{\baseball}[2]{\mathrm{B}_{#2}\left(#1\right)}
\newcommand{\basehalfball}[3]{\mathrm{B}_{#3}^{#1}\left(#2\right)}
\newcommand{\cyl}[2]{\mathbf{C}\!\left(#1, #2\right)}
\newcommand{\cyltilted}[3]{\mathbf{C}\!\left(#1, #2, #3\right)}

\newcommand{\cc}{\Subset}

\newcommand{\slice}[3]{\langle #1, #2, #3 \rangle}

\newcommand{\fdist}[2]{\mathbf{d}\left(#1, #2\right)}
\newcommand{\fdistin}[3]{\mathbf{d}_{#1}\left(#2, #3\right)}
\newcommand{\massnorm}[1]{\mathbb{M}\left(#1\right)}

\counterwithin{equation}{section}

\theoremstyle{definition}
\newtheorem{definition}{Definition}[section]

\newtheorem{remark}[definition]{Remark}
\newtheorem{assumption}{Assumption}
\theoremstyle{plain}
\newtheorem{theorem}[definition]{Theorem}
\newtheorem{proposition}[definition]{Proposition}
\newtheorem{corollary}[definition]{Corollary}
\newtheorem{lemma}[definition]{Lemma}
\newtheorem*{result}{Theorem}

\title{\bf{Density of the boundary regular set of $2d$ area minimizing currents with arbitrary codimension and multiplicity}}
\newcommand\blfootnote[1]{%
  \begingroup
  \renewcommand\thefootnote{}\footnote{#1}%
  \addtocounter{footnote}{-1}%
  \endgroup
}

\begin{document}

\date{}

\maketitle

\author{
\begin{center}
\begin{tabular}{c c c}
 Stefano Nardulli & Reinaldo Resende* \\ 
 \textit{Universidade Federal do ABC} & \textit{Carnegie Mellon University}\\
 \textit{stefano.nardulli@ufabc.edu.br} & \textit{rresende@andrew.cmu.edu}
\end{tabular}
\end{center}
}\blfootnote{*Corresponding author}
\bigskip

\begin{abstract}
In the present work, we consider area minimizing currents in the general setting of arbitrary codimension and arbitrary boundary multiplicity. We study the boundary regularity of $2d$ area minimizing currents, beyond that, several results are stated in the more general context of $(C_0, \alpha_0,r_0)$-almost area minimizing currents of arbitrary dimension $m$ and arbitrary codimension taking the boundary with arbitrary multiplicity. Furthermore, we do not consider any type of convex barrier assumption on the boundary, in our main regularity result which states that the regular set, which includes one-sided and two-sided points, of any $2d$ area minimizing current $T$ is an open dense set in the boundary. 
\end{abstract}

\setcounter{tocdepth}{2}
\tableofcontents

\section{Introduction}

\subsection{Historical overview}

The regularity of area-minimizing currents has been extensively investigated by numerous eminent mathematicians, leading to the emergence of significant findings pertaining to both interior and boundary regularity. Using the framework of integral currents to study the problem of minimizing the area functional is very adequate thanks to the foundational existence theory developed by Federer and Fleming in their celebrated paper \cite{FF}. 

\subsubsection{Regularity results in codimension 1}

Consider $T$ to be an $m$-dimensional integral current that minimizes area in $\R^{m+1}$. The regularity theory within this context is well-established. In fact, due to the effort of several mathematicians, namely De Giorgi, Simons, Federer, Almgren, and Fleming, it has been proven that the \emph{interior} singular set of $T$ has Hausdorff dimension at most $m-7$. Moreover, this dimensional bound is proven to be sharp, it can be verified with the so-called Simons' cone (introduced by Simons in \cite{Simons} where its stationarity and stability are established, while the proof of its minimality property were proven in the celebrated \cite{BDeGG} by Bombieri, De Giorgi, and Giusti) which is a $7$-dimensional area minimizing current in $\R^8$ such that the \emph{interior} singular set is a singleton.

The boundary regularity in codimension $1$ is as well completely developed. In fact, if $T$ has a contour $\Gamma\in C^{1,\alpha}, \alpha\in (0,1)$, then we have full regularity for $T$ at the boundary, namely, the \emph{boundary} singular set is empty. This result is proven by Hardt and Simon in \cite{HS} for multiplicity $1$ boundaries. In codimension $1$, it is not a very difficult problem to extend it for higher multiplicity boundaries though. Indeed, White in \cite{white1983regularity} provided an elegant and insightful decomposition argument that enables us to consistently reduce the problem to the case of multiplicity 1. Consequently, regardless of its boundary multiplicity, $T$ does not have \emph{boundary} singularities in codimension $1$. The generalization of \cite{HS} to integral currents minimizing area in a Riemannian submanifold of some Euclidean space, is due to Steinbruechel in \cite{steinbruechel2020boundary}.

\subsubsection{Interior regularity in high codimension}

The results above rely on several powerful features that are exclusive to the case of codimension $1$. A comprehensive elucidation of these properties can be found in \cite{delellis2015size}.

We now consider $T$ to be an $m$-dimensional integral current that minimizes area in $\Sigma\subset\R^{m+n}$ with $n>1$ and $\Sigma$ is a $(m+l)$-dimensional submanifold of $\R^{m+n}$. In this setting, more involved objects can emerge, for instance, the famous Federer's examples of complex varieties, \cite[Section 5.4.19]{Fed}. This example consists of an $2$-dimensional area minimizing current in $\R^4$ with a single \emph{interior} singular point. Notably, the dimensional bound proven in codimension $1$ dramatically fails. Nevertheless, Almgren's masterpiece \cite{Alm} proves that, if $\Sigma$ is $C^5$, the Hausdorff dimension of the \emph{interior} singular set of $T$ cannot exceed $m-2$ which, in view of Federer's example, is an optimal result. 

While Almgren's program offers profound insights and powerful tools, it is extremely long and intricate. Aiming at simplifying the proofs, De Lellis and Spadaro, in a series of papers \cite{DS1,DS2,DS3,DS4,DS5}, gave a shorter and significantly simpler proof for Almgren's result. Furthermore, they extended to encompass cases where $\Sigma\in C^{3,\alpha}$. 

Before walking through the results for the boundary regularity in high codimension, we set the notation and basic definitions that will be used throughout this manuscript.

\subsection{Definitions and notation}

For basic definitions and standard notations, we refer the reader to the textbooks \cite{Fed} and \cite{simon2014introduction}. Let us set up some notation that will be used in this work. 

We denote by $\mathscr{R}_m^{loc}(U)$ the space of $m$-dimensional integer rectifiable currents in $U\subset\R^{m+n}$, $\mathbf{I}_m(U)$ the space of $m$-dimensional integral currents in $U\subset\R^{m+n}$. We define the $m$-density of a current $T$ at $p\in\R^{m+n}$ as follows:
$$\Theta^m(T,p) := \lim_{r\to 0}\frac{\|T\|(\ball{p}{r})}{\omega_m r^m},$$
where $\omega_m$ is the $m$-dimensional Hausdorff measure of a $m$-dimensional ball with radius $1$ and $\|\cdot\|$ denotes the mass of the current $T$.

We also fix the notation of the flat distance between two $m$-dimensional integer rectifiable currents $T$ and $S$, i.e., $T, S\in\mathscr{R}_m^{loc}(U)$, $U$ open and $A\cc U$ as follows:
$$\fdistin{A}{T}{S} = \infimo{\|R\|(A)+\|\tilde{T}\|(A)}{T-S=R+\partial\tilde{T} \ \text{with} \ R\in\textbf{I}_m(U) \ \text{and} \ \tilde{T}\in\textbf{I}_{m+1}(U)}.$$

Given a $m$-rectifiable set $M\subset\R^{m+n}$, it naturally induces an $m$-current which we always denote by $\a{M}$. We say that the boundary of an $m$-current $T$ is taken with multiplicity $Q^{\star}$ if $T=Q^{\star}\a{\Gamma}$ for some $(m-1)$-rectifiable set $\Gamma$. 

\begin{definition}\label{def:almot-min}
Given three real numbers $C_{0}\ge0, r_{0}, \alpha_{0}>0$, we say that an $m$-dimensional integer rectifiable current $T \ ($i.e. $T\in\mathscr{R}_m^{loc}(\R^{m+n}))$ with $\partial T={Q^\star}\a{\Gamma}, {Q^\star}\in\N\setminus\{0\},$ is $(C_0, r_0, \alpha_0)$-\textbf{almost area minimizing at $x \in \operatorname{spt}(T)$}, if we have
\begin{equation}\label{Eq:AlmostMinimality}
 \|T\|\left(\ball{x}{r}\right) \leq\left(1+C_{0} r^{\alpha_{0}}\right)\|T+\partial\tilde{T}\|\left(\ball{x}{r}\right),
\end{equation}
for all $0<r<r_{0}$ and all integral $(m+1)$-dimensional currents $\tilde{T}$ supported in $\ball{x}{r}$, i.e., for all $\tilde{T}\in\mathbf{I}_{m+1}(\ball{x}{r})$. The current is called $(C_{0}, r_{0}, \alpha_{0})$-\textbf{almost area minimizing in $\R^{m+n}$}, if the current $T$ is $(C_{0}, r_{0}, \alpha_{0})$-almost area minimizing at each $x \in \operatorname{spt}(T)$ with the same constants $C_{0}, r_{0}, \alpha_{0}>0$ independently of the point $x$. If $T$ is $(0, r_0, \alpha_0)$-almost area minimizing, we say that $T$ is \textbf{area minimizing}.
\end{definition}

We now define the right notions of regular boundary points in the arbitrary boundary multiplicity setting. 

\begin{definition}[Regular and singular one-sided boundary points, Definition 0.1 of \cite{DNS}]\label{defi:one-sided-reg-sing-points}
Let $T$ be a $2$-dimensional integer rectifiable current $($i.e. $T\in\mathscr{R}_2^{loc}(\R^{2+n}))$ with $\partial T={Q^\star}\a{\Gamma}, {Q^\star}\in\N\setminus\{0\}, p\in\Gamma$ and $\Theta^2(T,p)=\frac{Q^\star}{2}$. Then $p$ is called a {\textbf{regular one-sided boundary point}} if $T$ consists, in a neighborhood $U$ of $p$, of the union of finitely many surfaces with boundary $\Gamma$, counted with multiplicities, which meet at $\Gamma$ transversally. More precisely, if there are:
\begin{enumerate}[(i):]
    \item a finite number $\Sigma_{1}, \ldots, \Sigma_{J}$ of oriented embedded surfaces in $U$,
    
    \item and a finite number of positive integers $k_{1}, \ldots, k_{J}$ such that:
        \begin{enumerate}
            \item $\partial \Sigma_{j} \cap U=\Gamma \cap U=\Gamma_{i} \cap U$ (in the sense of differential topology) for every $j$,
            
            \item $\Sigma_{j} \cap \Sigma_{l}=\Gamma \cap U$ for every $j \neq l$,
            
            \item for all $j \neq l$ and at each $q \in \Gamma$ the tangent cones to $\Lambda_{j}$ and $\Lambda_{l}$ are distinct,
            
            \item $T \res U=\sum_{j} k_{j} \a{\Sigma_{j}}$ and $\sum_{j} k_{j}=Q^\star$.            
        \end{enumerate}
\end{enumerate}
The set $\operatorname{Reg}_{b}^1(T)$ of {\textbf{regular boundary one-sided points}} is a relatively open subset of $\Gamma$.
\end{definition}

\begin{figure}[H]
\centering
\captionsetup{width=.7\linewidth}\includegraphics[width=.7\linewidth]{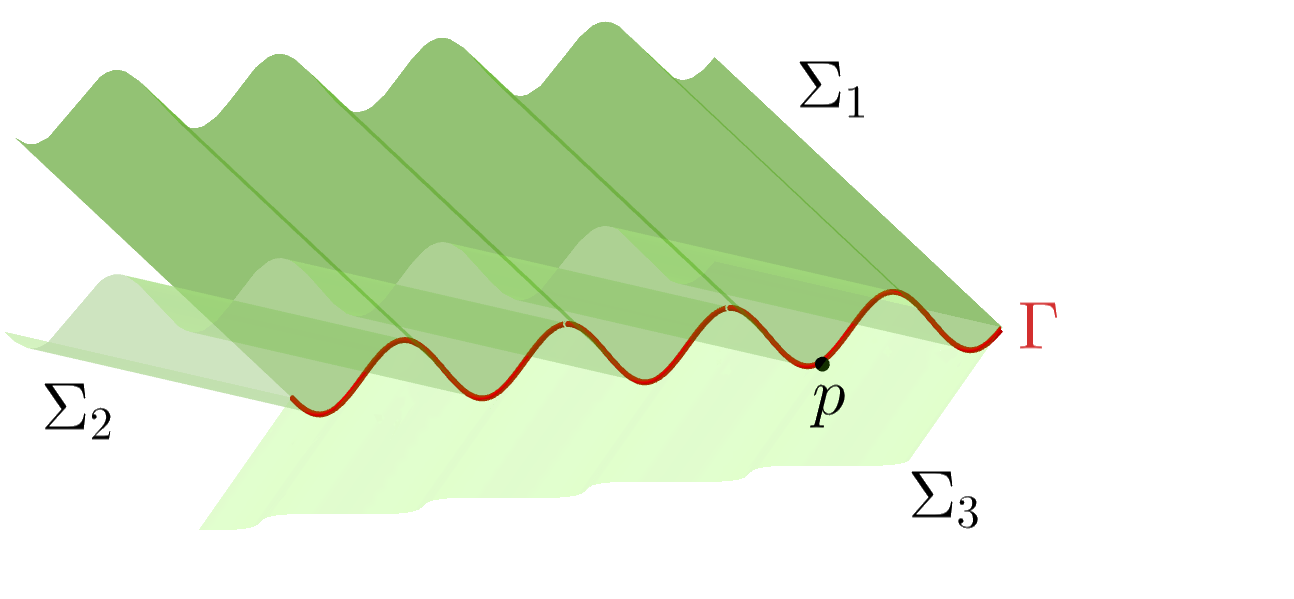}
\caption{Here $J=3$ and the current is given by $T = \sum_{j=1}^{3}k_j\a{\Sigma_j}$, then $p$ is a regular one-sided boundary point of $T$. Note that, each surface $\Sigma_j$ is taken with an integer multiplicity $k_j$ and the boundary $\partial T$ has multiplicity $Q^\star = k_1+k_2+k_3$.}
\label{}
\end{figure}

\begin{definition}[Regular and singular two-sided boundary points, Definition 1.1 of \cite{DDHM}]\label{defi:two-sided-reg-sing-points}
Let $T$ be a $m$-dimensional integer rectifiable current $($i.e. $T\in\mathscr{R}_2^{loc}(\R^{2+n}))$ with $\partial T={Q^\star}\a{\Gamma}, {Q^\star}\in\N\setminus\{0\}, p\in\Gamma$ and $\Theta^2(T,p)>\frac{Q^\star}{2}$. 

\begin{enumerate}[(i):]
    \item We say that $p$ is a {\textbf{regular boundary two-sided point for $T$}} if there exist a neighborhood $U \ni p$ and a surface $\Sigma \subset U \cap \Sigma$ such that $\operatorname{spt}(T) \cap U \subset \Sigma$. The set of such points will be denoted by $\operatorname{Reg}_{\mathrm{b}}^2(T)$,
    
    \item We also denote $\operatorname{Reg}_b(T):=\operatorname{Reg}_{\mathrm{b}}^1(T)\mathring{\cup}\operatorname{Reg}_{\mathrm{b}}^2(T)$, $\operatorname{Sing}_b(T):= \Gamma\setminus (\operatorname{Reg}_{\mathrm{b}}^1(T)\mathring{\cup}\operatorname{Reg}_{\mathrm{b}}^2(T))$

    \item We will say that $p \in \operatorname{Sing}_{\mathrm{b}}(T)$ is of {\textbf{crossing type}} if there is a neighborhood $U$ of $p$ and two currents $T_{1}$ and $T_{2}$ in $U$ with the properties that:
        \begin{enumerate}
            \item $T_{1}+T_{2}=T$ and $\partial T_{1}=0$,
            
            \item $p \in \operatorname{Reg}_{\mathrm{b}}\left(T_{2}\right)$.
        \end{enumerate} 

    \item If $p \in \operatorname{Sing}_{\mathrm{b}}(T)$ is not of crossing type, we will then say that $p$ is a {\textbf{genuine boundary singularity point of $T$}}.
\end{enumerate}
\end{definition} 

\begin{figure}[H]
\centering
\captionsetup{width=.7\linewidth}\includegraphics[width=.5\linewidth]{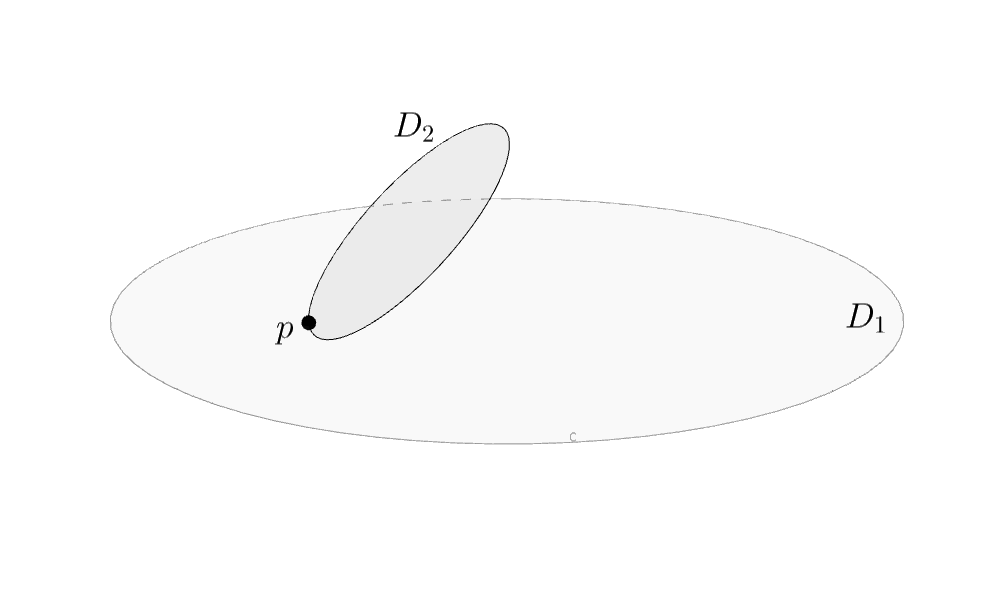}
\caption{Let $T=\a{D_1}+\a{D_2}$ and $p\in\partial D_2\cap\operatorname{int}(D_1)$. It is easy to see that $p$ is a crossing type singularity to the $2d$ current $T$.}
\label{}
\end{figure}

\subsection{Boundary regularity in high codimension}

\subsubsection{Boundary with multiplicity 1}

In \cite{AllB}, Allard proved that every \emph{boundary} point of an $m$-dimensional area minimizing current $T$ in $\R^{m+n}$ taking the boundary with multiplicity $1$, i.e., $\partial T = \a{\Gamma}$, $\Gamma\in C^{1,1}$, is regular, however, he needed to impose a crucial condition on $\operatorname{spt}\left(\partial T\right)$ which is that $\operatorname{spt}(\partial T)$ is contained in the boundary of a uniformly convex set, we call such condition \emph{convex barrier}. Indeed, what Allard proved is that boundary points with density close to $1/2$ are regular, afterwards he showed that every boundary point with the convex barrier assumption has density close to $1/2$ and hence is regular. 

Since Allard's result which concerns \emph{one-sided} points, no one successfully attacked the problem of removing the convex barrier condition (i.e., considering \emph{two-sided} points) until the recent work of De Lellis, De Philippis, Hirsch, and Massaccesi, \cite{DDHM}. In this article, the authors prove the following result:

\begin{result}[\cite{DDHM}]
Let $T$ be an $m$-dimensional area minimizing current in $\Sigma\subset\R^{m+n}$ with $\partial T = \a{\Gamma}$, $\Sigma$ is a $C^{3,\alpha}$-submanifold, and $\Gamma\subset\Sigma$ of class $C^{3,\alpha}$, then $\breg(T)$ is an open dense set in $\Gamma$. 
\end{result}

We cannot hope to prove a Hausdorff dimension bound for $\operatorname{Sing}_b(T)$, due to the fact that there is a $2$-dimensional area minimizing current in $\R^4$ with $\partial T = \a{\Gamma}$, $\Gamma$ smooth, and the \emph{boundary} singular set has Hausdorff dimension equal to $1$, see \cite[Theorem 1.8]{DDHM}. There are some further conjectures on dimensional bounds for crossing type and genuine singularities of $T$, we refer the reader to \cite{DDHM}.

To achieve the regularity result above, the authors in \cite{DDHM} have tailored the framework introduced by Almgren, and sharpened by De Lellis and Spadaro, to the boundary case which is even more involved and needs the introduction of a lot of highly nontrivial new ideas. Moreover, the authors also rely on Allard's result to be able to use that one-sided points are regular.

\subsubsection{Boundary with multiplicity \texorpdfstring{$Q^{\star}$}{Lg} and arbitrary codimension}

In codimension bigger than $1$, White's decomposition argument (\cite{white1983regularity}) does not work. It is then a more difficult task to pass from boundaries with multiplicity one to boundaries with higher multiplicity. In fact, very recently, De Lellis, Nardulli, and Steinbruechel, in \cite{DNS} proved the following result:

\begin{result}[\cite{DNS}]
Let $T$ be a $2$-dimensional area minimizing current in $\R^{2+n}$ with $\partial T = Q^{\star}\a{\Gamma}$, $\Gamma$ is a $C^{3,\alpha}$ arc. Then every point $p\in\Gamma$ with $\Theta^2(T,p) < \frac{Q^{\star}+1}{2}$ is a regular point for $T$.
\end{result}

In this result, the authors also modified the framework introduced by Almgren, and sharpened by De Lellis and Spadaro, to deal with one-sided points. It considerably differs from the techniques used in \cite{DDHM} suited for two-sided points.

It is important to mention that the restriction for $2$-dimensional currents is due to the fact that $2d$ area minimizing cones with arbitrary boundary multiplicity and codimension are classified in \cite{DNS2}, very little is known about tangent cones with arbitrary boundary multiplicity, dimension and codimension though.

Our main theorem will be proved under the following assumptions.

\begin{assumption}\label{assump:intro}
Let $\alpha \in (0,1]$ and an integer $Q^\star\geq 1$. Consider $\Gamma \subset \R^{2+n}$ a $C^{3,\alpha}$ oriented curve without boundary. Let $T$ be an integral $2$-dimensional area minimizing current in $\ball{0}{2}\subset\R^{2+n}$ with boundary $\partial T \res \ball{0}{2}={Q^\star}\a{\Gamma \cap \ball{0}{2}}$.
\end{assumption}

We can now state the main theorem which gives the density of the regular set $\breg(T)$ in $\Gamma$, where the regular set allows the existence of both one-sided and two-sided points.

\begin{theorem}\label{thm:main-theorem}
Let $T$ and $\Gamma$ be as in Assumptions \ref{assump:intro}. Then $\breg(T)$ is an open dense set in $\Gamma$.
\end{theorem}

The main theorem of the present work fills one of the gaps in the theory providing a boundary regularity theorem to $2d$ area minimizing currents with arbitrary boundary multiplicity, arbitrary codimension, and without any convex barrier assumption. To prove it, we rely on the regularity result for one-sided points (\cite{DNS}) in the same fashion that the authors in \cite{DDHM} relied on \cite{AllB} for the regularity of one-sided points.

A crucial tool for this manuscript is the characterization of $2d$ area minimizing cones given in \cite{DNS2} which allows us to proceed in various steps below. In fact, this characterization is strongly used to perform a stratification argument that is part of the proof of the density of $\operatorname{Reg}_b(T)$. This is thoroughly discussed and carefully proved in Section \ref{s:stratification}. This stratification argument will reduce the proof of Theorem \ref{thm:main-theorem} to proving the following statement:

\begin{result}[Theorem \ref{thm:collapsed-is-regular}]
Let $T$ and $\Gamma$ be as in Assumption \ref{assump:intro}. Then any two-sided collapsed point of $T$ is a two-sided regular point of $T$.
\end{result}
\begin{remark}
We actually prove Theorem \ref{thm:collapsed-is-regular} for $T$ of any dimension $m$. 
\end{remark}

Two-sided collapsed points (see Definition \ref{defi:collapsed}) are points of minimum density in some of their neighborhoods, and at these points, $T$ blows up to a cone whose support is an $m$-plane.

\subsection{Outline of the proof of Theorem \ref{thm:main-theorem}}

Our primary concern is the existence of two kinds of points that behave completely differently. These points are defined as one-sided points and two-sided points. It means that the current can be decomposed into pieces which are contained in half spaces of the ambient space in the one-sided case (see \cite[Thm 4.3 and 4.4]{DNS2}) and, in the case of two-sided points, the current is not contained in any half-space of the ambient space. 

In \cite{DDHM}, the authors consider $Q^\star = 1$ and thus $Q-\frac{1}{2}$ is the density of the two-sided regular points in the boundary. They also allow for the existence of both one-sided and two-sided points. However, their concept of regularity aligns with the definition provided by Allard (in the case where $Q = 1$) for one-sided points. Allard's notion involves what is termed a \emph{convex barrier}. This equivalence in definitions critically relies on the fact that the boundary multiplicity is $1$. If $Q^\star >1$, as we have defined, the notions of regularity for one-sided and two-sided points given in Definitions \ref{defi:one-sided-reg-sing-points} and \ref{defi:two-sided-reg-sing-points} do not coincide. This discrepancy arises specifically due to the presence of higher multiplicity settings at one-sided points, where we encounter regular currents with multiple sheets which are regular. Such currents are regular and Definition \ref{defi:two-sided-reg-sing-points} (when $Q=Q^\star$) does not cover them. Similarly to how the authors in \cite{DDHM} rely on Allard's regularity result to reduce the analysis to two-sided points, we employ \cite{DNS} to focus our analysis on two-sided points. In \cite{DNS}, the authors have proved the analog of Allard's theorem for higher multiplicity and currents of dimension $2$, i.e., assuming that $\Gamma$ belongs to a \emph{convex barrier}, $2d$ area minimizing currents with $T=Q^\star\a{\Gamma}$ are completely regular at the boundary. This reduction is allowed by Theorem \ref{theorem:not-dense-collapsed-singular} where we state, for $2$-dimensional $(C_0, r_0,\alpha_0)$-almost area minimizing currents, that, if every two-sided collapsed point (Definition \ref{defi:collapsed}) is regular, then $\breg(T)$ is dense. It is essential to emphasize that $\breg(T)$ encompasses both one-sided and two-sided points.

The principal objective now becomes establishing the regularity of two-sided collapsed points, Theorem \ref{thm:collapsed-is-regular}. To accomplish this, we follow the framework presented in \cite{DDHM}. Furthermore, this proof is done for area minimizing currents with arbitrary dimension $m$, codimension $n$ and arbitrary multiplicity $Q^\star$. Initially, we construct a linear theory for the pair $(f^+, f^-)$ where we define these pairs as $(Q-\frac{Q^\star}{2})$-valued functions (Definition \ref{def:function-interface}). We investigate the regularity that can be inferred under the assumption that this pair minimizes the Dirichlet energy. Given $\Omega\subset\R^m$ an open set and $\gamma$ a $(m-1)$-submanifold, called \emph{interface}, which splits $\Omega$ into $\Omega^+$ and $\Omega^-$, we define $(Q-\frac{Q^\star}{2})$-valued functions as a $Q$-valued function $f^+$ defined in $\Omega^+$ and a $(Q-Q^\star)$-valued function $f^-$ defined in $\Omega^-$ in the sense of Almgren. When these functions glue at $\gamma$ we say that $(f^+, f^-)$ collapses at the interface as in Definition \ref{def:function-interface}. If we assume that $(f^+, f^-)$ is a $(Q-\frac{Q^\star}{2})$-Dir minimizer which collapses at the interface, then it is given by $Q$ copies of $\kappa$ in $\Omega^+$ and $(Q-Q^\star)$ copies of $\kappa$ in $\Omega^-$ where $\kappa$ is a classical harmonic function, Theorem \ref{thm:harmonic-regularity-(Q-Q^*/2)-minimizing}.

The subsequent step involves approximating the current $T$ using Lipschitz maps. This approximation can be performed without presuming any minimizing property of the current, Theorem \ref{thm:first-lip-approx-and-good-set}. Furthermore, when we add the condition that $T$ is area minimizing, the approximation becomes minimizers of the Dirichlet energy and thus, by Theorem \ref{thm:harmonic-regularity-(Q-Q^*/2)-minimizing}, we obtain an harmonic approximation of the current $T$, Theorem \ref{thm:first-lip-approx-by-harmonic-sheet}. This harmonic approximation serves as a pivotal element for establishing a milder decay of excess for the area-minimizing current $T$, as indicated in Lemma \ref{lemma:milder-excess-decay}. Subsequently, we employ a straightforward argument to achieve a superlinear decay in the excess of $T$, an achievement outlined in Theorem \ref{thm:improved-excess-decay}. The superlinear decay is an important tool to prove the uniqueness of tangent cones at two-sided collapsed points of $T$ (Theorem \ref{thm:uniqueness-dimension-m}). Additionally, it enhances our Lipschitz approximation by providing better estimates, i.e. superlinear estimates, on the errors of the approximation, Theorem \ref{thm:lip-approx-superlinear-decay}. These constructions provide the foundation for the construction of the center manifolds.

Aiming at the construction of the center manifold, we first construct the Whitney decomposition with suitable stopping conditions paving the way for the establishment of the so-called $C^{3,\kappa}$ center manifolds $\mathcal{M}^+$ and $\mathcal{M}^-$, see Theorem \ref{thm:center_manifold}. Afterwards, in Theorem \ref{thm:center-manifold-app}, we introduce the $\mathcal{M}$-normal approximation. This approximation consists of multivalued Lipschitz maps $\mathcal{N}^+$ and $\mathcal{N}^-$ defined on the center manifold $\mathcal{M}^+$ and $\mathcal{M}^-$, respectively, taking values in the normal bundle of $\mathcal{M}$ which approximate the $m$-current $T$ in the desired manner. 

Subsequently, we present a blowup argument in Subsection \ref{blowup}, which serves to establish that $\mathcal{N}^\pm \equiv 0$. As a consequence, the $m$-dimensional area minimizing current $T$ has to coincide with $\mathcal{M}^+$ in the right portion and with $\mathcal{M}^-$ in the left portion. This outcome establishes that any two-sided collapsed point is in fact a two-sided regular point, which is formalized as Theorem \ref{thm:collapsed-is-regular}. This conclusion finishes the proof of our main theorem as aforementioned, i.e., Theorem \ref{thm:main-theorem}.

\subsection{Fundamental results}

We firstly define the blowups.

\begin{definition} We define for $x\in\euc{m+n}$ and $r>0$ the function $\iota_{x, r}(y):= \frac{y-x}{r}$. For any current $T\in\mathcal{D}_m^\prime(\R^{m+n})$ let us define the \textbf{rescaled current $T$ at $x$ at scale $r$} as ${\iota_{x, r}}_{\sharp} T=:T_{x,r}$ and $T_r =: {\iota_{0,r}}_{\sharp} T$. We call a current $T_{x}$ a \textbf{blowup of $T$ at $x$}, if there exists a sequence of radii $r_j\to0$ such that $T_{x,r_j}\to T_x$ in the weak topology. 
\end{definition} 

We now set the assumptions with the almost minimality condition for which we can prove the stratification theorems in Section \ref{s:stratification}. 

\begin{assumption}\label{assumptions}
Let $\left.\left. \alpha \in \right]0,1\right]$ and integers $m\geq 2$, $Q^\star\geq 1$. Consider $\Gamma \subset \R^{m+n}$ a $C^{3,\alpha}$ oriented $(m-1)$-submanifold without boundary. Let $T$ be an integral $m$-dimensional $(C_0, r_0, \alpha_0)$-almost area minimizing current in $\ball{0}{2}$ with boundary $\partial T \res \ball{0}{2}={Q^\star}\a{\Gamma \cap \ball{0}{2}}$.
\end{assumption}

Let us enunciate the almost monotonicity formula which will often used in what follows. The proof of this formula can be found in \cite[Lemma 2.1]{hirsch2019uniqueness}. Although, it is done for boundary multiplicity one currents, the proof can be readily adapted to the higher multiplicity case as it is done in \cite{DNS2}. 

\begin{proposition}[Almost monotonicity formula for boundary points, Proposition 3.3 of \cite{DNS2}]\label{prop:almost-monotonicity}
Let $T$ and $\Gamma$ be as in Assumption \ref{assumptions}, and $x\in\operatorname{spt}{T}\cap\Gamma$. Set $\alpha_1:=\min\{\alpha_0,\alpha\}$, $0<r_1<\min\{r_0, r'\}$. Then there is a constant $C_1=C_1(m, n, C_0, r_0, r', \alpha_0, \alpha, \theta, \| \Gamma \|_{1,\alpha})>0$, such that
\begin{equation}
e^{C_1 r^{\alpha_1}} \frac{\|T\|\left(\ball{x}{r}\right)}{r^{m}}-e^{C_1 s^{\alpha_1}} \frac{\|T\|\left(\ball{x}{s}\right)}{s^{m}} \geq \int_{\ball{x}{r} \setminus \ball{x}{s}} e^{C_1|z-x|^{\alpha_1}} \frac{\left|(z-x)^{\perp}\right|^{2}}{2|z-x|^{m+2}} d\|T\|(z), 
\end{equation}
for every $0<s<r<r_1$.
\end{proposition}

The upper semicontinuity of the density function is well known when restricted to either the interior or the boundary of an area minimizing current, we give a short proof of the validity of this fact at the boundary of almost area minimizing currents. We would like to remark that the upper semicontinuity holds for the restriction of the density function to boundary or interior points, however, it does not hold when it is considered defined on the whole $\spt(T)$, therefore, to get around that we state (iii) as B. White does in the context of area minimizing currents, c.f. \cite{White97}.

\begin{proposition}[Upper semicontinuity of the density function]\label{prop:upper-semicont-density}
Let $\Gamma$ and $T$ be as in Assumption \ref{assumptions}. Then
\begin{enumerate}[\upshape (i)]
    \item\label{item:prop:upper-semicont-density:bordo} The function $x\mapsto \Theta^m(T,x)$ is upper semicontinuous in $\spt(T)\cap\Gamma$,
    
    \item\label{item:prop:upper-semicont-density:interior} The function $x\mapsto \Theta^m(T,x)$ is upper semicontinuous in $\spt(T)\setminus\Gamma$,
    
    \item The function $x\mapsto \left\{\begin{matrix}\begin{aligned}
&\Theta^m(T,x), &x\notin\Gamma \\ 
&2\Theta^m(T,x), &x\in\Gamma
\end{aligned}\end{matrix}\right.$ is upper semicontinuous.
\end{enumerate}
\end{proposition}
\begin{proof}
We define $f(x, r) = e^{C_1r^{\alpha_1}}r^{-m}\|T\|(\ball{x}{r})$ and take $x_i \to x, x_i \in \spt(T)\cap\Gamma,$ and assume $\ball{x_i}{r}\subset\ball{x}{r+\varepsilon}, \forall i\in\N, \varepsilon>0.$ Applying the almost monotonocity formula (Proposition \ref{prop:almost-monotonicity}), we already know that $f(x_i,\cdot)$ is monotone nondecreasing and so we obtain 
\begin{align}
    f(x_i, t) &\leq f(x_i, r)\le  e^{C_1r^{\alpha_1}}r^{-m}\|T\|(\ball{x}{r+\varepsilon}) = f(x,r+\varepsilon)\overbrace{\left[e^{-C_1\varepsilon^{\alpha_1}} \left( 1 + \frac{\varepsilon}{r} \right)^m \right]}^{I} ,\nonumber
\end{align}
for any $0<t<r$, $i\ge i_\varepsilon$ and $\varepsilon>0$. First let $t\to 0$ we get that for every fixed positive $\varepsilon>0$  it holds
\begin{align}
    f(x_i, 0^+)= \Theta^m(T, x_i)&\le f(x, r+\varepsilon)\overbrace{\left[e^{-C_1\varepsilon^{\alpha_1}} \left( 1 + \frac{\varepsilon}{r} \right)^m \right]}^{I} ,\forall i\ge i_\varepsilon.\nonumber
\end{align}
In the last inequality, taking in first the limit with respect to $i$, in second with respect to $\varepsilon$ and finally with respect to $r$ leads to  $$\limsup_{i\to +\infty}\Theta^m(T, x_i)\leq\lim_{r\to 0}\limsup_{i\to +\infty}f(x,r) = \Theta^m(T,x).$$
We mention that to prove (ii), the upper semicontinuity at the interior, one may use this very same argument but now using the almost monotonicity formula given in \cite[Proposition 2.1]{DSS1}. Let us turn to the proof of (iii), it is enough to prove that 
$$ 2\Theta(T,x) \geq {\lim\sup}_{i\to+\infty}\Theta(T, x_i), $$
where $\{x_i\}_{i\in\N}\subset\R^{m+n}\setminus\Gamma$ converges to $x\in\Gamma$. Fix $y_i\in\Gamma$ to be the nearest point to $x_i$ and then define $T_i^\prime = \frac{T-y_i}{|x_i-y_i|}$ and $x_i^\prime = \frac{x_i-y_i}{|x_i-y_i|}$. Up to subsequences, we can set $x^\prime := \lim x_i^\prime$ and $T^\prime := \lim T_i^\prime$, notice that $T^{\prime}$ exists thanks to classical compactness theorems for integral currents, e.g. \cite[Theorem 4.2.17]{Fed}. By the almost monotonocity formula (Proposition \ref{prop:almost-monotonicity}), for every $r\in(0, r_1)$, the quantity $e^{C_1r^{\alpha_1}}\frac{\|T\|(\ball{y_i}{r})}{r^m}$ is monotone nondecreasing. Thus, for any $\rho\in\R_{+}$, 

$$\begin{aligned}
e^{C_1r^{\alpha_1}}\frac{\|T\|(\ball{x}{r})}{r^m} \overset{\eqref{item:prop:upper-semicont-density:bordo}}&{\geq} {\lim\sup}_{i\to+\infty}e^{C_1r^{\alpha_1}}\frac{\|T\|(\ball{y_i}{r})}{r^m} \geq   {\lim\sup}_{i\to+\infty}\frac{\|T\|(\ball{y_i}{\rho|x_i-y_i|})}{(\rho |x_i-y_i|)^m}  \\
&= {\lim\sup}_{i\to+\infty}\frac{\|T_i^\prime\|(\ball{0}{\rho})}{\rho^m} \geq \frac{\|T^\prime\|(\ball{0}{\rho})}{\rho^m}.
\end{aligned}
$$
According to \cite[Reflection Principle, 3.2]{AllB}, we can reflect $T^\prime$ w.r.t. $T_x\Gamma$ to obtain a stationary varifold $V^{\prime\prime}$ in $\R^{m+n}$. We let $r\to 0$ and fix $\rho>0$ in the latter equation, thus we obtain

$$
\begin{aligned}
\Theta^m(T,x) &\geq   \frac{\|T^\prime\|(\ball{0}{\rho})}{\rho^m} 
\overset{\text{def of $V^{\prime\prime}$}}{=} \frac{1}{2}  \frac{\|V^{\prime\prime}\|(\ball{0}{\rho})}{\rho^m} = \frac{1}{2}  \frac{\|V^{\prime\prime}\|(\ball{x^\prime}{\rho})}{\rho^m} 
\overset{\eqref{item:prop:upper-semicont-density:interior}}{\geq} \frac{1}{2}\Theta^m(V^{\prime\prime},x^\prime) \\
\overset{\text{def of $V^{\prime\prime}$}}&{=} \Theta^m(T^{\prime},x^\prime) \overset{\eqref{item:prop:upper-semicont-density:interior}}{\geq} {\lim\sup}_{i\to\infty}\Theta^m(T^{\prime}_i,x^\prime_i) =  {\lim\sup}_{i\to\infty}\Theta^m(T,x_i),
\end{aligned}
$$
where in the last equality we do a standard rescaling argument.
\end{proof}
 
Let us state the existence of area minimizing tangent cones at boundary points of an almost area minimizing current. Although this particular result may be considered standard, we prove it here for the sake of completeness.

\begin{proposition}[Existence of area minimizing tangent cones]\label{prop:existence-minimizing-cones}
Let $T$ and $\Gamma$ be as in Assumptions \ref{assumptions}. If $T$ is $(C_0, r_0, \alpha_0)$-almost area minimizing in a neighbourhood $U$ of $x\in\Gamma$, then, for any sequence $r_k\to 0,$ there exists a blowup $\lim_{k\to+\infty}T_{x, r_k} = T_0\in\emph{\textbf{I}}^{loc}_m(\euc{m+n})$ such that:

\begin{enumerate}[\upshape (i)]
\item\label{item-convergence} $\| T_{x, r_k} \| \to \|T_0\|$ as $k\to +\infty$, in the sense of measures,

\item\label{item-area-minimizing} $T_0$ is area minimizing,

\item\label{item-density} $\| T_0 \| (\ball{0}{r}) = \Theta^m(T,x)\omega_{m}r^m, \forall r>0$,

\item\label{item-tangent-cone} $T_0$ is a tangent cone to $T$ at $x$.
\end{enumerate}

\end{proposition}
\begin{remark}
We only need to assume that $\partial T={Q^\star}\a{\Gamma}$ and $\Gamma$ is an $(m-1)$-submanifold of class $C^{1, \alpha}$, $\alpha\in(0,1)$, in order to deduce that $\partial T_0 = {Q^\star}\a{T_x\Gamma}$.
\end{remark}
\begin{proof}
Assume that $x=0$. By the almost monotonicity formula at the boundary and at the interior $($i.e., Proposition \ref{prop:almost-monotonicity} and \cite[Proposition 2.1]{DSS1}$)$, we have that 
$$ {\lim\sup}_{k\to+\infty}\|T_{r_k}\|(\ball{y}{1}) < +\infty, $$
for every $y\in U.$ Thus
$$ {\lim\sup}_{k\to+\infty}\|T_{r_k}\|(K) < +\infty, $$
for all compact set $K\subset\euc{m+n}$. Since the boundary of $T_{r_k}$ is ${Q^\star}\a{\Gamma/r_k}$ and $\Gamma$ is the graph of $u\in C^{1,\alpha}, u(0) = 0, Du(0) = 0$, we have 
\begin{equation*}
\begin{aligned}
\| \partial T_{r_k} \|(K) &= \| {\iota_{0, r_k}}_\sharp ({Q^\star}\a{\Gamma}) \|(K)  =
 \frac{1}{r_k^{m-1}}\|{Q^\star}\a{\Gamma}\|(r_k K) = \frac{{Q^\star}}{r_k^{m-1}}\mathcal{H}^{m-1}((r_k K)\cap \Gamma)\\
 & \leq  \frac{{Q^\star}}{r_k^{m-1}}\int_{\operatorname{proj}(r_k K)}\sqrt{1+|Du(z)|^2}dz \leq \frac{{Q^\star}}{r_k^{m-1}}\mathcal{H}^{m-1}(\operatorname{proj}(r_k K))\sqrt{1+C_{\Gamma}(\operatorname{diam}(K)r_k)^{2\alpha}}\\
 & \leq  {Q^\star}\omega_{m-1}\operatorname{diam}(K)^{m-1}\sqrt{1+C_\Gamma \operatorname{diam}(K)^{2\alpha}} =: C(\Gamma, K, \alpha, m, {Q^\star}),
\end{aligned}
\end{equation*}
and thus we can bound uniformly the mass of the boundary of $T_{r_k}$ in $K$. Therefore, we can use standard compactness results (one could consult \cite[Section 4.2]{Fed}) to ensure the existence of $T_0\in\textbf{I}^{loc}_m(\euc{m+n})$ such that $T_{r_k} \to T_0$, up to a subsequence, in the flat norm.

\textbf{Proof of $($\ref{item-convergence}$)$:} Let us write $T_{r_k} - T_0 = R_{r_k} + \partial \tilde{T}_{r_k}$ in $\ball{0}{R+2}$ with 
$${\lim\sup}_{k\to+\infty}\biggl(\|R_{r_k}\| (\ball{0}{R+1}) + \|\tilde{T}_{r_k}\| (\ball{0}{R+1})\biggr) = 0.$$
Thus, since all the measures involved are Radon measures, for almost every $s \in (R, R+1)$, it follows that
\begin{equation}\label{zero-mass-at-infinity}
{\lim\sup}_{k\to+\infty}\| R_{r_k}\| (\ball{0}{s}) = 0
\end{equation}  
and 
\begin{equation}\label{zeromassofslice}
{\lim\sup}_{k\to+\infty}\massnorm{\slice{\tilde{T}_{r_k}}{d}{s}} = 0.
\end{equation}
Note that (\ref{zeromassofslice}) follows directly from the formula of the slice and the fact that $T_{r_k}$ converges to $T_0$ in the flat norm. We may use again the slice formula to get
\begin{equation}\label{decomposition-of-blowup}
T_{r_k}\res {\ball{0}{s}} = T_0\res {\ball{0}{s}} + R_{r_k}\res {\ball{0}{s}} - \slice{\tilde{T}_{r_k}}{d}{s} + \partial (\tilde{T}_{r_k}\res {\ball{0}{s}}).
\end{equation}
The almost minimality condition gives
$$ \| T_{r_k} \|(\ball{0}{s}) \leq \left(1+C_0 (rs)^{\alpha_0}\right)\|T_{r_k}+\partial \tilde{T}_{r_k}\|\left(\ball{0}{s}\right).$$
Putting into account the latter inequality, the triangle inequality and (\ref{decomposition-of-blowup}), we obtain that
$$ \| T_{r_k} \|(\ball{0}{s}) \leq \left(1+C_0 (rs)^{\alpha_0}\right)\biggl( \| T_0\|(\ball{0}{s}) + \|R_{r_k}\|(\ball{0}{s}) + \massnorm{\slice{\tilde{T}_{r_k}}{d}{s}} + 2\| \partial \tilde{T}_{r_k} \|(\ball{0}{s}) \biggr).$$ 
Note that, by our construction, it follows that $\| \partial \tilde{T}_{r_k} \|(\ball{0}{s}) \to 0$ as $k\to+\infty$. Finally, by the lower semicontinuity of the mass, (\ref{zero-mass-at-infinity}), (\ref{zeromassofslice}) and the last equation passed through $\lim_{k\to+\infty}$, we conclude the proof of $($\ref{item-convergence}$)$. 

\textbf{Proof of $($\ref{item-area-minimizing}$)$:}
Fix $R\in (0, +\infty)$, by the lower semicontinuity of the mass, for all $\tilde{T}\in\mathbf{I}_{m+1}(\ball{0}{R})$ we have that
$$\begin{aligned}
\|T_0\|(\ball{0}{R}) & \leq {\lim\inf}_{k\to+\infty}\|T_{r_k}\|(\ball{0}{R})\\
& \leq  {\lim\inf}_{k\to+\infty}\left(1+C_0 ({r_k}R)^{\alpha_0}\right)\|T_{r_k}+\partial \tilde{T}\|\left(\ball{0}{R}\right) \\
\overset{\eqref{item-convergence}}&{=} \|T_0+\partial \tilde{T}\|\left(\ball{0}{R}\right),
\end{aligned}$$
for all $R>0$.

\textbf{Proof of $($\ref{item-density}$)$:} 
From $($\ref{item-convergence}$)$ and the almost monotonicity formula, we know that $\Theta^m(T,0)$ exists and, $\forall r > 0,$ we have that 
$$\begin{aligned}
\| T_0 \| (\ball{0}{r})  & =    \lim_{k\to+\infty}\| T_{r_k} \| (\ball{0}{r}) = \lim_{k\to+\infty}\frac{\| T \| (\ball{0}{r_kr})}{r_k^m} = \lim_{k\to+\infty}\frac{\omega_mr^m\| T \| (\ball{0}{r_kr})}{\omega_m (r_k r)^m}\\
& = \Theta^m(T,0)\omega_{m}r^m.
\end{aligned}$$

\textbf{Proof of $($\ref{item-tangent-cone}$)$:} 
By following closely the argument given in \cite[Theorem 3.1, Chapter 7]{simon2014introduction} using \eqref{item-density}, we prove that $T_0$ is in fact a cone. Indeed, by \cite[Theorem 3.2]{DNS} applied for $T_0$ which is an area minimizing cone, we know that 
\begin{eqnarray*}
{Q^\star}\int_{T_0\Gamma\cap\ball{0}{\rho}}(x-p) \cdot \vec{n}(x) d \mathcal{H}^{m-1}(x) = 0,
\end{eqnarray*}
since $\vec{n} \perp T_0{\Gamma}$. Using \eqref{item-area-minimizing}, we obtain that $\vec{H}_{T}$ is zero a.e., so, by \eqref{item-density}, we obtain the constancy of the mass ratio and, then using again \cite[Theorem 3.2]{DNS}, we get that
\begin{eqnarray*}
\int_{\ball{0}{r} \backslash \ball{0}{s}} \frac{\left|x^{\perp}\right|^{2}}{|x|^{m+2}} d\|T\|(x) = 0.
\end{eqnarray*}
Then, if we fix a cone $C$ such that $\partial C = \partial T_0$, notice that since $\partial T_0 = {Q^\star}\a{T_0\Gamma}$ we may choose for instance $C$ as a half subspace, we can apply \cite[Lemma 2.33, Chapter 6]{simon2014introduction} for $T_0 -C$ and conclude that it is a cone and thus $T_0 = T_0 - C + C$ is a cone.
\end{proof}

\section{Stratification for \texorpdfstring{$(C_0, r_0, \alpha_0)$}{Lg}-almost area minimizing currents}\label{s:stratification}

\subsection{Stratification theorem}

\begin{definition}[Conelike functions]\label{def:conelike-functions}
An upper semicontinuous function $g:\R^{m+n}\to\R_{+}$ is called \textbf{conelike} provided:
\begin{enumerate}[\upshape (i)]
    \item $g(\lambda x) = g(x)$ for all $\lambda >0$ and for all $x\in\mathbb{R}^{m+n}$,
    
    \item If $g(x) = g(0)$, then $g(x+\lambda v) = g(x+v)$ for all $\lambda > 0$ and $v\in\R^{m+n}$.
\end{enumerate}
If $g$ is conelike we also define the \textbf{spine of $g$} as the set 
$$ \operatorname{spine}(g) := \{x\in\mathbb{R}^{m+n}: g(x) = g(0)\}.$$
\end{definition}

By (i) in the last definition and upper semicontinuity, we have that $g(z)\leq g(0)$ for all $z$. Note that, by \cite[Theorem 3.1]{White97}, $\operatorname{spine}(g)$ is a vector subspace and 
$$ \operatorname{spine}(g) = \{x\in\mathbb{R}^{m+n}: g(x+v) = g(v), \ \forall v\in\mathbb{R}^{m+n}\}.$$
Fix $T$ and $\Gamma$ as in Assumption \ref{assumptions}, and set the class of functions $$\mathscr{G}(p):= \{g_{p, T_0}: {T_0} \ \text{is a tangent cone to} \  T \ \text{at} \ p\},$$

where 

$$f_T(x) := \left\{\begin{matrix}\begin{aligned}
&\Theta^m(T,x), &x\notin\Gamma \\ 
&2\Theta^m(T,x) + 1, &x\in\Gamma
\end{aligned}\end{matrix}\right. ,\ \ \text{and} \ \ g_{p, T_0}(x) := \left\{\begin{matrix}\begin{aligned}
&\Theta^m({T_0},x), &x\notin T_p\Gamma \\ 
&2\Theta^m({T_0},x) + 1, &x\in T_p\Gamma
\end{aligned}\end{matrix}\right. ,$$
then $f$ and each $g_{p, T_0}$ are upper semicontinuous from Proprosition \ref{prop:upper-semicont-density}.

\begin{definition}[Spine of a cone]\label{defi:spine}
Let $p\in\Gamma$, $T_0$ be an oriented tangent cone with $\partial T_0 = {Q^\star}\a{T_p\Gamma}$, we define the {\textbf{spine of $T_0$}},  and denote it by $\spine(T_0)$, to be the set of vectors $v\in T_p\Gamma$ such that $(\tau_v)_\sharp T_0 = T_0$ where $\tau_v(w) = w + v$. Clearly the $\operatorname{spine}(T_0)$ is always a subspace of $T_p\Gamma$. 
\end{definition}

We now provide an equivalence for the definition of spine of an oriented tangent cone which follows the ideas furnished by Almgren in his stratification process, \cite[Theorem 2.26]{Alm}, and used in \cite{White97}, \cite{LeonSimon}, \cite{DDHM}.

\begin{lemma}[Spine as constant density set]\label{lemma:equiv-def-of-spine}
Let $T$ and $\Gamma$ be as in Assumption \ref{assumptions} and $T_0$ be an oriented tangent cone to $T$ at $p\in\Gamma$. Then $\operatorname{spine}(T_0) = \{x\in T_p\Gamma :  \Theta^m(T_0, x)= \Theta^m(T_0, 0)\}.$
\end{lemma}
\begin{proof}
Take $x\in\operatorname{spine}(T_0)$, by the definition of spine, we have the third equality below

$$\begin{aligned}
\Theta^m(T_0, x) = \lim_{r\to  0}\frac{\|T_0\|(\ball{x}{r})}{\omega_mr^m} = \lim_{r\to  0}\frac{\|(\tau_x)_\sharp T_0\|(\ball{0}{r})}{\omega_mr^m} = \lim_{r\to  0}\frac{\|T_0\|(\ball{0}{r})}{\omega_mr^m} = \Theta^m(T_0, 0) .
\end{aligned} 
$$

On the other hand, consider $x\in T_p\Gamma$ such that $\Theta^m(T_0,x) = \Theta^m(T_0,0)$, we claim that $x\in\operatorname{spine}(T_0)$. To prove this claim, we apply the monotonicity formula \cite[Thm 3.2]{DDHM} (which works mutatis mutandis for the higher multiplicity case) to the cone $T_0$, which is area minimizing and $\partial T_0 = Q^\star\a{T_p\Gamma}$, to obtain, for $0<s<r,$ 
\begin{equation*}
\begin{aligned}
    \frac{\| T_0 \|(\ball{x}{r})}{ r^m} - \frac{\| T_0 \|(\ball{x}{s})}{ s^m} &= \int_{\ball{x}{r}\setminus\ball{x}{s}}\frac{|(z-x)^\perp |^2}{|z-x|^{m+2}}\operatorname{d}\|T_0\|(z) \\
    &\quad + \int_s^r\rho^{m-1}\left( \int_{\ball{x}{\rho}}(z-x)^{\perp}\cdot H_{T_0}(z)\mathrm{d}\|T_0\|(z) \right. \\
    &\quad\quad\quad\quad\quad\quad \left.+  Q^{\star}\int_{T_p\Gamma\cap\ball{x}{\rho}}(z-x)\cdot \vec{n}_{T_0}(z)\mathrm{d}\mathcal{H}^{m-1}(z)\right)\mathrm{d}\rho.
\end{aligned}
\end{equation*}  
Note that, since $T_0$ is a cone and area miniminzing, we get $H_{T_0}$ vanishes $\|T_0\|$-a.e. and $(x-p)\cdot\vec{n}_{T_0}(x) =0$ for all $x\in T_p\Gamma\cap\ball{p}{\rho}$. Hence, by the last displayed equation, we obtain
\begin{equation}\label{eq:lemma-spine-constant-density-set}
    \frac{\| T_0 \|(\ball{x}{r})}{ r^m} - \frac{\| T_0 \|(\ball{x}{s})}{ s^m} = \int_{\ball{x}{r}\setminus\ball{x}{s}}\frac{|(z-x)^\perp |^2}{|z-x|^{m+2}}\operatorname{d}\|T_0\|(z).
\end{equation} 
Notice that the following holds true: $r_k^{m}r^{-m}\|T\| (\ball{x}{r_k^{-1}r}) = \|T\|(\ball{r_k x}{r})/r^m$ implies 

$$ \limsup_{r_k \to 0}\frac{\|T\|\left(\ball{x}{\frac{r}{r_k}}\right)}{\left(\frac{r}{r_k}\right)^m} \leq \frac{\|T\|(\ball{0}{r+\eta})}{(r+\eta)^m} = \Theta^m(T,0), $$

for any $\eta>0$. The last inequality and again using that $r \mapsto \|T\|(\ball{x}{r})/(\omega_m r^m)$ is nondecreasing, which follows from \eqref{eq:lemma-spine-constant-density-set}, ensure that

$$\Theta^m(T,x) \leq  \frac{\|T\|(\ball{x}{r})}{\omega_m r^m} \leq \limsup_{r_k \to 0}\frac{\|T\|\left(\ball{x}{\frac{r}{r_k}}\right)}{\left(\frac{r}{r_k}\right)^m} \leq\Theta^m(T,0) = \Theta^m(T,x),$$ for every $r>0$.

Hence, the mass ratio does not depend on $r$. Therefore, we have that $\|({\tau_{x}})_\sharp T_0\|(\ball{0}{r}) = \|T_0\|(\ball{x}{r}) = \| T_0\|(\ball{0}{r})$ thus, by measure theory, we get that $\| T_0 \| = \| ({\tau_x})_\sharp T_0 \|$ as measures, which in turn ensures that $T_0 = \pm ({\tau_x})_\sharp T_0$. We know that $\partial T_0 = Q^\star\a{T_p\Gamma}$ and $\partial ({\tau_x})_\sharp T_0 = ({\tau_x})_\sharp \partial T_0 =Q^\star\a{T_p\Gamma}$, then we conclude that $T_0 = ({\tau_x})_\sharp T_0$ which shows that $x\in\operatorname{spine}(T_0)$.
\end{proof} 

We would like to show that we are in position to apply \cite[Theorem 3.2]{White97}. To that end, we prove the following lemma.

\begin{lemma}\label{lemma:densities-conelike}
For each $p\in\Gamma$ and each oriented tangent cone $T_0$ to $T$ at $p$, $g_{p, T_0}$ is conelike.
\end{lemma}
\begin{proof}
Property (i) in Definition \ref{def:conelike-functions} is a direct consequence of the scaling invariance of ${T_0}$. To what concerns property (ii), if $g_{p, T_0}(x)=g_{p, T_0}(0),$ and $x\in T_p\Gamma,$ we have that $x\in \operatorname{spine}(T_0)$ and thus

$$ 
\begin{aligned}
\Theta^m(T_0, x+v) &= \Theta^m(({\tau_x})_\sharp T_0, v) = \Theta^m(T_0, v) = \Theta^m(T_0, \lambda v) = \Theta^m(({\tau_x})_\sharp T_0, \lambda v)\\
&= \Theta^m(T_0, x+ \lambda v),
\end{aligned}
$$

for any $\lambda > 0$ and $v\in\R^{m+n}$.

If $g_{p, T_0}(x)=g_{p, T_0}(0)$ and $x\notin T_p\Gamma$, then, by definition of $g_{p, T_0}$, we have $\Theta^m(T_0, x) = 2\Theta^m(T_0, 0) + 1$. Since $T_0$ is a cone we have that $\Theta^m(T_0, x)=\Theta^m(T_0, \lambda x)$, for every $\lambda>0$. Then $\Theta^m(T_0, \lambda x) = 2\Theta^m(T_0, 0) + 1$, for every $\lambda>0$. Taking the limsup and recalling Proposition \ref{prop:upper-semicont-density} (iii) we get $\limsup_{\lambda\to0^+}\Theta^m(T_0, \lambda x) = 2\Theta^m(T_0, 0) + 1\le2\Theta^m(T_0, 0)$, which is a contradiction. 
\end{proof}

\begin{remark}{
Note that, a simple consequence of the Lemma \ref{lemma:equiv-def-of-spine} and the proof of Lemma \ref{lemma:densities-conelike} is that  $\operatorname{spine}(g_{p, T_0}) = \operatorname{spine}({T_0}).$ }
\end{remark}

\begin{definition}[Stratum]\label{defi:stratum}
Let $p\in\Gamma$ and $T$ be a $m$-current with $\partial T = {Q^\star}\a{\Gamma}$, we define the {\textbf{$j$-stratum of $\Gamma$ with respect to $T$}} as the set
$$\mathscr{P}_j(T, \Gamma) = \{p\in\Gamma : \dim(\spine(T_0))\leq j, \ \text{for all tangent cone} \ T_0 \ \text{to} \ T \ \text{at} \ p\}  .$$
\end{definition}

Now, we shall directly apply Proposition \ref{prop:existence-minimizing-cones} to check the conditions (1) and (2) of \cite[Theorem 3.2]{White97} which in turn furnishes

\begin{theorem}[Stratification Theorem, Theorem 3.2 of \cite{White97}]\label{thm:stratification}
For $T$ and $\Gamma$ as in Assumption \ref{assumptions}, let
$$ \Sigma_i := \{x: f_T(x) >0 \ \text{and} \ \sup\{\operatorname{dim}(\operatorname{spine}(g)): g\in\mathscr{G}(x)\} \leq i\},$$
then the Hausdorff dimension of $\Sigma_i$ is at most $i$ and $\Sigma_0$ is at most countable. In particular, we have the same statements for the stratum $\mathscr{P}_i(T, \Gamma)$, i.e., the Hausdorff dimension of $\mathscr{P}_i(T, \Gamma)$ is at most $i$, $\mathscr{P}_0(T, \Gamma)$ is at most countable, and 
$$ \mathscr{P}_0(T, \Gamma)\subset\mathscr{P}_1(T, \Gamma)\subset\dots\subset\mathscr{P}_{m-1}(T, \Gamma) = \Gamma.$$
\end{theorem}

\subsection{Open books, flat cones, one-sided and two-sided points}

The characterization of tangent cones is an important tool in the subsequent theory, in fact, when dealing with two dimensional area minimizing cones we have general structure results, see for instance \cite[Lemma 3.1]{hirsch2019uniqueness} or \cite[Proposition 4.1]{DNS2}. If we consider arbitrary dimensions, there is no general structure theorem for area minimizing tangent cones, however, assuming that the density of the cone is constant along the boundary of the cone which is equivalent to assume that the spine has maximal dimension, as showed in Lemma \ref{lemma:equiv-def-of-spine}, we can characterize tangent cones as in \cite[Theorem 5.1]{brothers1977existence}, or \cite[Lemma 3.17]{DDHM}. We enunciate and prove in Lemma \ref{lemma:max-spine-dim} rigorous statement of the assertions just mentioned. To go further in our treatment we need the following definitions.

\begin{definition}[Open books]\label{def:open-book}
Let $T_0 \in \mathbf{I}_{m}^{\operatorname{loc}}\left(\euc{m+n}\right)$ be an oriented cone and $V$ is an oriented $(m-1)$-dimensional linear subspace of $\euc{m+n}$. We say that $T_0$ is an {\textbf{open book with boundary $\a{V}$ and multiplicity ${Q^\star}$}}, if $\partial T_0 = {Q^\star}\a{V}$ and there exist $N\in\N\setminus\{0\}$, ${Q}_1,\dots, {Q}_N\in\N\setminus\{0\}$ and $\pi_1,\dots, \pi_N$ distinct $m$-dimensional half-planes such that 
\begin{enumerate}[\upshape (i)]
\item $\partial\a{\pi_i} = \a{V}, \forall i\in\{1, \dots, N\}$,

\item $T_0 = \sum_{i=1}^{N}{Q}_i\a{\pi_i}$ with ${Q^\star} = \sum_{i=1}^{N}{Q}_i$.
\end{enumerate}
If there exist $i\neq j$ such that $\pi_i\neq\pi_j$, we say that $T_0$ is a {\textbf{genuine open book with boundary $\a{V}$ and multiplicity ${Q^\star}$}}.
\end{definition}

\begin{definition}[Flat cones]\label{def:flat-cone}
Let $T_0 \in \mathbf{I}_{m}^{\operatorname{loc}}\left(\euc{m+n}\right)$ be an oriented cone and $V$ is an oriented $(m-1)$-dimensional linear subspace of $\euc{m+n}$. We say that $T_0$ is a {\textbf{flat cone with boundary $\a{V}$ and multiplicity ${Q^\star}$}}, if $\partial T_0 = {Q^\star}\a{V}$ and there exist a $m$-dimensional closed plane $\pi$, $Q^{\textup{int}}\in\N$ and $Q,Q^{\star}\in\N\setminus\{0\}, Q\ge Q^\star$, such that 
\begin{enumerate}[\upshape (i)]
    \item $\spt T_0 = \pi$ is an $m$-dimensional subspace,

    \item $\partial\a{\pi^+} = -\partial\a{\pi^-} = \a{V}$,

    \item $T_0 = Q^{\textup{int}}\a{\pi} + Q\a{\pi^+} + (-1)^k(Q-{Q^\star})\a{\pi^-}$, $k\in\{0,1\}$.
\end{enumerate}
If either
\begin{enumerate}[\upshape (a)]
    \item $Q^{\textup{int}}=0$ and $Q=Q^\star$, or

    \item $Q^{\textup{int}}=0$ and $k=1$,
\end{enumerate}
we call $T_0$ an {\textbf{one-sided boundary flat cone with multiplicity $Q^\star$}}. If $k=0$ and $T_0$ is not an one-sided boundary flat cone, we say that $T_0$ is a {\textbf{two-sided boundary flat cone with multiplicity $Q^\star$}}.
\end{definition}

Note that: if $T_0$ is an open book which is not genuine, then $T_0$ is an one-sided boundary flat cone.

\begin{definition}
Let $T_0 \in \mathbf{I}_{m}^{\operatorname{loc}}\left(\euc{m+n}\right)$ be an oriented cone and $V$ is an oriented $(m-1)$-dimensional linear subspace of $\euc{m+n}$ such that $\partial T_0 = Q^\star\a{V}$. If $p\in\spt(\partial T_0)$, we say that
\begin{enumerate}[\upshape (i)]
    \item $p$ is a {\textbf{boundary flat point}} provided $T_0$ is a flat cone,
    
    \item $p$ is a {\textbf{one-sided boundary flat point}} provided $T_0$ is an open book which is non genuine,
    
    \item $p$ is a {\textbf{two-sided boundary flat point}} provided $T_0$ is a two-sided boundary flat cone.
\end{enumerate}
\end{definition}

\begin{lemma}[The set of one-sided points is open]\label{lemma:neighbourhood-density-bounded-above}
Let $T, \Gamma$ and $p\in\Gamma$ as in Assumption \ref{assumptions}. If $\Theta^m(T,p) < \frac{Q^\star + 1}{2}$, then there exists a neighbourhood $U$ of $p$ such that $\Theta^m(T,q)<\frac{Q^\star + 1}{2}$ for every $q\in U\cap\Gamma$.
\end{lemma}
\begin{proof}
It follows directly from the upper semicontinuity of the density function, see Proposition \ref{prop:upper-semicont-density}.
\end{proof}

Note that, if $m=2$ and the tangent cone $T_p$ is a two-sided boundary flat cone, Theorem \ref{thm:uniqueness-tang-cones-decay} ensures that the least possible density on $p$ is $\frac{Q^\star}{2} +1$. The following lemma is a generalization of \cite[Lemma 3.17]{DDHM} to the case of higher multiplicity with essentially the same proof.

\begin{lemma}\label{lemma:max-spine-dim}
Let $T, \Gamma$ and $p\in\Gamma$ as in Assumption \ref{assumptions} with $C_0=0$. If $T_0$ is an oriented tangent cone to $T$ at $p$ with $\dim(\spine(T_0)) = m-1$, then
\begin{enumerate}[\upshape (i)]
    \item If $\Theta^m(T_0,0) = \frac{Q^\star}{2}$, $T_0$ is an open book,
    
    \item If $\Theta^m(T_0,0) > \frac{Q^\star}{2}$, $T_0$ is a two-sided boundary flat cone.
\end{enumerate}
\end{lemma}
\begin{remark}
{We also mention that the assumption that the $\spine$ has maximal dimension in $\Gamma$ was assumed in a similar fashion (see Lemma \ref{eq:lemma-spine-constant-density-set}) in \cite[Theorem 5.1]{brothers1977existence}.}
\end{remark}
\begin{proof}
By the assumption, we have that $\spine(T_0) = T_p\Gamma$. By \cite[Theorem 2.2 (3)]{Alm}, there exists an one-dimensional area minimizing current ${T_0}_1$ in $(T_p\Gamma)^\perp$ such that $T_0 = \a{T_p\Gamma}\times {T_0}_1$. This fact that Almgren proved is an application of \cite[Theorem 5.4.8]{Fed} and \cite[Section 4.3.15]{Fed}, using that $T_0$ is an oriented cylinder with direction $v$ for any $v\in T_p\Gamma = \operatorname{spine}(T_0)$. Thus, since $\partial T_0 = {Q^\star}\a{T_p\Gamma}$, we obtain that $\partial {T_0}_1 = (-1)^{m-1}{Q^\star}\a{0}$, the fact that ${T_0}_1$ is invariant under homotheties allows us to write, for some $Q\in\N\setminus\{0\}$,
$$ (-1)^{m-1}{T_0}_1 = \sum_{i=1}^{{Q}}\a{\ell_i^+} + \sum_{j=1}^{Q-Q^\star}\a{\ell_i^-}, \ \|{T_0}_1\| = \sum_{i=1}^{{Q}}\|\a{\ell_i^+}\| + \sum_{j=1}^{Q-Q^\star}\|\a{\ell_i^-}\|,$$
where $\ell_i^+$, $\ell_j^-$ are all oriented half-lines such that $\partial\a{\ell_i^+}=\a{0}$ and $\partial\a{\ell_j^-}=-\a{0}$. In particular, we have 
\begin{equation}\label{eq:bdim-maximal-implie-flatness}
\spt(T_0) \subset \biggl(\bigcup_{i=1}^{Q} T_p\Gamma + \ell_i^+\biggr)\cup\biggl(\bigcup_{i=1}^{Q-Q^\star} T_p\Gamma + \ell_i^-\biggr).
\end{equation}

Note that, if $Q>Q^\star$, $\partial (\a{\ell_i^+}+\a{\ell_j^-}) = 0$ and $\a{\ell_i^+}+\a{\ell_j^-}$ is area minimizing for any choice of $i$ and $j$ which ensures that the support of $\a{\ell_i^+}+\a{\ell_j^-}$ is a straight line $\ell_{ij}$. Since the choice of $i$ and  $j$ is arbitrary, then we have $\spt(\a{\ell_i^+}+\a{\ell_j^-}) \subset \ell$, where $\ell$ is a straight line which, by \eqref{eq:bdim-maximal-implie-flatness}, concludes the proof of (ii). If $Q = Q^\star$, we have that ${T_0}_1$ is a sum of lines which might be distinct, and this concludes the proof of (i).
\end{proof}

\subsection{Two dimensional case}\label{two-dim-case}

In Lemma \ref{lemma:max-spine-dim}, $T_0$ has dimension $m$ and we assume that the dimension of the spine is maximal. Nevertheless, if $m=2$, we can drop the hypothesis on the spine since we have a full characterization of tangent cones with boundary being a subspace as stated in the proposition below.

\begin{proposition}[Proposition 4.1, \cite{DNS2}]\label{prop:2d-tangent-cone}
Let ${T_0}$ be a $2$-dimensional area minimizing cone in $\mathbb{R}^{2+n}$ with $\partial {T_0}=Q^\star \a{\ell}$ for some positive integer $Q^\star$ and a straight line $\ell$ containing the origin. Then we can decompose ${T_0}={T_0}^{\emph{int}}+{T_0}^{\flat}$ into two area minimizing cones with supports intersecting only at the origin which satisfy
\begin{enumerate}[\upshape (i)]
    \item $\partial {T_0}^{\emph{int}}=0$ and thus ${T_0}^{\emph{int }}=\sum_{i=1}^{N} Q_{i} \a{ \pi_{i} }$ where $Q_{1}, \ldots, Q_{N}$ are positive integers and $\pi_{1}, \ldots, \pi_{N}$ are distinct oriented $2$-dimensional planes such that $\pi_{i} \cap \pi_{j}=\{0\}$ for all $i \neq j$,
    
    \item ${T_0}^{\flat}$ is either a two-sided boundary flat cone or an open book.
\end{enumerate}
\end{proposition}

Let us also recall two pivotal results in the theory which will be used in this work. We also denote $\dist_{H}$ for the Hausdorff distance between closed sets and we denote by $e(p, r)$ the \textbf{spherical excess of a current $T$}, namely
$$
e(p, r):=\frac{\|T\|\left(\ball{p}{r}\right)}{\pi r^{2}}-\Theta(T, p),
$$

we also define the H\"older seminorm used to measure the regularity of $\Gamma$, for any open set $U$,

$$
[\Gamma]_{0, \alpha, U}:=\sup _{q \neq p \in \Gamma \cap U} \frac{\left|T_{p} \Gamma-T_{q} \Gamma\right|}{|p-q|^{\alpha}} .
$$

We have the following decay properties.

\begin{theorem}[Uniqueness of tangent cones and speed of convergence, Theorem 2.1, \cite{DNS2}]\label{thm:uniqueness-tang-cones-decay}
Let $T$ and $\Gamma$ be as in Assumption \ref{assumptions} with $m=2$. Then there are positive constants $\varepsilon_{0}$, $C$ and $\beta$ with the following property. If $p \in \Gamma$ and $e(p, r) \leq \varepsilon_{0}^{2}$ for some $r \leq \operatorname{dist}\left(p, \partial \ball{0}{1}\right)$, then there exists a unique tangent cone $T_p$ to $T$ at $p$ which, for every $\rho\in (0, r]$, satisfies:
$$ \begin{aligned}
|e(p, \rho)| &\leq C|e(p, r)|\left(\frac{\rho}{r}\right)^{2 \beta}+C\left(C_0^{2}+[\Gamma]_{0, \beta, \ball{p}{r}}^{2}\right)\left(\frac{\rho}{r}\right)^{2 \beta}, \\
\fdistin{\ball{0}{1}}{T_{p, \rho}}{T_p} &\leq C|e(p, r)|^{\frac{1}{2}}\left(\frac{\rho}{r}\right)^{\beta}+C\left(C_0+[\Gamma]_{0, \beta, \ball{p}{r}}\right)\left(\frac{\rho}{r}\right)^{\beta},\\
\operatorname{dist}_{H}\left(\operatorname{spt}\left(T_{p, \rho}\right) \cap \ball{0}{1}, \spt(T_p) \cap \ball{0}{1}\right) &\leq C|e(p, r)|^{\frac{1}{2}}\left(\frac{\rho}{r}\right)^{\beta}+C\left(C_0+[\Gamma]_{0, \beta, \ball{p}{r}}\right)\left(\frac{\rho}{r}\right)^{\beta}.
\end{aligned} $$
\end{theorem}

We also state the H\"older continuity of the map that to each point $p\in\Gamma$ assigns its unique tangent cone $T_p$.

\begin{lemma}[H\"older continuity]\label{lemma:Holder-continuity-tang-cones}
Let $T,p, r$ be as in Theorem \ref{thm:uniqueness-tang-cones-decay} and $q\in\Gamma\cap\ball{p}{r}$. Then the functions $q\mapsto T_q$ is H\"older continuous, i.e., it holds
\begin{equation}
\fdist{T_q\res\ball{0}{1}}{T_p\res\ball{0}{1}} \leq C|q-p|^{\beta}, \quad \forall q \in \ball{p}{r}.
\end{equation}
\end{lemma}
\begin{proof}
For the proof we refer the reader directly to \cite[Equation 4.7]{DNS} which can be readily adjusted to the almost area minimizing setting.
\end{proof}

\begin{definition}[Two-sided collapsed points]\label{defi:collapsed}
Let $T$ and $\Gamma$ be as in Assumption \ref{assumptions}. A point $p\in\Gamma$ will be called {\textbf{two-sided collapsed point of $T$}} if
\begin{enumerate}[\upshape (i)]
\item\label{def:collapsed-point-tangent-cone} there exists a tangent cone $T_0$ to $T$ at $p$ which is a two-sided boundary flat cone,

\item\label{def:collapsed-point-density} there exists a neighbourhood $U$ of $p$ such that $\Theta(T, q)\geq\Theta(T, p)$ for every $q\in\Gamma\cap U$.
\end{enumerate}
\end{definition}

\begin{figure}[H]
\centering
\captionsetup{width=.7\linewidth}\includegraphics[width=.5\linewidth]{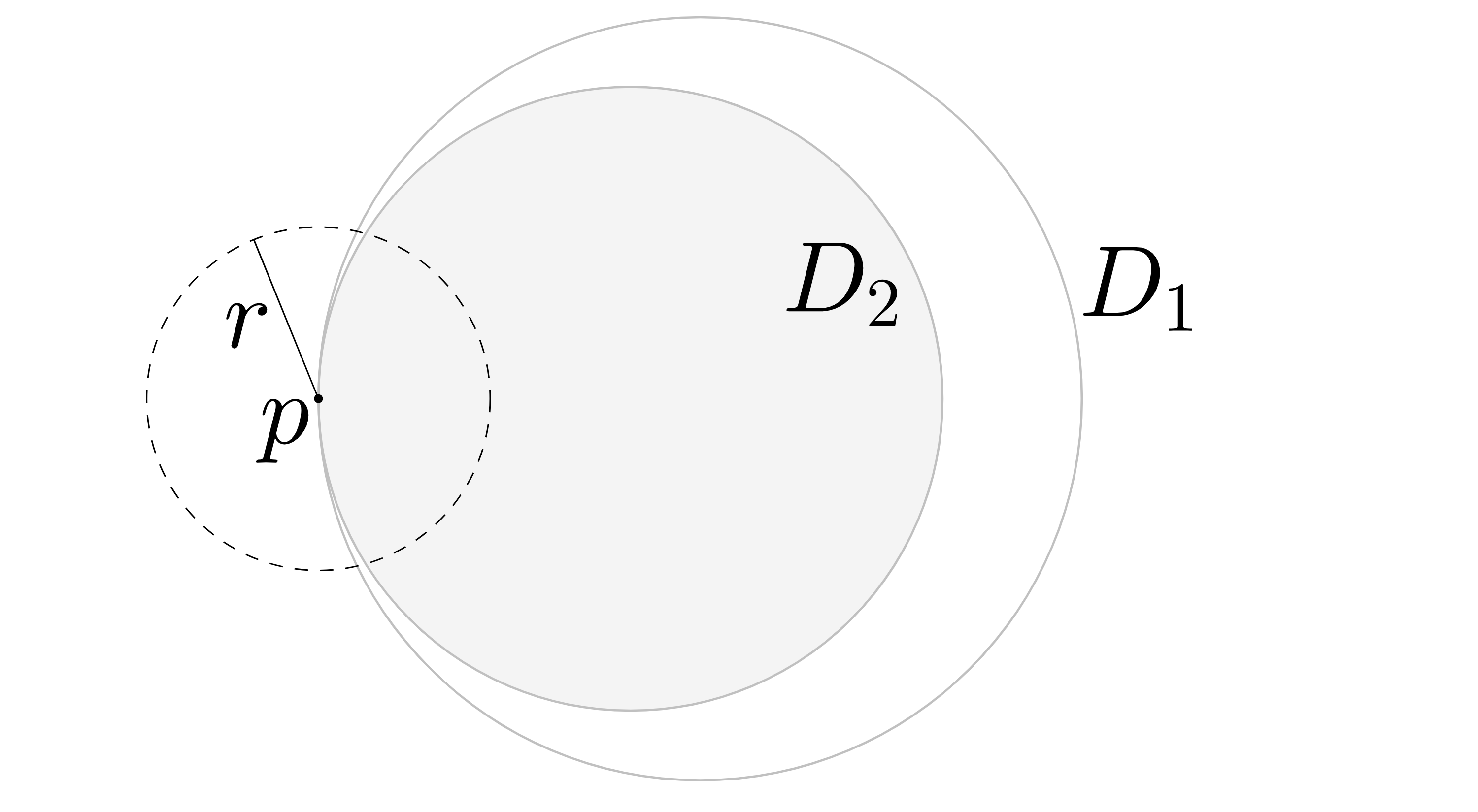}
\caption{Condition \eqref{def:collapsed-point-density} essentially excludes points of tangential intersection of connected parts of $\spt T$, i.e., it forbid the existence of one-sided points arbitrarily close to two-sided points. For instance as in the picture $p$ does not verify \eqref{def:collapsed-point-density}, let $T = \a{D_1}+\a{D_2}$ with $D_1$ and $D_2$ being tangential circles at $p$ then $\Theta(T, p) = 1 \geq\frac{1}{2} = \Theta(T, q)$ for all $q\neq p$ which belongs to the outer circumference.}
\label{}
\end{figure}

\begin{lemma}[The set of two-sided collapsed points is open]\label{lemma:collapsed-set-is-open}
Let $T$ and $p$ as in Theorem \ref{thm:uniqueness-tang-cones-decay}. Assume that $p \in \Gamma$ is a two-sided collapsed point, then there is $\rho >0$ such that $\Theta^2(T,q) = \Theta^2(T,p)$ for all $q\in \ball{p}{\rho}\cap \Gamma$. In particular, every such $q$ is two-sided collapsed.
\end{lemma}
\begin{proof}
Fix $Q\in\N, Q>Q^\star$ such that $\Theta^2(T,p) = Q - \frac{Q^\star}{2}$ and the unique tangent cone to $T$ at $p$ is $T_p = Q\a{\pi^+}+(Q-Q^\star)\a{\pi^-}$. If we choose $r>0$ small enough we can assume that 
$$ \frac{\|T\|(\ball{p}{r})}{\pi r^2} \leq Q- \frac{Q^\star}{2} + \frac{1}{8}.$$
Now, we choose $s\in (0, r)$, in order to hold
$$ \|T\|(\ball{q}{r-s}) \leq \|T\|(\ball{p}{r}) \leq \pi r^2\biggl( Q - \frac{Q^\star}{2} + \frac{1}{8}\biggr) \leq \pi (r-s)^2\biggl( Q - \frac{Q^\star}{2}+\frac{3}{16}\biggr) ,$$
for every $q\in\ball{p}{s}\cap\Gamma$. For every $\sigma \in (0, r-s)$, by the almost monotonicity formula, Proposition \ref{prop:almost-monotonicity}, we have that 
\begin{equation}\label{eq:lemma-collapsed-are-interior:bound-mass-ratio}
    \begin{aligned}
\frac{\|T\|(\ball{q}{\sigma})}{\pi\sigma^2} &\leq e^{C_1((r-s)^{\beta_1} - \sigma^{\beta_1})}\frac{\|T\|(\ball{q}{r-s})}{\pi (r-s)^2} \\
&\leq e^{C_1(r-s)^{\beta_1}}\biggl( Q - \frac{Q^\star}{2} + \frac{3}{16} \biggr), 
\end{aligned}
\end{equation}
for every $q\in\ball{p}{s}\cap\Gamma$. Then take $q\in\ball{p}{s}\cap\Gamma$ where this ball is chosen to be a subset of the neighbourhood $U$ given by the definition of two-sided collapsed points. Hence, by Proposition \ref{prop:2d-tangent-cone} and Lemma \ref{lemma:max-spine-dim}, the tangent cone to $T$ at $q$ has to be of the form $T_q = Q^\prime\a{\pi_q^+}+(Q^\prime-Q^\star)\a{\pi_q^-}$, for some integer $Q^\prime>Q^\star$, thus, letting $r\to 0$ in \eqref{eq:lemma-collapsed-are-interior:bound-mass-ratio} we obtain

$$ Q-\frac{Q^\star}{2} = \Theta^2(T,p) \leq \Theta^2(T,q) = Q^\prime - \frac{Q^\star}{2} \overset{\eqref{eq:lemma-collapsed-are-interior:bound-mass-ratio}}{\leq}Q - \frac{Q^\star}{2} + \frac{3}{16}.$$
\end{proof}

The following theorem allows us to reduce the proof of Theorem \ref{thm:main-theorem} to the proof that any two-sided collapsed point is regular.

\begin{theorem}\label{theorem:not-dense-collapsed-singular}
Let $T$ and $p$ as in Theorem \ref{thm:uniqueness-tang-cones-decay} and assume that $C_0=0$. If $\breg(T)$ is not dense in $\Gamma$, then there exists a two-sided collapsed singular point $p\in\Gamma$ with $\Theta^2(T,p)>\frac{Q^\star}{2}$.
\end{theorem}
\begin{remark}
When we consider the setting of \cite{DNS}, i.e., when $\Gamma$ belongs to a $C^{3,\alpha}$ convex barrier, $\alpha\in (0,1)$, we have that two-sided points do not exist. In particular, $\breg^2(T)=\emptyset$ and, by \cite[Theorem 0.2]{DNS}, we know that $\breg^1(T) = \Gamma$. In other words, the authors in \cite{DNS} proved the full regularity of the current at the boundary.
\end{remark}
\begin{proof}
Assume that $\bsing(T)$ has no empty interior, then we can define

$$ C_i := \biggl \{p\in\Gamma: \Theta^2(T, p)\geq i - \frac{1}{2}\biggr \}\cap \interior(\bsing(T)).$$

By Proposition \ref{prop:upper-semicont-density}, the density restricted to the boundary is upper semi-continuous, then $C_i$ is relatively closed in $\interior(\bsing(T))$. Let $D_i$ be the topological interior of $C_i$ and $E_i$ be the relatively open set $D_i\setminus C_{i+1}$ in $\interior(\bsing(T))$. We fix $p\in\Gamma$ and the natural number $i$ such that 

\begin{equation}\label{eq:stratification:density-estimate}
i - \frac12 \leq \Theta^2(T,p) < i +\frac12.
\end{equation}

Assume that $p\notin \cup_{i\geq 1}E_i$, by the latter inequalities, we have $p\in C_i\setminus D_i$ which leads to 

$$
\interior(\bsing(T))\setminus \cup_{i}E_i\subset \cup_{i}C_i\setminus D_i.
$$

Observe that $C_i\setminus D_i$ is relatively closed in $\interior(\bsing(T))$ and then $\interior(\bsing(T))\setminus \cup_{i}E_i$ is the union of countably many closed subsets of $\interior(\bsing(T))$ which guarantees, by the Baire Category Theorem, that $\cup_i E_i$ cannot be empty. So, there is $E_i\neq\emptyset$ relatively open in $\Gamma$, hence, in view of \cite[Theorem 0.2]{DNS}, since $E_i$ contains only singular points, any $p\in E_i$ satisfies $\Theta^2(T,p)\geq \frac{Q^\star +1}{2}$ and, from Proposition \ref{prop:2d-tangent-cone} and Lemma \ref{lemma:max-spine-dim}, $p$ is two-sided boundary flat point. Fix $p\in E_i$, we know that there exists $Q\in\N, Q>Q^*$ such that $\Theta^2(T,p) = Q -\frac{Q^\star}{2}$, thus
\begin{itemize}
    \item if $\Theta^2(T,p) \in \N$, we get $\Theta^2(T,p)=i$. Now, assume by contradiction that there is $q\in E_i$ such that $\Theta^2(T, q)<\Theta^2(T,p)$, then, by \eqref{eq:stratification:density-estimate}, we necessarily have $i - \frac12 \leq \Theta^2(T,q) < i+\frac12$ which ensures, by the classification of tangent cones, $\Theta^2(T,q) = Q - \frac{Q^\star +1}{2}$. Since $\Theta^2(T, q) = Q^\prime - \frac{Q^\star}{2}$ for some $Q^\prime\in\N$, we obtain $Q^\prime = Q - \frac{1}{2}$ which is a contradiction. We then conclude that $p$ is a singular point which is also two-sided collapsed.
    
    \item if $\Theta^2(T,p) \notin \N$, then $\Theta^2(T, p) = i - \frac{1}{2}$ and, since $E_i$ is relatively open, there is a relatively open, in $\Gamma$, neighborhood $U\subset E_i$ of $p$. By definition of $E_i$, for every $q\in U$, $\Theta^2(T, q)\geq i-\frac{1}{2} = \Theta^2(T, p)$ which ensures that $p$ is two-sided collapsed.
\end{itemize}
\end{proof}

We have now reduced our situation to prove the following theorem. 

\begin{theorem}[Two-sided collapsed points are regular]\label{thm:collapsed-is-regular}
Let $T$ and $\Gamma$ be as in Assumption \ref{assumptions} with $C_0 = 0$. Then any two-sided collapsed point of $T$ is a two-sided regular point of $T$.
\end{theorem}

The rest of the paper is devoted to prove Theorem \ref{thm:collapsed-is-regular}, we also mention that we prove this to the general setting of $m$-dimensional area minimizing currents with boundary multiplicity $Q^\star\geq 1$. 

However, our main result (Theorem \ref{thm:main-theorem}) is stated for $2d$ area minimizing currents with boundary multiplicity $Q^\star\geq 1$ because we need to apply Theorem \ref{theorem:not-dense-collapsed-singular} which we proved only in this setting. We recall that Theorem \ref{theorem:not-dense-collapsed-singular} is restricted to the $2d$ case due to the fact that the classification of area minimizing cones (Proposition \ref{prop:2d-tangent-cone}) is only known for $2d$ cones.

\section{Approximations of currents by multi-valued collapsed Dirichlet minimizers}\label{sec:linear}

\subsection{Definitions and regularity of collapsed Dirichlet minimizers}\label{sec:linear1}
We refer the reader to \cite{DS3} and \cite{DS2} for standard definitions and notations about the theory of multiple valued functions. Throughout all this section we will consider an open set $\Omega \subset\R^{m}$ together with a $(m-1)$-submanifold $\gamma$ of class $C^{3,\alpha}$ dividing $\Omega$ in two disjoint open sets $\Omega^{+}$ and $\Omega^{-}$.

\begin{definition}\label{def:function-interface} 
Let $\varphi\in W^{\frac12,2}(\gamma, \Is{Q^{\star}})$, $Q,Q^{\star}\in\N$, $Q\ge Q^{\star}\ge1$. A \textbf{$(Q-\frac{Q^{\star}}2)$-valued function with interface $(\gamma,\varphi)$}, consists of a pair $(f^+, f^-)$ satisfying the following properties 
\begin{enumerate}[\upshape (i)]
          \item $f^+\in W^{1,2}(\Omega^+, \Is{Q})$, $f^-\in W^{1,2}(\Omega^-, \Is{Q-Q^{\star}})$,
          \item $f^+\restrict\gamma=f^-\restrict\gamma+\varphi$.
\end{enumerate} 
We define the \textbf{Dirichlet energy of $(f^+,f^-)$} as $\D(f^+,f^-,\Omega):=\D(f^+,\Omega^+)+\D(f^-,\Omega^-)$. Such a pair will be called \textbf{$\D$-minimizing in $\Omega$}, if for all $\left(Q-\frac{Q^{\star}}{2}\right)$-valued function $(g^+,g^-)$ with interface $(\gamma, \varphi)$  which agrees with $\left(f^{+}, f^{-}\right)$ outside of a compact set $K \subset\subset\Omega$ satisfies $\D(f^+,f^-,\Omega)\le\D(g^+,g^-, \Omega)$.
\end{definition}

\begin{remark}{Note that when $Q^{\star}$ is an even number, we have $Q-\frac{Q^{\star}}{2} =: k\in\N$ which shows that unfortunately we have an overlapping between the nomenclatures of Definition \ref{def:function-interface} and Almgren's definition of $k$-valued functions. However, since it will not cause any confusion in what follows, we will use this abuse of notation. }
\end{remark}

The interesting case to be treated here is when $Q > Q^{\star} > 1$. When $Q=Q^{\star}=1$, the pair $\left(f^{+}, f^{-}\right)$ consists of a single-valued function $f^{+}$ and its Dir-minimality is equivalent to the harmonicity of $f^{+} .$ The case $Q>Q^{\star} =1$ is studied in \cite[Section 4]{DDHM}. The \emph{one-sided} case, i.e. $Q=Q^{\star}$, has to be treated differently and it is done in dimension $2$ in \cite{DNS}. 

\begin{definition}\label{def:collapses-interface}
Let $(f^+,f^-)$ be a $\left(Q-\frac{Q^{\star}}{2}\right)$-valued function with interface $(\gamma, \varphi)$ and $\varphi = Q^{\star}\a{\hat{\varphi}}$ for a single valued function $\hat{\varphi}$. We say that  $(f^+,f^-)$ \textbf{collapses at the interface}, if $f^{+}\restrict\gamma = Q\a{\hat{\varphi}}$.
\end{definition}

Notice that, $(f^+, f^-)$ satisfy the properties of the preceding definition, if and only if, $f^-\restrict\gamma=(Q-Q^{\star})\a{\hat{\varphi}}$.

\begin{figure}[H]
\centering
\begin{subfigure}{.49\linewidth}
\centering
\captionsetup{width=.9\linewidth}
\includegraphics[width=.99\linewidth]{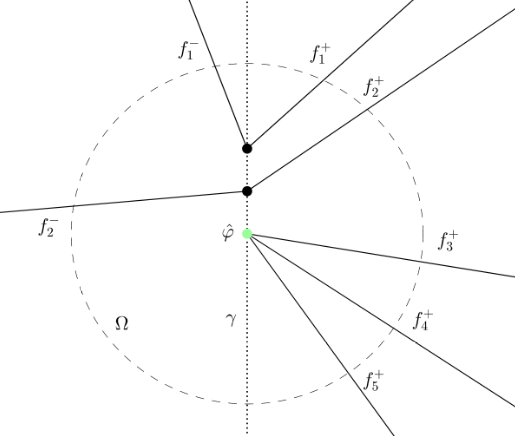}
\caption{$n=1, Q=5, Q^{\star}=3, f^+(x) = \sum_{i=1}^5\a{f_i^+(x)}$ and $f^+(x) = \sum_{i=1}^2\a{f_i^-(x)}$, so that the $(Q-\frac{Q^{\star}}{2})$-valued function $(f^+, f^-)$ has interface $(\gamma, \hat{\varphi})$ where $\gamma = \{x=0\}$ and $(x,\hat{\varphi}(x))$ is constantly equal to the green point.}
\end{subfigure}
\begin{subfigure}{.49\linewidth}
\centering
\captionsetup{width=.9\linewidth}
\includegraphics[width=.99\linewidth]{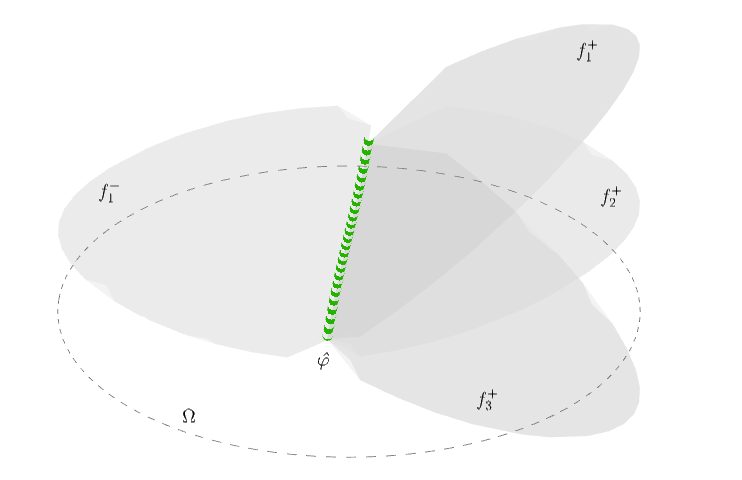}
\caption{Assume $\Omega\subset\R$ and $n=1, Q=3, Q^{\star}=2, f^+(x) = \sum_{i=1}^3\a{f_i^+(x)}$ and $f^+(x) = \a{f_1^-(x)}$, so that the $(Q-\frac{Q^{\star}}{2})$-valued function $(f^+, f^-)$ collapses at the interface $(\gamma, \hat{\varphi})$ where $\gamma = \{x=0\}$ and $\hat{\varphi}$ is represented by the green curve.}
\end{subfigure}
\end{figure} 

With these definitions settled, we aim to prove the harmonic regularity of collapsed $(Q-\frac{Q^\star}{2})$-valued maps along the same lines for $Q^\star =1$ as it is done in \cite[Theorem 4.5]{DDHM}. As we mentioned above, this part of the linear theory in our setting is true only when we consider $Q>Q^\star$ since we will construct some competitors in the arguments which need the existence of multi-valued functions defined in both sides of $\gamma$.

\begin{theorem}[Regularity of collapsing $(Q-\frac{Q^{\star}}{2})$-Dir minimizers]\label{thm:harmonic-regularity-(Q-Q^*/2)-minimizing}
Let $\varphi:\gamma\to\Is{Q^{\star}}$, where $\varphi=Q^{\star}\a{\hat{\varphi}}$ for some $\hat{\varphi}\in C^{1, \alpha}(\gamma, \R^n)$, $\gamma$ be a $(m-1)$-submanifold of class $C^{3}$ in $\R^m$, $Q > Q^{\star}\geq 1$, and $\left(f^{+}, f^{-}\right)$ be $a\left(Q-\frac{Q^{\star}}{2}\right)$-valued Dir-minimizer with interface $(\gamma, \varphi) .$ If $\left(f^{+}, f^{-}\right)$ collapses at the interface, then there is a single-valued harmonic function $h: \Omega\to\R^n$ such that $f^{+}=Q\a{h|_{\Omega^{+}}}$ and $f^{-}=(Q-Q^{\star})\a{h|_{\Omega^-}}$.
\end{theorem}

If we do not assume that the pair $(f^+, f^-)$ is collapsed and impose that $\gamma$ is real analytic, we can obtain that the singular set of the pair is discrete when $Q^\star = 1$, see \cite[Theorem 1.6]{delellis2019dirichlet}. Let us turn to the proof of Theorem \ref{thm:harmonic-regularity-(Q-Q^*/2)-minimizing}, firstly we define the tangent function and then we characterize these tangent function.

\begin{definition}[Tangent function]\label{defi:tangent-function}
Let $\left(f^{+}, f^{-}\right)$ be a $\left(Q-\frac{Q^\star}{2}\right)$-valued function with interface
$(\gamma, Q^\star\a{0})$. Fix $p \in \gamma$ and define a {\textbf{blowup of $f$ at $p$ at scale $r$}} as follows
$$
f_{p, r}^{\pm}(x):=\frac{f^{\pm}(p+r x)}{\sqrt{r^{2-m}\operatorname{Dir}(f^+, f^-, \baseball{p}{r})}}, \forall r>0,
$$
where we assume that $\left(f^{+}, f^{-}\right)$ is not identically $(Q\a{0},(Q-Q^\star)\a{0})$ in every ball $\baseball{0}{r}$. For any sequence $r_k\to 0$, if the limit exists, we say that $g^\pm = \lim_{k\to +\infty}f_{p,r_k}^\pm$ is a {\textbf{tangent function at $p$ to $f$}}.
\end{definition}

\begin{lemma}\label{lemma:linearized-tangent-function}
Let $Q>Q^{\star}$, $\left(f^{+}, f^{-}\right)$ be a $\left(Q-\frac{Q^{\star}}{2}\right)$ Dir-minimizer which collapses at the interface
$(\gamma, Q^{\star}\a{0}),$ where $\gamma$ is a $C^{3}$ $(m-1)$-submanifold in $\R^m$, and fix $p \in \gamma$. Consider a tangent function $\left(h_0^{+}, h_0^{-}\right)$ to $\left(f^{+}, f^{-}\right)$ at $p$ and $\{e_1,\cdots, e_{m-1}\}$ a base of $T_p\gamma$. For each $i\in\{1,\cdots,m-1\}$, we define $(h_i^+, h_i^-)$ to be a tangent to $(h_{i-1}^+, h_{i-1}^-)$ at $e_i$. Then $\left(h^{+}, h^{-}\right):= (h_{m-1}^+, h_{m-1}^-)$ is given by $(Q\a{L},(Q-Q^{\star})\a{L})$
where $L$ is a nonzero linear function which vanishes on $T_{p} \gamma.$
\end{lemma}
\begin{proof}
Assume $T_{p} \gamma=\left\{x: x_{m}=0\right\} .$ The consequences of \cite[Lemma 4.29, Remark 4.31]{DDHM} readily holds in our higher multiplicity case, so we have the following properties:
\begin{enumerate}[\upshape (A)]
\item $\left(h^{+}, h^{-}\right)$ is a $\left(Q-\frac{Q^{\star}}{2}\right)$ Dir-minimizer which collapses at the interface $(T_{p} \gamma, Q^\star\a{0}),$

\item $\left(h^{+}, h^{-}\right)$ depends only on $x_m$ namely there exist $Q$-valued function $\alpha^{+}:\R _{+} \rightarrow$
$\Is{Q}$ and a $(Q-Q^{\star})$-valued function $\alpha^{-}:\R_{-} \rightarrow\Is{Q-Q^{\star}}$ such that $h^{\pm}(x)= \alpha^{\pm}\left(x_m\right),$

\item $\left(h^{+}, h^{-}\right)$ is an $I$-homogeneous function for some $I>0,$ namely there is a $Q$-point $P$ and a $(Q-Q^{\star})$-point $P^{\prime}$ such that $\alpha^{+}\left(x_m\right)=x_m^{I} P$ and $\alpha^{-}\left(x_m\right)=\left(-x_m\right)^{I} P^{\prime},$

\item $\operatorname{Dir}\left(h^{+}, \baseball{0}{1}\right)+\operatorname{Dir}\left(h^{-}, \baseball{0}{1}\right)=1.$
\end{enumerate}
Since $\left(h^{+}, h^{-}\right)$ is a Dir-minimizer and both $h^{+}$ and $h^{-}$ are $C^2$, both $h^{+}$ and $h^{-}$ are classical harmonic functions, therefore, since they depend only upon one variable, we necessarily have that $I=1 .$ So there are coefficients $\beta_{1}^{+}, \ldots, \beta_{Q}^{+}$ and $\beta_{1}^{-}, \ldots, \beta_{Q-Q^{\star}}^{-}$ such that
$$ h^{+}(x)=\sum_{i=1}^{Q}\a{\beta_{i}^{+} x_m}, \text{if} \ x_m > 0, \ \ \text{and} \ \
h^{-}(x)=\sum_{i=1}^{Q-Q^{\star}}\a{\beta_{i}^{-} x_m}, \text{if} \ x_m < 0.$$

If $Q > Q^{\star} > 1,$ then we can assume that $\beta_{j_0}^+ \neq \beta_{i_0}^-$. Now, we will construct a competitor of $(h^+, h^-)$ with less Dir-energy which is the desired contradiction. Note that, in order to construct a competitor, we have to assure that it has the same interface of $(h^+, h^-)$, i.e. it takes $\a{0}$ at least $Q^{\star}$ times at $T_p\gamma = \{ x_m = 0\}$. For $x=(x^\prime, x_m)$, define
$$ \hat{h}^+(x) = \left\{\begin{array}{ll} 
\a{\hat{\beta}x_m + c(|x^\prime|)} + \sum_{j=1, j\neq j_0}^{Q}\a{\beta_j^+ x_m},&\text{if } x\in \overline{\basehalfball{+}{0}{\frac{1}{2}}}, \\
h^+(x),&\text{if } x\in \basehalfball{+}{0}{1}\setminus\overline{\basehalfball{+}{0}{\frac{1}{2}}}.
\end{array}\right.
$$
$$ \hat{h}^-(x) = \left\{\begin{array}{ll} 
\a{\hat{\beta}x_m + c(|x^\prime|)} + \sum_{i=1, i\neq i_0}^{Q-Q^{\star}}\a{\beta_i^+ x_m},&\text{if } x\in \overline{\basehalfball{-}{0}{\frac{1}{2}}}, \\
h^-(x),&\text{if } x\in \basehalfball{-}{0}{1}\setminus\overline{\basehalfball{-}{0}{\frac{1}{2}}}.
\end{array}\right.
$$
where $\hat{\beta} = \frac{\beta_{j_0}^+ + \beta_{i_0}^-}{2}, c(|x^\prime|) = \bar{\beta}\sqrt{1/4 - |x^\prime|^2}$ and $\hat{\beta} = \frac{\beta_{j_0}^+ - \beta_{i_0}^-}{2}$. 

\begin{figure}[h]
\centering
\begin{subfigure}{.49\linewidth}
\centering
\captionsetup{width=.9\linewidth}
\includegraphics[width=\linewidth]{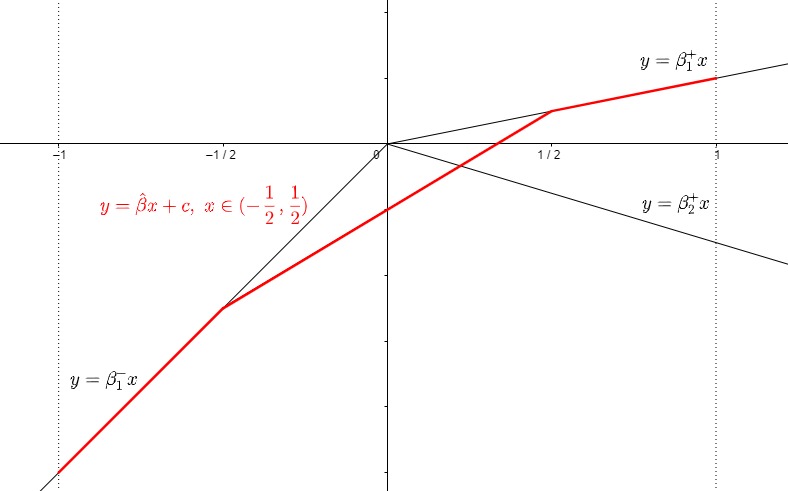}
\caption{if $\Omega\subset\R$ and $n=1, Q=2, Q^{\star}=1, h^+(x)=\a{\beta_1^+x}+\a{\beta_2^+x}$ and $h^-(x)=\a{\beta_1^-x}$, we define the competitor $(\hat{h}^+, \hat{h}^{-})$ as the red function which has the same interface of $(h^+, h^-)$ and less $\operatorname{Dir}$-energy.}
\end{subfigure}
\begin{subfigure}{.49\linewidth}
\centering
\captionsetup{width=.9\linewidth}
\includegraphics[width=\linewidth]{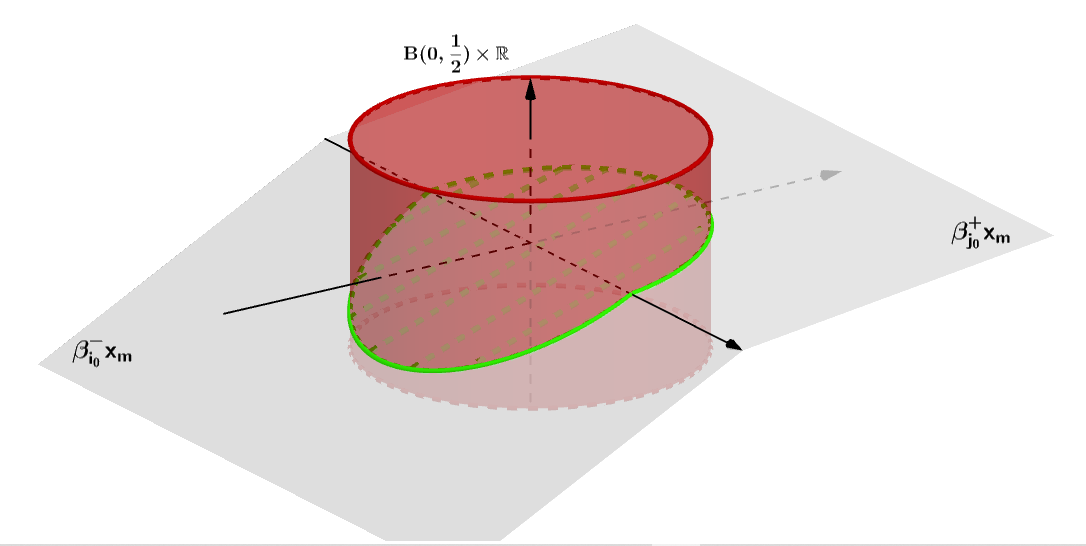}
\caption{With $m=2$ competitor $(\hat{h}^+, \hat{h}^-)$ represented by the green hypersurface is not linear inside the cylinder, but it also satisfies what we need, i.e. it has the same interface of $(h^+, h^-)$ and it has less $\operatorname{Dir}$-energy.}
\end{subfigure}
\end{figure}

By direct computation, we have
$$ \operatorname{Dir}(\hat{h}^+, \basehalfball{+}{0}{1/2}) = |\basehalfball{+}{0}{1/2}|\biggl[ \sum_{j=1, j\neq j_0}^{Q}|\beta_{j_0}^+|^2 + |\hat{\beta}|^2\biggr] + \int_{\mathrm{B}^{m-1}(0, 1 / 2)}\int_{0}^{\sqrt{\frac{1}{4} - |x^\prime|^2}}\frac{|\bar{\beta}|^2|x^\prime|^2}{\frac{1}{4}-|x^\prime|^2}dx'dx_m, $$
the integral on the right hand side can be bounded by $|\bar{\beta}|^2|\basehalfball{+}{0}{1/2}|$ since $\bar{\beta}\neq 0$ and the integrating function is radial. By the very same computation, we finally have that 

$$
\begin{aligned}
\operatorname{Dir}(\hat{h}^+, \hat{h}^-, \baseball{0}{1/2}) &< |\basehalfball{+}{0}{1/2}|\biggl[ \sum_{j=1, j\neq j_0}^{Q}|\beta_{j}^+|^2 + \sum_{i=1, i\neq i_0}^{Q-Q^{\star}}|\beta_{i}^+|^2 + 2|\hat{\beta}|^2 + 2|\bar{\beta}|^2\biggr] \\
&=\operatorname{Dir}(h^+, h^-, \baseball{0}{1/2}).
\end{aligned}
$$

By construction they have the same $\operatorname{Dir}$-energy outside $\baseball{0}{1/2}$, thus every $\beta_j^+$ has to coincide with $\beta_i^-$ and we finish the proof of the lemma.
\end{proof}

\begin{definition}
Let us denote $\boldsymbol{\eta}(P) = \frac{1}{Q}\sum_{i=1}^QP_i$ the center of the $Q$-point $P=\sum_{i=1}^Q\a{P_i}$.
\end{definition}

As a simple corollary of the above lemma we have:

\begin{corollary}\label{cor:center-zero}
Let $Q>Q^{\star}$ and assume $\left(f^{+}, f^{-}\right)$ is a $\left(Q-\frac{Q^{\star}}{2}\right)$ Dir-minimizer which collapses at $(\gamma, Q^\star\a{0}),$ where $\gamma$ is a $C^{3}$ $(m-1)$-submanifold in $\R^m$. If $\mathbf{\boldsymbol{\eta}}\circ f^{-} = \mathbf{\boldsymbol{\eta}}\circ f^{+}=0,$ then $f^{+}=Q\a{0}$ and $f^{-}=(Q-Q^{\star})\a{0}$.
\end{corollary}
\begin{proof}
If $\left(f^{+}, f^{-}\right)$ is identically $(Q\a{0},(Q-Q^{\star})\a{0})$ in a neighborhood $U$ of a point $p \in \gamma,$ then, by the interior regularity theory of Dir-minimizer (precisely, \cite[Proposition 3.22]{DS1}), $\left(f^{+}, f^{-}\right)$ is identically $(Q\a{0},(Q-Q^{\star})\a{0})$  in the connected component of the domain of  $\left(f^{+}, f^{-}\right)$ which contains
$p .$ Thus, if the corollary were false, then there would be a point $p\in\gamma$ such that 
$\operatorname{Dir}\left(f^{+}, \basehalfball{+}{p}{r}\right)+
\operatorname{Dir}\left(f^{-}, \basehalfball{-}{p}{r}\right)>0 \text { for every } r>0$ such that $\baseball{p}{r}\subset\Omega$. If we consider $\left(h^{+}, h^{-}\right)$ as in Lemma \ref{lemma:linearized-tangent-function}, we conclude that $\boldsymbol{\eta} \circ h^{+}=\boldsymbol{\eta} \circ h^{-}=0,$ since such property is inherited by each tangent map. But then the nonzero linear function $L$ of the conclusion of Lemma \ref{lemma:linearized-tangent-function} should be equal to $\boldsymbol{\eta} \circ h^{+}$ on $\left\{x_{m}>0\right\}$ and $\boldsymbol{\eta} \circ h^{-}$ on $\left\{x_{m} \leq 0\right\}$. Hence $L$ should vanish identically, contradicting Lemma \ref{lemma:linearized-tangent-function}.
\end{proof}

Before the proof of Theorem \ref{thm:harmonic-regularity-(Q-Q^*/2)-minimizing}, we introduce the following notation which will be used throughout the paper.
$$ f\oplus \zeta := \sum_{i}\a{f_i +\zeta}, $$
where $f$ is a $Q$-valued function with a measurable selection of single-valued functions $f_i$ and $\zeta$ is a single-valued functions both defined on the same domain.

\begin{proof}[Proof of Theorem \ref{thm:harmonic-regularity-(Q-Q^*/2)-minimizing}]
\textbf{The case $\hat{\varphi}\equiv 0$:} Firstly, using the regularity theory for harmonic functions, we obtain that the functions $\boldsymbol{\eta} \circ f^{\pm}$ are differentiable up to the boundary $\gamma$, i.e., belong to $C^{1}\left(\Omega^{\pm} \cup \gamma\right) .$ Let $\nu$ be the unit normal to $\gamma$. We claim that
\begin{equation}\label{Eq:4.34}
\partial_{\nu}\left( \mathbf{\boldsymbol{\eta}} \circ f^{+}\right)(p)=\partial_{\nu}\left( \boldsymbol{\eta} \circ f^{-}\right)(p) \quad \text { for all } p \in \gamma \cap \Omega.
\end{equation}
In fact, assume by contradiction that, at some point $p \in \gamma \cap \Omega$, we have $\partial_{\nu}\left( \boldsymbol{\eta} \circ f^{+}\right)(p) \neq \partial_{\nu}\left(\boldsymbol{\eta} \circ f^{-}\right)(p)$ and consider a tangent function $\left(h^{+}, h^{-}\right)$ to $\left(f^{+}, f^{-}\right)$ at $p$
which is the limit of some rescaled sequence $\left(f_{p, \rho_{k}}^{+}, f_{p, \rho_{k}}^{-}\right) $, where we denote 

$$f_{p, \rho_{k}}^{\pm}(x):=\frac{f^{\pm}(p+\rho_k x)}{\sqrt{\rho_k^{2-m}\operatorname{Dir}(f^+, f^-, \baseball{p}{\rho_k})}}.$$

Observe that, since at least one among $\partial_{\nu}\left(\boldsymbol{\eta} \circ f^{+}\right)(p)$ and $\partial_{\nu}\left( \boldsymbol{\eta} \circ f^{-}\right)(p)$ differs from $0,$ we necessarily have

$$c_1\rho_k^m \geq \operatorname{Dir}\left(\boldsymbol{\eta}\circ f^{+}, \boldsymbol{\eta}\circ f^-, \baseball{p}{\rho_{k}}\right) \geq c_{0} \rho_{k}^{m},$$

for some constants $c_1 = c_1(\boldsymbol{\eta}\circ f^+,\boldsymbol{\eta}\circ  f^-)>0, c_{0}=c_0(\boldsymbol{\eta}\circ f^+,\boldsymbol{\eta}\circ f^-)>0.$ Thus, by rescaling, we obtain

\begin{equation}\label{eq:new-4.38}
\begin{aligned}
c_1\frac{\rho_k^m}{\operatorname{Dir}\left( f^{+},  f^-, \baseball{p}{\rho_k}\right)} &\geq \frac{\operatorname{Dir}\left(\boldsymbol{\eta}\circ f^{+}, \boldsymbol{\eta}\circ f^-, \baseball{p}{\rho_k}\right)}{\operatorname{Dir}\left( f^{+}, f^-, \baseball{p}{\rho_k}\right)}\\
&= \frac{\operatorname{Dir}\left(\boldsymbol{\eta}\circ f_{p,\rho_k}^{+}, \boldsymbol{\eta}\circ f_{p,\rho_k}^-, \baseball{0}{1}\right)}{\operatorname{Dir}\left( f_{p,\rho_k}^{+}, f_{p,\rho_k}^-, \baseball{0}{1}\right)} \\ &=\operatorname{Dir}\left(\boldsymbol{\eta}\circ f_{p,\rho_k}^{+}, \boldsymbol{\eta}\circ f_{p,\rho_k}^-, \baseball{0}{1}\right) \\
&\geq c_{0} \frac{\rho_{k}^{m}}{\operatorname{Dir}\left( f^{+},  f^-, \baseball{p}{\rho_k}\right)}.
\end{aligned}
\end{equation} 

Therefore, we have the following two alternatives:

\begin{enumerate}[\upshape (I)]
\item If $\lim \sup _{k}\left(\rho_{k}\right)^{-m}\operatorname{Dir}\left(f^{+}, f^-, \baseball{p}{\rho_{k}}\right)=+\infty$, by \eqref{eq:new-4.38}, denoting by $(h_0^+, h_0^-)$ the tangent function to $(f^+, f^-)$ at $p$, passing to the limit in \eqref{eq:new-4.38}, we have that $\operatorname{Dir}(\boldsymbol{\eta}\circ h_0^+, \boldsymbol{\eta}\circ h_0^-, \baseball{0}{1}) = 0$ and then $\boldsymbol{\eta}\circ h_0^\pm\equiv 0$. By Corollary \ref{cor:center-zero}, $\left(h_0^{+}, h_0^{-}\right)$ should be trivial. But this is not possible, because the energy of a tangent function satisfies $\operatorname{Dir}\left(h_0^{+}, h_0^-, \baseball{0}{1}\right) = 1$, see Lemma \ref{lemma:linearized-tangent-function}.

\item If $\lim \sup _{k}\left(\rho_{k}\right)^{-m}\operatorname{Dir}\left(f^{+}, f^-, \baseball{p}{\rho_{k}}\right)<+\infty$, by \eqref{eq:new-4.38}, we have $\operatorname{Dir}(\boldsymbol{\eta}\circ h_0^+, \boldsymbol{\eta}\circ h_0^-, \baseball{0}{1}) > 0$ and thus  $\boldsymbol{\eta} \circ h_0^{+}$ and $\boldsymbol{\eta} \circ h_0^{-}$ are distinct functions and at least one among them is a nontrivial function. Indeed, they are distinct follows from the fact that the blowup of a differentiable function coincides with its differential and we are assuming that $\partial_{\nu}\left(\boldsymbol{\eta} \circ f^{+}\right)(p)\neq \partial_{\nu}\left( \boldsymbol{\eta} \circ f^{-}\right)(p).$ Since case (I) never occurs, we can apply this argument iteratively until we reach the pair $(h^+, h^-)$ of Lemma \ref{lemma:linearized-tangent-function} and then conclude that $\boldsymbol{\eta} \circ h^{+}$ and $\boldsymbol{\eta} \circ h^{-}$ are two distinct linear functions with one of them being non trivial and this contradicts Lemma \ref{lemma:linearized-tangent-function}.
\end{enumerate}

We have verified the validity of \eqref{Eq:4.34} and it is enough to conclude the proof, indeed, it implies that the function
\begin{equation}\label{Eq:4.35}
\zeta :=\left\{\begin{array}{ll}
\boldsymbol{\eta} \circ f^{+} & \text {on } \Omega^{+}, \\
\boldsymbol{\eta} \circ f^{-} & \text {on } \Omega^{-}.
\end{array}\right.
\end{equation}
is an harmonic function defined on the entire $\Omega$. Using the notation above, we set

\begin{equation}\label{Eq:4.36}
\tilde{f}^{+}:=f^+\oplus(-\zeta) \ \text{and} \
\tilde{f}^{-}:=f^-\oplus(-\zeta).
\end{equation}

By \cite[Lemma 3.23]{DS1}, it is easy to see that $\left(\tilde{f}^{+}, \tilde{f}^{-}\right)$ is a $\left(Q-\frac{Q^{\star}}{2}\right)$ Dir-minimizer which collapses at the interface $(\gamma, \a{0})$ and that $\boldsymbol{\eta} \circ \tilde{f}^{+}=\boldsymbol{\eta} \circ \tilde{f}^{-}=0 .$ Thus we apply Corollary \ref{cor:center-zero} and conclude that $\tilde{f}^{+}=Q\a{0}$ and $\tilde{f}^{-}=(Q-Q^{\star})\a{0},$ which complete the proof.

\textbf{The general case:}
We fix $\nu$ as an unit normal to $\gamma$. As in the particular case $\hat{\varphi}\equiv 0$, we claim that $\partial_{\nu}\left( \boldsymbol{\eta} \circ f^{+}\right)=\partial_{\nu}\left( \boldsymbol{\eta} \circ f^{-}\right)$. With this claim, proceeding as in the former case, we can define $\zeta$ as in \eqref{Eq:4.35} and conclude that it is a harmonic function. We then define $\left(\tilde{f}^{+}, \tilde{f}^{-}\right)$ as in \eqref{Eq:4.36}. To this pair, we can apply the former case and conclude the proof of the theorem. To prove the claim, assume by contradiction that, for some $p \in \gamma,$ we have that $\partial_{\nu}(\boldsymbol{\eta} \circ f^{+}(p)) \neq \partial_{\nu}\left( \boldsymbol{\eta} \circ f^{-}\right)(p)$. Without loss of generality we can assume that $p=0, \hat{\varphi}(0)=0$
and $D \hat{\varphi}(0)=0.$ Since at least one between $\partial_{\nu}(\boldsymbol{\eta}\circ f^{\pm})(0)$ does not vanish, we must have

\begin{equation}\label{4.38}
\operatorname{Dir}\left(f^{+}, f^-, \baseball{0}{\rho}\right) \geq \operatorname{Dir}\left(\boldsymbol{\eta}\circ f^{+}, \boldsymbol{\eta}\circ f^-, \baseball{0}{\rho}\right) \geq c_{0} \rho^{m},
\end{equation}

for some positive constant $c_{0} .$ It also means that there exist a constant $\eta>0$ and a sequence $\rho_{k} \downarrow 0$ such that

$$
\operatorname{Dir}\left(f^{+}, f^-, \baseball{0}{\rho_{k}}\right) \geq \eta\left(\operatorname{Dir}\left(f^{+}, f^-, \baseball{0}{2\rho_{k}}\right)\right),
$$

otherwise we would contradict the lower bound \eqref{4.38}. We see that $f_{0, \rho_{k}}^{\pm}$ have finite energy on $\baseball{0}{2}$ and thus there is strong convergence of a subsequence to a $\left(Q-\frac{Q^{\star}}{2}\right)$ Dir-minimizer $\left(h^{+}, h^{-}\right)$ with interface $\left(T_{0} \gamma, Q^\star\a{0}\right) .$ The latter must then have Dirichlet energy 1 on $\baseball{0}{1}$. We then have two possibilities:
\begin{enumerate}[\upshape (I)]
\item $\lim \sup _{k}\left(\rho_{k}\right)^{-m}\operatorname{Dir}\left(f^{+}, f^-, \baseball{0}{\rho_{k}}\right)=+\infty$. Arguing as in (I) of the former case, this gives that $\boldsymbol{\eta} \circ h^{+}=\boldsymbol{\eta} \circ h^{-}=0$ and thus we conclude that $\left(h_0^{+}, h_0^{-}\right)$ is trivial, which is a contradiction,

\item $\lim \sup _{k}\left(\rho_{k}\right)^{-m}\operatorname{Dir}\left(f^{+}, f^-, \baseball{0}{\rho_{k}}\right)<+\infty$. Assuming in this case that
$T_{0} \gamma=\left\{x_{m}=0\right\},$ we conclude that $\left(h_0^{+}, h_0^{-}\right)$ is a $\left(Q-\frac{Q^{\star}}{2}\right)$ Dir-minimizer with flat interface $\left(T_{0} \gamma, Q^\star\a{0}\right),$ but also that $\boldsymbol{\eta} \circ h_0^{\pm}(x)= C_d \partial_{\nu}\left(\boldsymbol{\eta} \circ f^{\pm}\right)(0) x_{m}$ for some positive
constant $C_d$ that in general is not necessarily equal to one, because we use a normalization constant to have $\operatorname{Dir}(h_0^+,h_0^-, \baseball{0}{1})=1$. By the particular case $\hat{\varphi}\equiv 0$, we thus conclude that $\partial_{\nu}\left( \boldsymbol{\eta} \circ f^{+}\right)(0)=\partial_{\nu}\left( \boldsymbol{\eta} \circ f^{-}\right)(0)$.
\end{enumerate}
\end{proof}

\subsection{Harmonic approximations}

In this chapter we aim to approximate the area minimizing current $T$ by $Q$ copies of an harmonic function in the right side and $Q-Q^\star$ copies of the same harmonic function in the left side. To this end, we will first approximate the current $T$ by $(Q-\frac{Q^\star}{2})$-Lipschitz functions which, if we do not assume the minimality of $T$, will not be necessarily minimizers for the Dirichlet energy. Once we consider the minimality condition on $T$ we will be able to upgrade our approximations using the regularity theorem, see Theorem \ref{thm:harmonic-regularity-(Q-Q^*/2)-minimizing}, to furnish the desired harmonic approximations. 

For any $\pi,\pi_0$ belonging to $\mathrm{G}_{m,m+n}$, where $\mathrm{G}_{k,l}$ denotes the set of $k$-dimensional subspaces of $\R^{l}$, we introduce, for any $p\in\R^{m+n}$ the notation $\mathrm{B}_r(p, \pi)$ for the disks $\ball{p}{r} \cap(p+\pi)$, if $\pi$ is omitted, then we assume $\pi=\pi_{0}=\R^m\times\{0\}$, and $\cyltilted{p}{r}{\pi}$ for the cylinders $\mathrm{B}_r(p, \pi)+\pi^{\perp}$, we also fix $\cyl{x}{r} := \cyltilted{p}{r}{\pi_0}$.

\begin{definition}\label{def:excess-flat}
Let $\left.\left. \alpha \in \right]0,1\right]$ and integers $m\geq 2$, $Q^\star\geq 1$, and take $\Gamma$ any $(m-1)$-rectifiable set. Let $T$ be an $m$-dimensional integral current with boundary $\partial T ={Q^\star}\a{\Gamma}$ and assume that $p\in\Gamma$. Then
\begin{enumerate}[\upshape (i)]
    \item We call the \textbf{cylindrical excess relative to the plane $\pi$} the quantity 
    $$\bE (T, \cyl{p}{r}, \pi) := \frac{1}{\omega_m r^m} \int_{\cyl{p}{r}} \frac{|\vec{T} (x) - \vec{\pi}|^2}{2} \, d\|T\| (x),$$ 
    and the {\textbf{cylindrical excess}} the quantity 
    $$\bE (T, \cyl{p}{r}) := \min \{\bE (T, \cyl{p}{r}, \pi): \pi\subset\R^{m+n}\}.$$
    
    \item We call the \textbf{spherical excess relative to the plane $\pi$} the quantity 
    $$\bE (T, \ball{p}{r}, \pi) := \frac{1}{\omega_m r^m} \int_{\ball{p}{r}} \frac{|\vec{T} (x) - \vec{\pi}|^2}{2} \, d\|T\| (x),$$ 
    and the {\textbf{spherical excess}} the quantity 
    $$\bE (T, \ball{p}{r}) := \min \{\bE (T, \ball{p}{r}, \pi): \pi\subset\R^{m+n}\}.$$
    
    \item We say that the {\textbf{boundary spherical excess}} is the quantity $$\bE^\flat(T,\ball{p}{r}):=\min\left\{\bE(T,\ball{p}{r},\pi):\,T_{p}\Gamma\subset\pi\subset\R^{m+n}\right\} .$$
    
    \item The \textbf{height of $T$ in a set $G \subset \mathbb{R}^{m+n}$ with respect to a plane $\pi$} is defined as 
    $$ \mathbf{h}(T, G, \pi):=\operatorname{diam}(\bp_{\pi}^\perp(\spt(T)\cap G)) = \sup \left\{\left|\mathbf{p}_{\pi}^{\perp}(q-p)\right|: q, p \in \operatorname{spt}(T) \cap G\right\},$$
    where $\mathbf{p}_{\pi}^{\perp}$ denotes the orthogonal projection onto $\pi^\perp$.
    
    \item If $\spt(T)\subset\cyltilted{p}{r}{\pi}$, we define the \textbf{excess measure with respect to $\cyltilted{p}{r}{\pi}$} as the measure which to each $F\subset\mathrm{B}_r(p,\pi)$ gives
    $$ \mathbf{e}_{T}(F):=\frac{1}{2} \int_{F+\pi^{\perp}}|\vec{T}-\vec{\pi}|^{2} d\|T\|.$$
\end{enumerate}
\end{definition}

In this subsection we assume that $\pi_{0}=\mathbb{R}^{m} \times\{0\}$ and we use the notation $\mathbf{p}$ and $\mathbf{p}^{\perp}$ for the orthogonal projections onto $\pi_{0}$ and $\pi_{0}^{\perp}$ respectively, whereas $\mathbf{p}_\pi$ and $\mathbf{p}_{\pi}^{\perp}$ will denote, respectively, the orthogonal projections onto the plane $\pi$ and its orthogonal complement $\pi^{\perp}$. For the remaining part of this work, we will call \textbf{dimensional constants} those which depends only on $m,n, Q^\star$ and $Q$.

\begin{assumption}\label{assump:general.asump-first-lip}
Let $\alpha \in (0,1]$ and integers $m\geq 2$, $Q^\star\geq 1$. Consider $\Gamma$ a $C^{2,\alpha}$ oriented $(m-1)$-submanifold without boundary. Let $T$ be an $m$-dimensional integral current in $\ball{0}{2}$ with boundary $\partial T \res \ball{0}{2}={Q^\star}\a{\Gamma \cap \ball{0}{2}}$ and assume that $p\in\Gamma$. We also assume $T_{p} \Gamma=\mathbb{R}^{m-1} \times\{0\} \subset \pi_{0}, \psi_1:\R^{m-1}\to\R, \psi:\gamma\subset\R^{m}\times\{0\}\rightarrow\mathbb{R}^{n},\psi_2:\R^{m-1}\rightarrow\mathbb{R}^{n+1}, \psi_2(x) = (\psi_1(x), \psi(x, \psi_1(x)))$ with $\gr{\psi_1} = \gamma$ and $\Gamma=\gr{\psi_2}$ satisfying the bounds $\|D \psi_2\|_{0, \baseball{0}{2}} \leq c_{0}$ and $\bA:=\left\|A_{\Gamma}\right\|_{0,\ball{0}{2}} \leq c_{0}$, where $A_{\Gamma}$ denotes the second fundamental form of $\Gamma$ and $c_{0}$ is a positive small geometric constant. We assume that

\begin{enumerate}[\upshape (i)]
\item\label{assump-first-lip-apprx:item-density} $p\in \Gamma$ is a two-sided collapsed point with $Q-\frac{Q^{\star}}{2} = \Theta^m(T, p)$, for some integers $Q >Q^\star \geq 1$,

\item $\gamma=\mathbf{p}(\Gamma)$ divides $\baseball{p}{4 r}\subset \pi_0$ in two disjoint open sets $\Omega^{+}$ and $\Omega^{-} ,$

\item\label{def:item-projection-condition}  $\mathbf{p}_{\sharp} T = Q \a{\Omega^{+}}+(Q-Q^{\star})\a{\Omega^{-}}$.
\end{enumerate}
\end{assumption}

Observe that thanks to \eqref{def:item-projection-condition} we have the identities

\begin{equation}\label{eq:excess-and-excessmeasure-formula}
\begin{aligned}
\bE\left(T, \cyl{p}{4 r}\right) &=\frac{1}{\omega_{m} (4r)^{m}}\left(\|T\|\left(\cyl{p}{4 r}\right)-\left(Q\left|\Omega^{+}\right|+(Q-Q^{\star})\left|\Omega^{-}\right|\right)\right), \\
\be_{T}(F) &=\|T\|\left(F \times \mathbb{R}^{n}\right)-\left(Q\left|\Omega^{+} \cap F\right|+(Q-Q^{\star})\left|\Omega^{-} \cap F\right|\right).
\end{aligned}
\end{equation}

\begin{definition}\label{def:maximal-function-excessmeasure}
Given a current $T$ in a cylinder $\cyltilted{p}{4 r}{\pi},$ we introduce the noncentered maximal function of $\mathbf{e}_{T}$ as
$$
\mathbf{me}_{T}(y):=\sup _{y \in \mathrm{B}_s(z,\pi) \subset \mathrm{B}_{4 r}(p, \pi)} \frac{\mathbf{e}_{T}\left(\mathrm{B}_s(z, \pi)\right)}{\omega_{m} s^{m}}.
$$
\end{definition}

The following theorem allows us to approximate the current by a $(Q-\frac{Q^{\star}}{2})$-Lipschitz map which coincides with the current in a closed set $K$ which is called the \textbf{good set}. Moreover, we prove that the \textbf{bad set}, i.e. $\baseball{0}{3r}\setminus K$, has small measure. The tricky part of this theorem is to show that we can take such approximation collapsing at the interface. Notice that, no minimality condition are being assumed to prove this result.

\begin{theorem}[Lusin type weak Lipschitz approximation]\label{thm:first-lip-approx-and-good-set} 
There are positive geometric constants $C= C(m,n,Q,Q^\star)$ and $c_0=c_0(m,n,Q^\star,Q)$ with the following properties. Assume $T$ satisfies Assumption \ref{assump:general.asump-first-lip}, it is area minimizing, and $\bE\left(T, \cyl{p}{4 r}\right) \leq c_{0}$. Then, for any $\delta_{*} \in(0,1)$, there are a closed set $K \subset \baseball{\bp(p)}{3r}$ and a $\left(Q-\frac{Q^{\star}}{2}\right)$ valued function $\left(u^{+}, u^{-}\right)$ on $\baseball{\bp(p)}{3r}$ which collapses at the interface $(\gamma, Q^\star\a{\psi})$ satisfying the following properties:

\begin{align}
\operatorname{Lip}\left(u^{\pm}\right) & \leq C\left(\delta_{*}^{1 / 2}+r^{\frac{1}{2}} \mathbf{A}^{\frac{1}{2}}\right),\label{eq:thmlipaprox:bound-lipconstant} 
\\ 
\operatorname{osc}\left(u^{\pm}\right) & \leq C \mathbf{h}\left(T, \cyl{p}{4r}, \pi_{0}\right)+C r \bE\left(T, \cyl{p}{4 r}\right)^{1 / 2}+C r^{2}\mathbf{A} , \label{eq:thmlipaprox:bound-osc}  
\\
K & \subset \baseball{\bp(p)}{3r} \cap\left\{\mathbf{me}_{T} \leq \delta_{*}\right\},\label{eq:thmlipaprox:goodset-contained}
\\
\mathbf{G}_{u^\pm} \res\left[\left(K \cap \Omega^{\pm}\right) \times \mathbb{R}^{n}\right] &= T \res\left[\left(K \cap \Omega^{\pm}\right) \times \mathbb{R}^{n}\right], \label{eq:thmlipaprox:T=graph(u)}
\\
\left|\baseball{\bp(p)}{s} \backslash K\right| & \leq \frac{C}{\delta_{*}} \mathbf{e}_{T}\left(\left\{\mathbf{me}_{T}>\delta_{*}\right\} \cap \baseball{\bp(p)}{s+r_{1} r}\right), \quad \forall s \leq\left(3-r_{1}\right) r, \label{eq:thmlipaprox:badset-estimate}
\\
\frac{\left\|T-\mathbf{G}_{u^{+}}-\mathbf{G}_{u^-}\right\|\left(\cyl{p}{3 r}\right)}{r^{m}} & \leq \frac{C}{\delta_{*}} \bE\left(T, \cyl{p}{4 r}\right),\label{eq:thmlipaprox:T-graph(u)}
\end{align}

where $r_{1}=c_0 \sqrt[m]{\frac{\bE\left(T, \cyl{p}{4 r}\right)}{\delta_{*}}}$.
\end{theorem}
\begin{proof}
The proof of this theorem is a straightforward adaptation of the corresponding statement considering multiplicity $Q^\star = 1$, c.f. \cite[Theorem 5.5]{DDHM}.
\end{proof}

From now on the approximation of Theorem \ref{thm:first-lip-approx-and-good-set} is called the \textbf{$\delta_{*}^{\frac{1}{2}}$ -approximation of $T$ in $\cyl{p}{3 r}$}. If $E:=\bE\left(T, \cyl{p}{4 r}\right)$, actually in the sequel we will choose $\delta_{*}^{\frac{1}{2}}$ to be $E^{\beta}$ for a small suitable constant $\beta$. In the following theorem, we will use the minimality condition on $T$ to prove that, if $E$ is taken sufficiently small, then $(u^+, u^-)$ is close to a minimizer of the Dirichlet energy, i.e., a $\left(Q-\frac{Q^{\star}}{2}\right)$-Dir-minimizer, which collapses at its interface and thus, by Theorem \ref{thm:harmonic-regularity-(Q-Q^*/2)-minimizing}, consists of a single harmonic sheet.

\begin{theorem}[Harmonic approximation]\label{thm:first-lip-approx-by-harmonic-sheet}
For every $\eta_{*}>0$ and every $\beta \in\left(0, \frac{1}{4 m}\right)$ there exist constants $\varepsilon= \varepsilon(m,n,Q^\star,Q,\eta_*, \beta)>0$ and $C= C(m,n,Q^\star,Q,\eta_*, \beta)>0$ with the following property. Let $T$ and $\Gamma$ be as in Assumption \ref{assumptions} with $C_0=0$ and under the conditions of Theorem \ref{thm:first-lip-approx-and-good-set}, $E \leq c_0$, let $\left(u^{+}, u^{-}\right)$ be the $E^{\beta}$-approximation of $T$ in $\baseball{\bp(p)}{3r}$ and let $K$ be the good set satisfying all the properties \eqref{eq:thmlipaprox:bound-lipconstant}-\eqref{eq:thmlipaprox:T-graph(u)}. If $E\leq \varepsilon$ and $r \mathbf{A} \leq \varepsilon E^{\frac{1}{2}}$, then

\begin{equation}\label{eq:thm-first-lip-approx-by-harm-sheet:excees-measure-of-bad-set}
    \mathbf{e}_{T}\left(\baseball{\bp(p)}{5r / 2} \backslash K\right) \leq \eta_{*} E,
\end{equation}

and

\begin{equation}\label{eq:thm-first-lip-approx-by-harm-sheet:Dir-energy-on-bad-set}
    \operatorname{Dir}\left(u^{+}, u^{-}, \Omega \cap \baseball{\bp(p)}{2r} \backslash K\right) \leq C \eta_{*} E .
\end{equation}

Moreover, set $Q^+:=Q$ and $Q^-:=Q-Q^{\star}$, there exists an harmonic function $h: \baseball{\bp(p)}{2r} \rightarrow \mathbb{R}^{n}$ such that $\left.h\right|_{\{x_{m}=0\}} \equiv 0$ and satisfies the following inequalities:
\begin{align}
r^{-2} \int_{\baseball{\bp(p)}{2r} \cap \Omega^+} \mathcal{G}\left(u^{+}, Q^+ \a{h}\right)^{2}+\int_{\baseball{\bp(p)}{2r} \cap \Omega^{+}}\left(\left|D u^{+}\right|-\sqrt{Q^+}|D h|\right)^{2} &\leq \eta_{*} E r^{m}, \label{eq:thm-harm-aprox:estimate-u^+}
\\
r^{-2} \int_{\baseball{\bp(p)}{2r} \cap \Omega^{-}} \mathcal{G}\left(u^{-},Q^-\a{h}\right)^{2}+\int_{\baseball{\bp(p)}{2r} \cap \Omega^{-}}\left(\left|D u^{-}\right|-\sqrt{Q^-}|D h|\right)^{2} &\leq \eta_{*} E r^{m}, \label{eq:thm-harm-aprox:estimate-u^-}
\\
\int_{\baseball{\bp(p)}{2r} \cap \Omega^\pm}\left|D\left(\boldsymbol{\eta} \circ u^{\pm}\right)-D h\right|^{2} &\leq \eta_{*} E r^{m}. \label{eq:thm-harm-aprox:estimate-center-of-u}
\end{align}
\end{theorem}
\begin{proof}
Without loss of generality we assume that $p=0$, $r=1$, and $\psi(0)=0$.

\textbf{Proof of \eqref{eq:thm-first-lip-approx-by-harm-sheet:excees-measure-of-bad-set} and \eqref{eq:thm-first-lip-approx-by-harm-sheet:Dir-energy-on-bad-set}}. Firstly we want to note that \eqref{eq:thm-first-lip-approx-by-harm-sheet:Dir-energy-on-bad-set} is a consequence of \eqref{eq:thm-first-lip-approx-by-harm-sheet:excees-measure-of-bad-set}. Indeed, since, $\delta_{*} = E^{2\beta}$, use first \eqref{eq:thmlipaprox:goodset-contained}, \eqref{eq:thmlipaprox:badset-estimate} and \eqref{eq:thm-first-lip-approx-by-harm-sheet:excees-measure-of-bad-set} to estimate
$$
\left|\baseball{0}{2} \backslash K\right| \leq C \eta_{*} E^{1-2 \beta}.
$$
Since $\operatorname{Lip}\left(u^{\pm}\right) \leq C E^{2 \beta}$, \eqref{eq:thm-first-lip-approx-by-harm-sheet:Dir-energy-on-bad-set} follows easily. We will argue by contradiction to prove \eqref{eq:thm-first-lip-approx-by-harm-sheet:excees-measure-of-bad-set}. Assuming that the statement is false, there exist $\beta \in (0, 1/4m),\eta_{*}>0$ and a sequence of area minimizing currents $T_{k}$ and submanifolds $\Gamma_{k}$ as in Assumption \ref{assump:general.asump-first-lip} satisfying the following properties for all $k\in\N$:  

\begin{enumerate}[\upshape (i)]
    \item\label{eq:first-lip-approx:contradicition1} The cylindrical excesses $E_{k}:=\mathbf{E}\left(T_{k}, \cyl{0}{4}, \pi_{0}\right)$ satisfy $E_k\leq\frac{1}{k}$,  
    
    \item\label{eq:first-lip-approx:contradicition2} $\Gamma_{k}$ is the graph of the entire function ${\psi_2}_{k}: \mathbb{R}^{m-1} \rightarrow \mathbb{R}^{n+1}$ satisfying the bound
    \begin{equation}\label{eq:lip-approx-by-harm:contradictory-assump-geometry-of-manifolds} 
    \left\|{\psi_2}_{k}\right\|_{C^{2}\left(\baseball{0}{8}\right)} \leq C \mathbf{A}_{k} \leq \frac{C}{k} E_{k}^{1 / 2},
    \end{equation}
    
    \item\label{eq:first-lip-approx:contradicition3} The estimate \eqref{eq:thm-first-lip-approx-by-harm-sheet:excees-measure-of-bad-set} fails, i.e.,  for some positive $c_{2}>0$,
    
    \begin{equation}\label{eq:lip-approx-by-harm:contradiction-hyp-of-excess-measure}  \mathbf{e}_{T_{k}}\left(\baseball{0}{5 / 2} \backslash K_{k}\right)>\eta_{*} E_{k} = 3 c_{2} E_{k}.
    \end{equation}

\end{enumerate}

The pair $\left(f_{k}^{+}, f_{k}^{-}\right)$ are $\left(Q-\frac{Q^{\star}}{2}\right)$-valued maps defined on $\baseball{0}{3}$ which collapses at its interface $(\gamma_k, Q^\star\a{\psi_k})$ denotes the $E_{k}^{\beta}$-Lipschitz approximations of the current $T_{k}$ and $K_k$ the corresponding good set. We denote by $B^\pm_{k, r}$ the domains of the functions $f_{k}^{\pm}$ intersected with the ball $\baseball{0}{r} \subset \pi_{0}$. From the Taylor expansion of the area functional, arguing as in \cite[Remark 5.5]{DS3}, since $E_{k} \downarrow 0$, we conclude the following inequalities for every $s \in[5 / 2,3]:$

\begin{equation}\label{eq:lip-approx-by-harm:eq-that-gives-contradiction}
\begin{aligned}
\int_{K_{k} \cap B^+_{k, s}} \frac{\left|D f_{k}^{+}\right|^{2}}{2}+\int_{K_{k} \cap B_{k, s}^-} \frac{\left|D f_{k}^{-}\right|^{2}}{2} \overset{\textup{Taylor}}&{\leq}\left(1+C E_{k}^{2 \beta}\right) \mathbf{e}_{T_{k}}\left(K_{k} \cap \baseball{0}{s}\right) \\
\overset{\eqref{eq:lip-approx-by-harm:contradiction-hyp-of-excess-measure}}&{<}\left(1+C E_{k}^{2 \beta}\right)\left(\mathbf{e}_{T_{k}}\left(\baseball{0}{s}\right)-3 c_{2} E_{k}\right) \\
& \leq \mathbf{e}_{T_{k}}\left(\baseball{0}{s}\right)-2 c_{2} E_{k} .
\end{aligned}
\end{equation}

The last inequality holds when $E_k$ is sufficiently small, i.e., $k$ large enough. The rest of the proof is devoted to show that \eqref{eq:lip-approx-by-harm:eq-that-gives-contradiction} contradicts the minimizing property of $T_{k}$. We have

\begin{equation}\label{eq:lip-approx-by-harm:bound-dir-energy-f_k}
\begin{aligned}
    \operatorname{Dir}\left(f_{k}^{+},f_k^{-},B_{k,3}\right) &\leq  \operatorname{Dir}\left(f_{k}^{+},f_k^{-},B_{k,3}\cap K_k\right) + \operatorname{Dir}\left(f_{k}^{+},f_k^{-},B_{k,3}\setminus K_k\right)\\
    \overset{\eqref{eq:lip-approx-by-harm:eq-that-gives-contradiction}}&{\leq}  \mathbf{e}_{T_{k}}\left(\baseball{0}{3}\right)-2 c_{2} E_{k} + \operatorname{Dir}\left(f_{k}^{+},f_k^{-},B_{k,3}\setminus K_k\right)\\
    &\leq \mathbf{e}_{T_{k}}\left(\baseball{0}{3}\right)-2 c_{2} E_{k} + \frac{C}{2}E_k^{1-2\beta+2\beta},
\end{aligned}
\end{equation}

where in the last inequality we use the fact that $\operatorname{Lip}\left(f_{k}^{\pm}\right) \leq C E_{k}^{\beta}$ and $\left|\baseball{0}{3} \backslash K_{k}\right| \leq C E_{k}^{1-2 \beta}$. Now we define $\left(g_{k}^{+}, g_{k}^{-}\right)$ as $g_{k}^{\pm}:=E_{k}^{-\frac{1}{2}} f_{k}^{\pm}$ which also collapses at the interface
$\left(\gamma_{k}, E_k^{-\frac{1}{2}}Q^\star\a{\psi_{k}}\right)$. Then if we take $\gamma$ to be the plane $\left\{x_{m}=0\right\} \subset \pi_{0}$, by \eqref{eq:lip-approx-by-harm:contradictory-assump-geometry-of-manifolds}, we obtain the following convergences

\begin{equation}\label{eq:first-lip-approx:C^1-convergeneces}
    \gamma_{k} \overset{C^1}{\longrightarrow} \gamma, \quad \psi_{k} \overset{C^1}{\longrightarrow} 0.
\end{equation}

We now want to do an argument based on the interpolation of the sequence $g_k^\pm$ using the Interpolation Lemma, c.f. \cite[Lemma 4.9]{DDHM} (we mention that this theorems works line by line in higher multiplicity), but unfortunately they do not have the same interface. To overcome this difficulty, we do the following construction. For each $k$, we let $\Phi_{k}$ be a diffeomorphism which maps $\baseball{0}{3}$ onto itself and $\gamma_{k} \cap \baseball{0}{3}$ onto $\gamma \cap \baseball{0}{3}$. By the $C^1$ convergences above for $k$ large enough, it is not difficult to see that we can construct without loss of generality 

$$\left\|\Phi_{k}-\mathrm{Id}\right\|_{C^{1}} \rightarrow 0, \quad \Phi_{k}\left(\partial \baseball{0}{r}\right)=\partial \baseball{0}{r}, \quad \forall r \in[2,3].$$ 

Furthermore, we have that $\left\|\psi_{k} \circ \Phi_{k}^{-1}\right\|_{C^{1}\left(\gamma\right)} \rightarrow 0$. Now consider, for every $x=(x',x_m)\in\R^{m}$, then we define $\kappa_k\in C^1(\baseball{0}{3})$ as follows $\kappa_k(x) = (\psi_{k} \circ \Phi_{k}^{-1})(x',0) + x_m$ then we have that $\kappa_{k} := \psi_{k} \circ \Phi_{k}^{-1}$ on $\gamma$ and $\left\|\kappa_k\right\|_{C^{1}\left(\baseball{0}{3}\right)} \rightarrow 0.$ We fix a measurable selection $f_{k}^{\pm}(x)=\sum_{i} \a{\left(f_{k}^{\pm}\right)_{i}(x)}$, and we set 

$$\hat{g}_{k}^{\pm}:=\left(g_{k}^{\pm}\circ \Phi_{k}^{-1}\right)\oplus(-\kappa_{k}),$$

thus $\left(\hat{g}_{k}^{+}, \hat{g}_{k}^{-}\right)$ are $\left(Q-\frac{Q^{\star}}{2}\right)$-valued maps which collapses at the same interface $(\gamma, Q^\star\a{0})$ and by straightforward computations

\begin{equation}\label{eq:lip-approx-by-harm:relation-Dir-energy-error1}
    \operatorname{Dir}\left(\hat{g}_{k}^{\pm}, \Phi_{k}^{-1}(A) \cap B_{k,3}^{\pm}\right)=(1+o(1))\left(\operatorname{Dir}\left(g_{k}^{+}, A \cap B_{k,3}^{\pm}\right)+\operatorname{Dir}\left(g_{k}^{-}, A \cap B_{k,3}^{\pm}\right)\right)+o(1),
\end{equation}

for all measurable $A \subset \baseball{0}{3}$ and $o(1)$ is independent of the set $A$. From the very definition of $g_k^\pm$ and \eqref{eq:lip-approx-by-harm:bound-dir-energy-f_k}, we conclude that the Dirichlet energy of $\left(\hat{g}_{k}^{+}, \hat{g}_{k}^{-}\right)$ is uniformly bounded. From this bound and \eqref{eq:first-lip-approx:C^1-convergeneces} we may apply the compactness theorem, see \cite[Theorem 4.8]{DDHM}, we can find a not relabeled subsequence and a $\left(Q-\frac{Q^{\star}}{2}\right)$-valued map $\left(g^{+}, g^{-}\right)$ with interface $(\gamma, Q^\star\a{0})$ such that $\left\|\mathcal{G}\left(\hat{g}_{k}^{\pm}, g^{\pm}\right)\right\|_{L^{2}\left(\basehalfball{\pm}{0}{3}\right)} \rightarrow 0$ and

$$
\operatorname{Dir}\left(g^{+}, g^{-}\right) \leq \liminf _{k \rightarrow \infty} \operatorname{Dir}\left(\hat{g}_{k}^{+}, \hat{g}_{k}^{-}\right) \overset{\eqref{eq:lip-approx-by-harm:relation-Dir-energy-error1}}{=} \liminf _{k \rightarrow \infty}\operatorname{Dir}\left(g_{k}^{+}, g_{k}^{-}\right).
$$

Moreover, up to extracting a subsequence, we can assume that $\left|D \hat{g}_{k}^{\pm}\right| \rightharpoonup G^{\pm}$ weakly in $L^{2}\left(\baseball{0}{3}\right)$. Once can then easily check, see for instance the proof of \cite[Proposition 4.3]{DS3}, that
$$
\left|D g^{\pm}\right| \leq G^{\pm}, \text{ a.e. in } \baseball{0}{3}.
$$

In particular, since $E_k\to 0$ and it bounds the size of the bad set, we have $\left|\baseball{0}{3} \backslash K_{k}\right| \rightarrow 0$, hence for every $s\in(2,3)$:

\begin{equation}\label{eq:lip-approx-by-harm:bounds-to-interpolate1}
\begin{aligned}
\operatorname{Dir}\left(g^{\pm}, \basehalfball{\pm}{0}{s}\right) & \leq \liminf _{k \rightarrow \infty} \int_{\basehalfball{\pm}{0}{s} \cap \Phi_{k}(K_{k})}\left(G^{\pm}\right)^{2} \leq \liminf _{k \rightarrow \infty} \operatorname{Dir}\left(\hat{g}_{k}^{\pm}, \basehalfball{\pm}{0}{s} \cap \Phi_{k}(K_{k})\right) \\
\overset{\eqref{eq:lip-approx-by-harm:relation-Dir-energy-error1}}&{\leq} \liminf _{k \rightarrow \infty} \operatorname{Dir}\left(g_{k}^{\pm}, \basehalfball{\pm}{0}{s} \cap K_{k}\right)= \liminf _{k \rightarrow \infty} \operatorname{Dir}\left(g_{k}^{\pm}, \basehalfball{\pm}{0}{s} \cap K_{k}\right).
\end{aligned}
\end{equation}

Let $\varepsilon>0$ be a small parameter to be chosen later, we apply \cite[Lemma 5.8]{DDHM} to $\left.\left(g^{+}, g^{-}\right)\right|_{\baseball{0}{3}}$ with such an $\varepsilon$ to produce a $(Q-\frac{Q^\star}{2})$-Lipschitz multivalued function $\left(g_{\varepsilon}^{+}, g_{\varepsilon}^{-}\right)$ satisfying:

\begin{equation}\label{Eq:LipApproximation:epsilon-estimates}
\int_{\basehalfball{\pm}{0}{3}} \mathcal{G}\left(g^{\pm}, g_{\varepsilon}^{\pm}\right)^{2}+\int_{\basehalfball{\pm}{0}{3}}\left(\left|D g^{\pm}\right|-\left|D g_{\varepsilon}^{\pm}\right|\right)^{2}+\int_{\basehalfball{\pm}{0}{3}}\left|D\left(\boldsymbol{\eta} \circ g^{\pm}\right)-D\left(\boldsymbol{\eta} \circ g_{\varepsilon}^{\pm}\right)\right|^{2} \leq \varepsilon ,
\end{equation}
\begin{equation*}
\int_{\partial \basehalfball{\pm}{0}{3}} \mathcal{G}\left(g^{\pm}, g_{\varepsilon}^{\pm}\right)^{2}+\int_{\partial \basehalfball{\pm}{0}{3}}\left(\left|D g^{\pm}\right|-\left|D g_{\varepsilon}^{\pm}\right|\right)^{2} \leq \varepsilon.
\end{equation*}

Additionally, we would like to interpolate without increasing too much Dirichlet energy in the transition region. To solve this problem, let us define the Radon measures
$$
\mu_{k}(A)=\int_{A \cap \basehalfball{+}{0}{3}}\left|D \hat{g}_{k}^{+}\right|^{2}+\int_{A \cap \basehalfball{-}{0}{3}}\left|D \hat{g}_{k}^{-}\right|^{2}, \quad A \subset \baseball{0}{3}.
$$
It is easy to check using \eqref{eq:lip-approx-by-harm:bound-dir-energy-f_k} that $\mu_{k}(A)\le C$ where $C$ is independent of $k$ and $A$. So, up to a subsequence, we can assume that $\mu_{k} \overset{*}{\rightharpoonup} \mu$ for some Radon measure $\mu$. We now choose $r \in(5 / 2,3)$ and a subsequence, not relabeled, such that
\begin{enumerate}[\upshape (A)]
    \item $\mu\left(\partial \baseball{0}{r}\right)=0,$
    
    \item $\mathbf{M}\left(\left\langle T_{k}-\left(\mathbf{G}_{f_{k}^{+}}+\mathbf{G}_{f_{k}^-}\right),|\mathbf{p}|, r\right\rangle\right) \leq C E_{k}^{1-2 \beta}$, where the map $|\mathbf{p}|$ is given by $\pi_{0} \times \pi_{0}^{\perp} \ni$ $(x, y)\mapsto|x|.$
\end{enumerate}

Indeed, by standard measure theory arguments, (A) is true for all but countably many radii while (B) can be obtained from the estimate \eqref{eq:thmlipaprox:T-graph(u)} through the slicing theory for currents. In particular, by (A) and the properties of weak convergence of measures, we have

$$
\begin{aligned}
\limsup _{s \rightarrow r} \limsup _{k \rightarrow \infty}& \biggl[ \int_{\basehalfball{+}{0}{r} \backslash \basehalfball{+}{0}{s}}\left|D \hat{g}_{k}^{+}\right|^{2}+\int_{A \cap \basehalfball{-}{0}{r} \backslash \basehalfball{-}{0}{s}}\left|D \hat{g}_{k}^{-}\right|^{2} \biggr]\\
& \leq \limsup _{s \rightarrow r} \mu\left(\overline{\baseball{0}{r}} \backslash \baseball{0}{s}\right)=0 .
\end{aligned}
$$

Hence, given $r \in(5 / 2,3)$ satisfying (A) and (B) above, we can now choose $s \in(5 / 2,3)$ such that

\begin{equation}\label{eq:lip-approx-by-harm:bounds-to-interpolate2}
   \limsup _{k \rightarrow \infty} \int_{\basehalfball{+}{0}{r} \backslash \basehalfball{+}{0}{s}}\left|D \hat{g}_{k}^{+}\right|^{2}+\int_{\basehalfball{-}{0}{r} \backslash \basehalfball{-}{0}{s}}\left|D \hat{g}_{k}^{-}\right|^{2} \leq \frac{c_{2}}{3} . 
\end{equation}

Finally, as aforementioned, we interpolate the pairs $\left(\hat{g}_{k}^{+}, \hat{g}_{k}^{-}\right)$ and $\left(g_{\varepsilon}^{+}, g_{\varepsilon}^{-}\right)$ which we can do now because all of them have the same interface $(\gamma, Q^\star\a{0})$ and we have control on their Dirichlet energy. We finally apply, for each $k$, the interpolation lemma to connect the functions $\left(\hat{g}_{k}^{+}, \hat{g}_{k}^{-}\right)$ and $\left(g_{\varepsilon}^{+}, g_{\varepsilon}^{-}\right)$ on the annulus $\baseball{0}{r} \backslash \baseball{0}{s}.$ This gives sets $\overline{\baseball{0}{s}} \subset V_{\lambda, \varepsilon}^{k} \subset W_{\lambda, \varepsilon}^{k} \subset \baseball{0}{r}$ and a $\left(Q-\frac{Q^{\star}}{2}\right)$ valued interpolation map $\left(\zeta_{k, \varepsilon}^{+}, \zeta_{k, \varepsilon}^-\right)$ with

$$
\begin{aligned}
\int_{\left(W_{\lambda, \varepsilon}^{k}\right)^\pm \setminus V_{\lambda, \varepsilon}^{k}}\left|D \zeta_{k, \varepsilon}^{\pm}\right|^{2} &\leq  C \lambda \int_{\left(W_{\lambda, \varepsilon}^{k}\right)^\pm \setminus V_{\lambda, \varepsilon}^{k}}\left(\left|D \hat{g}_{k}^{\pm}\right|^{2}+\left|D g_{\varepsilon}^{\pm}\right|^{2}\right)+\frac{C}{\lambda} \int_{\left(W_{\lambda, \varepsilon}^{k}\right)^\pm \setminus V_{\lambda, \varepsilon}^{k}} \mathcal{G}\left(\hat{g}_{k}^{\pm}, g_{\varepsilon}^{\pm}\right)^{2} \\
\leq C \lambda \int_{\left(W_{\lambda, \varepsilon}^{k}\right)^ \pm\setminus V_{\lambda, \varepsilon}^{k}}&\left(\left|D \hat{g}_{k}^{\pm}\right|^{2}+\left|D g_{\varepsilon}^{\pm}\right|^{2}\right)+\frac{C}{\lambda} \int_{\left(W_{\lambda, \varepsilon}^{k}\right)^\pm\setminus V_{\lambda, \varepsilon}^{k}} \left(\mathcal{G}\left(\hat{g}_{k}^{\pm}, g^{\pm}\right)^{2}+\mathcal{G}\left(g^{\pm}, g_{\varepsilon}^{\pm}\right)^{2}\right).
\end{aligned}
$$

Hence
$$
\underset{\lambda \rightarrow 0}{\limsup } \limsup _{\varepsilon \rightarrow 0} \limsup _{k \rightarrow \infty} \int_{\left(W_{\lambda,\varepsilon}^{k}\right)^\pm\setminus V_{\lambda,\varepsilon}^{k}}\left|D \zeta_{k, \varepsilon}^{\pm}\right|^{2}=0.
$$
Thus we can find $\lambda, \varepsilon>0$ sufficiently small such that
\begin{equation}\label{eq:lip-approx-by-harm:bounds-to-interpolate3}
  \limsup _{k \rightarrow \infty} \int_{\left(W_{\lambda, \varepsilon}^{k}\right)^\pm \backslash V_{\lambda, \varepsilon}^{k}}\left|D \zeta_{k, \varepsilon}^{\pm}\right|^{2}<\frac{c_{2}}{3}.  
\end{equation}
Moreover, up to further reduce $\varepsilon$, by \eqref{Eq:LipApproximation:epsilon-estimates} we can also assume that
\begin{equation}\label{eq:lip-approx-by-harm:bounds-to-interpolate4}
    \int_{\baseball{0}{r}^{\pm}}\left|D g_{\varepsilon}^{\pm}\right|^{2} \leq \int_{\baseball{0}{r}^{\pm}}\left|D g^{\pm}\right|^{2}+\frac{c_{2}}{6} .
\end{equation}

Now that we have interpolated the functions without adding too much energy, we define the $(Q-\frac{Q^\star}{2})$-Lipschitz function on $\baseball{0}{r}$ with interface $(\gamma, Q^\star\a{0})$ by

$$
\hat{h}_{k,\lambda,\varepsilon}^{\pm}:= \begin{cases}\hat{g}_{k}^{\pm} & \text {on } \baseball{0}{r} \backslash\left(W_{\lambda,\varepsilon}^{k}\right)^{\pm}, \\ \zeta_{k, \varepsilon}^{\pm} & \text {on }\left(W_{\lambda, \varepsilon}^{k}\right)^{\pm} \backslash V_{\lambda,\varepsilon}^{k} \\ g_{\varepsilon}^{\pm} & \text {on }\left(V_{\lambda, \varepsilon}^{k}\right)^{\pm}.
\end{cases},
$$

Let us then consider $\left(Q-\frac{Q^{\star}}{2}\right)$-valued map $\left(h_{k,\lambda,\varepsilon}^{+}, h_{k,\lambda,\varepsilon}^{-}\right)$ defined on $B_{k, 3}^{\pm}$ with interface $\left(\gamma_{k}, Q^\star\a{\psi_{k}}\right)$ given by 

$$h_{k,\lambda,\varepsilon}^{\pm}:=\left(\hat{h}_{k,\lambda,\varepsilon}^{\pm} \circ \Phi_{k}\right)\oplus\left(\kappa_{k} \circ \Phi_{k}\right),$$

which satisfies

\begin{equation}\label{eq:thm-lip-approx-by-harm:bounds-on-varying-sequences-energy}
\begin{aligned}
\liminf _{k \rightarrow \infty}\operatorname{Dir}\left(h_{k,\lambda,\varepsilon}^{+}, h_{k,\lambda,\varepsilon}^{-}, B_{k, r}\right) &= \liminf _{k \rightarrow \infty}\operatorname{Dir}\left(\hat{h}_{k,\lambda,\varepsilon}^{+}, \hat{h}_{k,\lambda,\varepsilon}^{-}, \baseball{0}{r}\right) \\
&\leq \operatorname{Dir}\left(g_{\varepsilon}^{+}, g_{\varepsilon}^{-}, \baseball{0}{r}\right) +\limsup _{k \rightarrow \infty}\operatorname{Dir}\left(\zeta_{k,\varepsilon}^{+}, \zeta_{k,\varepsilon}^{-},\left(W_{\lambda, \varepsilon}^{k}\right)\backslash V_{\lambda, \varepsilon}^{k}\right)\\
&\quad +\limsup _{k \rightarrow \infty}\operatorname{Dir}\left(\hat{g}_{k}^{+}, \hat{g}_{k}^{-}, \baseball{0}{r}\backslash \baseball{0}{s}\right) \\
\overset{\eqref{eq:lip-approx-by-harm:bounds-to-interpolate2},\eqref{eq:lip-approx-by-harm:bounds-to-interpolate3},\eqref{eq:lip-approx-by-harm:bounds-to-interpolate4}}&{<} \operatorname{Dir}\left(g^{+}, g^{-}, \baseball{0}{r}\right)+c_{2} \\
\overset{\eqref{eq:lip-approx-by-harm:bounds-to-interpolate1}}&{<}  \liminf _{k \rightarrow \infty}\left(\operatorname{Dir}\left(g_{k}^{+}, g_{k}^{-}, \baseball{0}{r}\cap K_{k}\right)\right)+2c_{2}.
\end{aligned}
\end{equation}

Let us consider the function $w_{k,\lambda,\varepsilon}^{\pm}(x):=E_{k}^{1 / 2} h_{k,\lambda,\varepsilon}^{\pm}(x).$ 
Observe that, by the constructions of $\hat{g}_k^\pm$, $w_{k,\lambda,\varepsilon}^{\pm}|_{\partial \baseball{0}{r}}=f_{k}^{\pm}|_{\partial \baseball{0}{r}}$ and $\operatorname{Lip}\left(w_{k,\lambda,\varepsilon}^{\pm}\right) \leq C E_{k}^{\beta}$. We are now ready to construct a sequence of competitors one for each $T_k$ which for large $k$ will contradict the almost minimality of the sequence $T_{k}$. First of all, by the isoperimetric inequality, \cite[Section 4.4.2]{Fed}, there is a current $S_{k}$ such that
\begin{equation}\label{eq:first-lip-approx:bound-competitor}
    \partial S_{k,r}:=\left\langle T_{k}-\left(\mathbf{G}_{f_{k}^{+}}+\mathbf{G}_{f_{k}^{-}}\right),|\mathbf{p}|, r\right\rangle \quad \text { and } \quad \mathbf{M}\left(S_{k,r}\right) \leq C\left(E_{k}^{1-2 \beta}\right)^{\frac{m}{m-1}} \overset{(\beta<\frac{1}{2m})}{=} o\left(E_{k}\right).
\end{equation}

Let $Z_{k,r}:= \bG_{w_k^+} \res \cyl{0}{r} + \bG_{w_k^-} \res \cyl{0}{r} + S_{k,r}$. Since $w_{k,\lambda,\varepsilon}^{\pm}|_{\partial \baseball{0}{r}}=f_{k}^{\pm}|_{\partial \baseball{0}{r}}$, we can see that the boundary of $Z_k$ matches that of $T_k \res \cyl{0}{r}$, thus it is an admissible competitor and will furnishes the desired contradiction. To that end, first we compare the Dirichlet energies of $w_{k,\lambda,\varepsilon}^+$ and $f_k^+$. To begin with this comparison, we note that, up to a subsequence not relabeled, it holds

\begin{equation}\label{e:bus}
    \begin{aligned}
    \D (w^+_{k,\lambda,\varepsilon}, w^-_{k,\lambda,\varepsilon}, B_{k,r}) &= E_k\D (h^+_{k,\lambda,\varepsilon}, h^-_{k,\lambda,\varepsilon}, B_{k,r}) \\
    \overset{\eqref{eq:thm-lip-approx-by-harm:bounds-on-varying-sequences-energy}}&{<} E_k\D (g^+_k, g^-_k, \baseball{0}{r}\cap K_k) +c_2E_k\\
    &= \D (f^+_k, f^-_k, B_{k,r}\cap K_k) + 2c_2E_k,
    \end{aligned}
\end{equation}
for $k$ large enough.
In addition, the latter inequality combined with the second inequality in \eqref{eq:lip-approx-by-harm:eq-that-gives-contradiction} implies for $k$ large enough that

\begin{equation}\label{e:improv3}
\D (w_{k,\lambda,\varepsilon}^+, w_{k,\lambda,\varepsilon}^-, B_{k,r}) < \be_{T_k} (\baseball{0}{r}) -  c_2 E_k + o(E_k).
\end{equation}

Finally, we estimate

\begin{equation}\label{eq:thm:first-lip-approx:better-competitor}
\begin{aligned}
\bM(Z_k) - \bM(T_k) &\leq \bM(\bG_{w_k^+} \res \cyl{0}{r}) + \bM(\bG_{w_k^-} \res \cyl{0}{r}) + \bM(S_k) - \bM(T_k) \\
\overset{\textup{Taylor}}&{\leq}  Q|B_{k,r}^+| + (Q-Q^\star)|B_{k,r}^-| + \D (w_{k,\lambda,\varepsilon}^+, w_{k,\lambda,\varepsilon}^-, B_{k,r}) + o(E_k) - \bM(T_k) \\
&\leq Q|B_{k,r}^+| + (Q-Q^\star)|B_{k,r}^-| + \D (w_{k,\lambda,\varepsilon}^+, w_{k,\lambda,\varepsilon}^-, B_{k,r}) + o(E_k) \\
&\quad -Q|B_{k,r}^+| - (Q-Q^\star)|B_{k,r}^-| - \be_{T_k}(\baseball{0}{r}) \\
\overset{\eqref{e:improv3}}&{<} -c_2E_k +o(E_k),
\end{aligned}    
\end{equation}

the expression is negative when $k$ is large enough. In particular $Z_k$ is a competitor with less mass than $T_k$ and this completes the proof of \eqref{eq:thm-first-lip-approx-by-harm-sheet:excees-measure-of-bad-set}, we recall that \eqref{eq:thm-first-lip-approx-by-harm-sheet:Dir-energy-on-bad-set} follows from \eqref{eq:thm-first-lip-approx-by-harm-sheet:excees-measure-of-bad-set} as mentioned at the beginning of the proof.

\textbf{Proof of {\eqref{eq:thm-harm-aprox:estimate-u^+}}, {\eqref{eq:thm-harm-aprox:estimate-u^-}} and {\eqref{eq:thm-harm-aprox:estimate-center-of-u}}.}  As in the first part, we argue by contradiction assuming \eqref{eq:first-lip-approx:contradicition1}, \eqref{eq:first-lip-approx:contradicition2}, and \eqref{eq:first-lip-approx:contradicition3} becomes 

\begin{enumerate}[\upshape (iii)'] 
    \item The $E_k^\beta$-Lipschitz approximations $(f_k^+, f_k^-)$ fail to satisfy one among the estimates \eqref{eq:thm-harm-aprox:estimate-u^+},
    \eqref{eq:thm-harm-aprox:estimate-u^-} and \eqref{eq:thm-harm-aprox:estimate-center-of-u} for any choice of the function $h$.
\end{enumerate}

We use the same notations of the previous step for $\psi_k, g_k^{\pm}, \Phi_k, \kappa_k, \hat{g}_k^{\pm}$ and $g^{\pm}$. Therefore, we now claim that

\begin{enumerate}[\upshape (A)]
\item The convergence of $\hat g_k^\pm$ to $g^\pm$ is strong in $W^{1,2} (\baseball{0}{5/2})$, namely 
$$
\lim_{k\to \infty} \D (\hat g_k^+, \hat g_k^-, \baseball{0}{5/2}) = \D (g^+, g^-, \baseball{0}{5/2}) \, .
$$

\item $(g^+, g^-)$ is a $(Q-\frac{Q^{\star}}{2})$-minimizer in $\baseball{0}{5/2}$. 
\end{enumerate}

Recall that, by Theorem \ref{thm:first-lip-approx-and-good-set} and the construction, $(g^+, g^-)$ collapses at the interface $(\gamma,Q^\star\a{0})$, thus provided we assume that (A) and (B) are proved, from  Theorem \ref{thm:harmonic-regularity-(Q-Q^*/2)-minimizing} we would then infer the existence of a classical harmonic function $\hat{h}$ on $\baseball{0}{\sfrac{5}{2}}$ which vanishes identically
on $\{x_m =0\}$ such that $g^+ = Q \a{\hat{h}}$ and $g^- = (Q-Q^{\star}) \a{\hat{h}}$. If we set $h_k = E_k^{\sfrac{1}{2}}\hat{h}$, the following hold
\begin{align*}
\int_{B_{k,5/2 }^+}\mathcal{G}(f_k^+,Q\a{h_k})^2+\int_{B_{k,5/2}^+}\left(|Df_k^+|-\sqrt Q|D h_k|\right)^2 &= \; o (E_k) ,\\
\int_{B_{k,5/2 }^-}\mathcal{G}(f_k^-,(Q-Q^{\star})\a{h_k})^2+\int_{B_{k,5/2}^-}\left(|Df_k^-|-\sqrt {(Q-Q^{\star})}|D h_k|\right)^2 &= o(E_k), \\
\int_{B_{k,5/2}^\pm} \abs{D(\boldsymbol{\eta}\circ f_k^\pm) - D h_k}^2 &= o(E_k).
\end{align*}

But these estimates are incompatible with (iii)' above.  Hence,  at least one between (A) and (B) needs to fail. As in the previous section we will use this to contradict the minimality of $T_k$. Note that in both cases there exists a  $(Q-\frac{Q^{\star}}{2})$-valued function  $(\bar g^+, \bar g^-)$ with  interface \((\gamma, Q^\star\a{0})\), \(\gamma=\{x_m=0\}\), and a positive constant \(c_3>0\), such that 

\begin{equation}\label{eq:lip-approx-by-harm:bound-gbar-energy}
\D(\bar g^+,\bar g^-, \baseball{0}{s}) \le \liminf_{k\to \infty} \D(\hat g^+,\hat g^-, \baseball{0}{s})-2c_3
\end{equation}

for all \(s\in (5/2,3)\). Indeed this is  true with $(\bar g^+, \bar g^-)=(g^+,  g^-)$ if (A) fails, while if (B) fails we choose $(\bar g^+, \bar g^-)$ to be a $(Q-\frac{Q^{\star}}{2})$-minimizer with boundary data \(g^\pm\) on \(\partial \baseball{0}{5/2}\) extended to be equal to \(g^\pm\) on  \(\baseball{0}{3}\setminus \baseball{0}{5/2}\). We can now follow the exactly the same argument as in the previous step to find a radius \(r\in (5/2,3)\) and functions  \(\hat h_k^\pm\) such that
\[
\bM (\langle T_k - (\bG_{f_k^+}+ \bG_{f_k^-}), |\mathbf{p}|,r \rangle) \leq C E_k^{1-2\beta}
\]
and, arguing as we have done for \eqref{eq:thm-lip-approx-by-harm:bounds-on-varying-sequences-energy},

\begin{equation}
\begin{aligned}
\liminf_{k\to \infty} \D(h^+,h^-,B_{k,r})&\le \D(\bar{g}^+, \bar{g}^-,\baseball{0}{r})+c_3
\\
&\le \liminf_{k\to \infty} \D(g^+, g^-,B_{k,r})-c_3.\label{e:ddai}
\end{aligned}
\end{equation}

Defining $w_k^\pm$ as above, we again observe that
$w_k^\pm|_{\partial \basehalfball{\pm}{0}{r}} = f_k^\pm|_{\partial \basehalfball{\pm}{0}{r}}$. We then construct the same competitor currents to test the minimality of $T_k$. First we consider a current \(S_k\) supported in \(\R^{m+n}\) such that 

\begin{equation}\label{eq:first-lip-approx:competitor-to-harmonotonicity}
\partial S_k = \langle T_k - (\bG_{f_k^+}+ \bG_{f_k^-}), |\mathbf{p}|, r\rangle \mbox{ and }
\bM (S_k) \leq C (E_k^{1-2\beta})^{\frac{m}{m-1}} = o (E_k)\, ,
\end{equation}

where we again used $\beta<\frac{1}{4m}.$
Then we define, as before,  $Z_k := \bG_{w_k^+} \res \cyl{0}{r} + \bG_{w_k^-} \res \cyl{0}{r} + S_k$, for which the minimality condition guarantees 

$$\bM (Z_k) \geq \bM (T_k \res \cyl{0}{r})\, .
$$

Since we proved the first part of the theorem, we use it to show that 

$$\be_{T_k}(\baseball{0}{r})= \D (f_k^+, \basehalfball{+}{0}{r}) + \D (f_k^-, \basehalfball{-}{0}{r}) +  O(\eta_k E_k)\, .
$$
Observe that now we can choose \(\eta_k\to 0\) as \(k\to \infty\). Arguing as in \eqref{e:bus} and relying on \eqref{e:ddai} we also have 

$$\D (w_k^+, w_k^-, \baseball{0}{r}) \leq  \D (f_k^+, f_k^-, \baseball{0}{r}) - c_3 E_k + o (E_k)\, .
$$
Accordingly, the latter inequality combined with \eqref{eq:lip-approx-by-harm:eq-that-gives-contradiction} implies

$$\D (w_k^+, w_k^-, B_{k,r}) < \be_{T_k} (\baseball{0}{r}) -  c_3 E_k .
$$
As before, see \eqref{eq:thm:first-lip-approx:better-competitor}, we complete the proof.
\end{proof}

\section{Superlinear decays and Lusin type strong Lipschitz approximations}

The approximation furnished in the last subsection, Theorem \ref{thm:first-lip-approx-and-good-set}, have sublinear exponents bounding the size of the bad set among other important quantities. In the construction of the center manifold, in order to derive good properties of it we will need superlinear decays in some of the estimates of Theorem \ref{thm:first-lip-approx-and-good-set}. To that end, we need to improve our approximation. In fact, we need an accurate height bound and harmonic approximations to achieve a satisfactory excess decay and therefore provide stronger Lipschitz approximations that will appear in Theorem \ref{thm:lip-approx-superlinear-decay}, c.f. \cite[Chapter 6]{DDHM}, with the subtle exponents estimating the error of such approximations. So, as a first step we state the height bound of the current $T$.

\begin{lemma}[Height bound, Lemma 10.4, \cite{DNS}]\label{lemma:height-bound}
Let $T$, $\cyl{p}{4r}$, $\Gamma$ and $\pi_0:= \R^m \times \{0\}$ be as in Assumption \ref{assump:general.asump-first-lip} with $C_0 = 0$. Then, there exist positive constants $\eps_h=\eps_h(Q, Q^\star, m, n)$ and $C_h=C_h (Q, Q^\star, m, n)$ such that, if $ \bE (T, \cyl{p}{4r}) +\bA \leq \varepsilon_h$, then 
$$\bh (T, \cyl{p}{2r}, \pi_0) \leq C_h (r^{-1}\bE (T, \cyl{p}{4r})+ \bA)^{\frac{1}{2}} r^\frac{3}{2} .$$
\end{lemma}

\subsection{Improved excess estimate}

We will follow a very well known process that using the height bound provided by Lemma \ref{lemma:height-bound} allows to obtain our proof of the improved excess decay. To achieve this result we will firstly prove a milder decay in Lemma \ref{lemma:milder-excess-decay} and then iterate it to reach the needed superlinear excess estimate. To prove this milder statement which will be stated for the modified excess function introduced in Definition \ref{def:excess-flat}, we will reduce this milder decay of the current $T$ in some steps until we can rely on a similar decay for harmonic functions, this reduction will be possible thanks to Theorem \ref{thm:first-lip-approx-by-harmonic-sheet}. 

\begin{lemma}[Milder excess decay]\label{lemma:milder-excess-decay}
Let $T$ and $\Gamma$ be as in Assumptions \ref{assumptions} with $C_0 = 0$, $p\in \Gamma\cap U$ is a two-sided collapsed point where $U$ is the neighborhood in Definition \ref{defi:collapsed}. Then, for every $q\in U\cap\Gamma$ and $\varepsilon>0$, there is an $\eps_0=\varepsilon_0 (\varepsilon, Q, Q^\star,m,n)>0$ $($assume that $\varepsilon_h\leq \varepsilon_0^2)$ and a $M_0=M_0 (\varepsilon, Q, Q^\star,m,n)>0$ with the following property. We set $\theta(\sigma) := \max\{\bE^\flat (T, \ball{q}{\sigma}) ,M_0 \bA^2 \sigma^2\}$, and assume 

\begin{equation}\label{eq:lemma-milder-decay:small-excess-assump}
\bA^2 \sigma^2 + E := \|A_\Gamma\|^2 \sigma^2 + \bE^\flat (T, \ball{q}{4 \sigma}) < \varepsilon_0,
\end{equation}
\begin{equation}\label{eq:lemma:milder-decay:densities-p-q-close}
 \|T\|(\ball{q}{4\sigma}) \leq \left( \Theta^m(T,p)+\frac{1}{4}\right)\omega_m(4\sigma)^m.
\end{equation}

Then we have

\begin{equation}\label{eq:thm:milder-excess-decay}
\theta(\sigma) \leq \max \{ 2^{-4+4\varepsilon} \theta (4\sigma),
 2^{-2+2\varepsilon}  \theta(2\sigma)\}\, .
\end{equation}
\end{lemma}

This milder statement is not enough for our purposes, since the excess are considered in balls with the same center $q$. However, it facilitates a lot the proof of the improved excess decay which we enunciate below.

\begin{theorem}[Improved excess decay and height bound]\label{thm:improved-excess-decay}
Let $T$ and $\Gamma$ be as in Assumption \ref{assumptions} with $C_0=0$. If $p\in \Gamma \cap U$ is a two-sided collapsed point with density $\Theta (T, p) = Q - \frac{{Q^\star}}{2}$, $U$ is a neighborhood of Definition \ref{defi:collapsed}, then there exists \(r>0\) such that $\ball{p}{r}\subset U$, for all $q\in\ball{p}{r}\cap U$ there exists a $m$-dimensional plane $\pi(q)$ which $T_q\Gamma\subset\pi(q)$, and for all $\varepsilon >0$ there is a constant $C = C (m,n,Q^\star, Q, \varepsilon)>0$ with

\begin{equation}\label{eq:decay-improved-excess}
    \bE^\flat (T, \ball{q}{\rho}) \leq \bE^\flat (T, \ball{q}{\rho},\pi(q)) \leq C \left(\frac{\rho}{r}\right)^{2-2\varepsilon}  \bE^\flat (T, \ball{p}{2r}) + C \rho^{2-2\varepsilon}r^{2\varepsilon} \bA^2,
\end{equation}

for all $\rho\in (0, r)$. Moreover, if we take $\rho\in (0, {\textstyle{\frac{r}{2\sqrt{2}}}})$, then

\begin{equation}\label{eq:thm-improved-excess-decay:height-bound}
\bh (T, \ball{q}{\rho}, \pi (q)) \leq C (r^{-1} \bE^\flat (T, \ball{p}{2r}) + \bA)^{\frac{1}{2}} \rho^{\frac{3}{2}}, \ \forall q\in \Gamma \cap \ball{p}{r}.
\end{equation}
\end{theorem}
\begin{remark}
{We announce that, in Theorem \ref{thm:uniqueness-dimension-m}, we prove that $\pi(q)$ is in fact the support of the unique tangent cone to $T$ at $q$.}
\end{remark}

We begin with the proof of the Lemma \ref{lemma:milder-excess-decay} which will be used to prove the Theorem \ref{thm:improved-excess-decay}.

\begin{proof}[Proof of Lemma \ref{lemma:milder-excess-decay}]\label{subsec:exc_dec_prel}

Without loss of generality by scaling, translating and rotating, we can assume $\sigma=1$, $q=0$, $\bE^\flat(T,\ball{0}{2})=\bE (T,\ball{0}{2},\pi_0)$, where $\pi_0=\R^m\times\{0\}$, and $T_0\Gamma = \R^{m-1}\times \{0\}$. We begin assuming

\begin{equation}\label{e:Apiccolo}
\bE^\flat(T,\ball{0}{2})\ge 2^{-m} M_0 \bA^2 \quad \text{and}\quad \bE^\flat(T,\ball{0}{2})\ge 2^{-4-m}  \bE^\flat(T,\ball{0}{4}).
\end{equation}

Indeed, note that 

\[
\theta(1) = \max\{ M_0\bA^2, \bE^\flat (T, \ball{0}{1})\} \leq \max \{M_0 \bA^2, 2^{m}\bE^\flat(T,\ball{0}{2})\} \, .
\] 

So, if the first inequality in \eqref{e:Apiccolo} fails, by the latter inequality, we have

\begin{eqnarray*}
\theta(1) \leq M_0 \bA^2= 2^{-2} (2^2 M_0 \bA^2) \leq 2^{-2}\theta(2)\, ,
\end{eqnarray*}

whereas, if the second inequality in \eqref{e:Apiccolo} fails, then

\begin{eqnarray*}
\theta(1)\leq \max \{M_0 \bA^2, 2^{-4} \bE^\flat (T, \ball{0}{4})\} \leq 2^{-4} \theta(4)\, .
\end{eqnarray*}

Hence in both cases the conclusion should hold. Reiterating, under assumption \eqref{e:Apiccolo}, we need to show the decay estimate:

\begin{equation}\label{e:gollonzo}
\bE^\flat(T,\ball{0}{1})\le 2^{2\varepsilon-2}\bE^\flat(T,\ball{0}{2})\,.
\end{equation}

Let us now fix a positive $\eta<1$, to be chosen sufficiently small later, and consider the cylinder ${U_\eta}:=\mathrm{B}_{4-\eta} (0, \pi_0) + \mathrm{B}_{\sqrt\eta} (0, \pi_0^\perp)$, which by abuse of notation we denote by $B_{4-\eta}\times B^n_{\sqrt{\eta}}$. If $\varepsilon_0>0$ is sufficiently small, we claim that 
\begin{eqnarray}
\spt (T)\cap\partial {U_\eta} \subset \partial B_{4-\eta}\times B^n_{\sqrt\eta} \label{e:supp1}\\
\ball{0}{4-\eta}\cap\spt(T) \subset  {U_\eta}\,.\label{e:supp2}
\end{eqnarray} 
Otherwise, arguing by contradiction, we would have a sequence of currents $T_k$ satisfying the assumptions of the theorem with $\varepsilon_0=\frac 1k$, but violating either \eqref{e:supp1} or \eqref{e:supp2}.  Then $T_k$ would converge, in the sense of currents, to a current $T_\infty$ that is area minimizing whose excess w.r.t. $\pi_0$ is identically zero so its support is contained in the plane $\pi_0$, $\partial T_\infty=Q^\star\a{T_q\Gamma}$. Thus we are in position to apply the Constancy Lemma, \cite[4.1.17]{Fed}, to asserts 
\[
T_\infty:=Q' \a{B^+_4}+(Q'-Q^\star)\a{B^-_4}\, ,
\]

where $B_4^\pm=B_4 (0, \pi_0) \cap\{\pm x_m>0\}$ and $Q'\ge Q^\star$ is a positive integer. Since $\partial T_k = Q^\star\a{\R^{m-1}\times\{0\}}$, $\forall k\in\N$, we can use the area-minimizing property to obtain an uniform bound on $\|T_k\|(B_4)$ and thus be in place to apply \cite[Theorem 7.2, Chapter 6]{simon2014introduction} which says that, since $T_k$ converge to $T_\infty$ in the sense of currents, the supports of $T_k$ converge to either $\overline B_4$ in case $Q'>Q^\star$ or $\overline B_4^+$ otherwise, i.e. if $Q'=Q^\star$, in the Hausdorff sense in every compact subset of $\ball{0}{4}$. This is a contradiction because the following both inequalities hold

$$ \operatorname{dist}_H(\overline{\ball{0}{4-\eta}\setminus {U_\eta}}, \ball{0}{4}) > 0 \ \text{and} \ \operatorname{dist}_H(\partial {U_\eta}\setminus(\partial B_{4-\eta}\times B^n_{\sqrt \eta}), \ball{0}{4}) > 0 .$$

We have therefore proved \eqref{e:supp1} and \eqref{e:supp2}. 

We now let $T_0$ be a tangent cone of $T$ at $0$ and $\rho_k\to0^+$ a sequence such that $T_{0,\rho_k}\to T_0$. By a standard argument using the Constancy Lemma, we know that 
\begin{equation}\label{Eq:LemmaMilderExcessConstancyLemma}
    {\bp_{\pi_0}}_\sharp T_0 =Q'\a{\pi_0}+(Q'-Q^\star)\a{\pi_0},
\end{equation} for some natural number $Q'$. By the lower semicontinuity of the total variation and \eqref{eq:lemma:milder-decay:densities-p-q-close}, we further notice that we necessarily have
$\|{\bp_{\pi_0}}_\sharp T_0\| (\ball{0}{4}) \leq (Q-\frac{2Q^\star - 1}{4}) \omega_m 4^m$. Hence, by the monotonicity formula 
\begin{equation}\label{eq:lemma:milder-decay:T_infty-from-above}
\Theta^m({\bp_{\pi_0}}_\sharp T_0, 0) \leq Q-\frac{2Q^\star - 1}{4}.
\end{equation}
On the other hand, the upper semicontinuity of the density w.r.t. the convergence of area-minimizing currents \cite[Chapter 7, Section 3, Eq. (12)]{simon2014introduction} and the fact that $p$ is a two-sided collapsed point allow us to conclude 
\begin{equation}\label{eq:lemma:milder-decay:T_infty-from-below}
\Theta^m({\bp_{\pi_0}}_\sharp T_0, 0) \geq \limsup \Theta^m({\bp_{\pi_0}}_\sharp T_{0,\rho_k}, 0) = Q-\frac{Q^\star}{2}.
\end{equation}
By equations \eqref{eq:lemma:milder-decay:T_infty-from-above}, \eqref{eq:lemma:milder-decay:T_infty-from-below}, and \eqref{Eq:LemmaMilderExcessConstancyLemma}, we have 
$$Q-\frac{2Q^\star - 1}{4}\ge\Theta^m({\bp_{\pi_0}}_\sharp T_0, 0)=Q'-\frac{Q^\star}2\ge Q-\frac{Q^\star}{2},$$
since $Q'$ is an integer it turns out that $Q'=Q$. 
By \eqref{e:supp2}, we can straightforwardly check that 
$$ \operatorname{dist}_H(\spt(T),\spt(T_0))) < \eta, $$
which, provided $\eta$ and $\varepsilon_0$ are small enough, leads to 
$$ \|T\|(\ball{0}{r}) = \|T_0\|(\ball{0}{r}) + O(\eta^{m-1}), \ \ \forall r\in(1, 4-\frac{\eta}{2}). $$
Thus, we can state the following property:

\begin{itemize}
\item[(A)] the mass of $T$ in the ball $\ball{0}{r}$ is $\left(Q-{\frac{Q^\star}{2}}\right) \omega_mr^m + O(\eta^{m-1})$, for any radius $1\le r\le 4-\frac\eta 2$.
\end{itemize}

Next, let us define ${S_\eta}:=T\res {U_\eta}$. Observe that \eqref{e:supp1} and \eqref{e:supp2} imply:

\begin{itemize}
\item[(B)] $\partial {S_\eta}\res\cyl{0}{4-\eta}=Q^\star\a{\Gamma\cap\cyl{0}{4-\eta}}$;
\item[(C)] $T\res\ball{0}{4-\eta}={S_\eta}\res\ball{0}{4-\eta}$.
\end{itemize}

Choose a plane $\overline\pi$ which minimizes the boundary excess, i.e., which contains $T_0\Gamma$ and $\bE(T,\ball{0}{4},\overline\pi)=\bE^\flat(T,\ball{0}{4})$. Let us observe that, since \(\pi_0\) is the optimal plane for \(\bE^\flat(T,\ball{0}{2})\), we have

\begin{equation}\label{e:nonsisa0}
\begin{aligned}
|\overline\pi-\pi_0|^2\|T\|(\ball{0}{2}) &=\int_{\ball{0}{2}}|\overline\pi-\pi_0|^2\,d\|T\|\\
&\le 2\int_{\ball{0}{2}}|\vec{T}-\pi_0|^2\,d\|T\|+2\int_{\ball{0}{2}}|\vec T-\overline\pi|^2\,d\|T\|\\
&\le 2\cdot 2^m \omega_m\bE^\flat(T,\ball{0}{2})+2\cdot 4^m\omega_m\bE^\flat(T,\ball{0}{4})\\
&\le C\bE^\flat(T,\ball{0}{4}).
\end{aligned}  
\end{equation}

Moreover,

\begin{equation}\label{e:nonsisa}
\begin{aligned}
\bE ({S_\eta},\cyl{0}{4-\eta}) & \le\bE(T,\ball{0}{4-\sfrac\eta 2},\pi_0)\\
\overset{\textup{Triangular}}&{\le} \frac{2}{\omega_m(4-\sfrac{\eta}{2})^m}\left(\bE^\flat(T,\ball{0}{4}) + |\overline\pi-\pi_0|^2\|T\|(\ball{0}{4-\sfrac\eta 2})\right)\\
\overset{\text{(A)}}&{\le} 2\bE^\flat(T,\ball{0}{4}) + C|\overline\pi-\pi_0|^2\|T\|(\ball{0}{2})
\overset{\eqref{e:nonsisa0}}{\le} C\bE^\flat(T,\ball{0}{4}),
\end{aligned}  
\end{equation}

Moreover, recalling that ${\bp}:\R^{m+n}\to\pi_0$ is the orthogonal projection, by the Constancy Lemma (\cite[4.1.17]{Fed}), we have

\begin{itemize}
\item[(D)] ${\bp}_\sharp {S_\eta}={Q_\bp}\a{\Omega^+}+({Q_\bp}-Q^\star)\a{\Omega^-}$, where ${Q_\bp}$ is a positive integer and $\Omega^\pm$ are the regions in which $\mathrm{B}_4(0, \pi_0)$ is divided by ${\bp} (\Gamma)$; in particular  
\[
\partial\a{\Omega^+}\res\cyl{0}{4-\eta}=-\partial\a{\Omega^-}\res\cyl{0}{4-\eta}={\bp}_\sharp\a{\Gamma}\res\cyl{0}{4-\eta}.
\]
\end{itemize}

Since ${S_\eta} = T \res {U_\eta}$ and ${U_\eta}\subset \ball{0}{4-\sfrac{\eta}{2}}$, clearly 

\begin{equation}\label{eq:bound-restricted-current-S-0}
    \|{S_\eta}\| (\cyl{0}{4-\eta}) \leq \|T\| (\ball{0}{4-\sfrac{\eta}{2}}).
\end{equation} 

Since projections do not increase mass, we obtain

\begin{equation}\label{eq:bound-restricted-current-S}
    \|{S_\eta}\| (\cyl{0}{4-\eta}) \geq \|{\bp}_\sharp{S_\eta}\| (\cyl{0}{4-\eta}).
\end{equation}

Assuming that the constant $\varepsilon_0$ in the assumption of the theorem is sufficiently small, we conclude that ${\bp}_\sharp\a{\Gamma}\res\cyl{0}{4-\eta}$ is close to $T_0\Gamma = \R^{m-1}\times \{0\}$. In particular, $|\Omega^\pm|$ is close to $|\basehalfball{\pm}{0}{4-\eta}|$ and thus 
${Q_\bp} |\Omega^+| + ({Q_\bp}-Q^\star) |\Omega^-|$ is close to $({Q_\bp}-\frac{Q^\star}{2}) \omega_m (4-\eta)^m$ too. Therefore, if $\varepsilon_0$ is smaller than a geometric constant, we infer from \eqref{eq:bound-restricted-current-S} that

$$\|{S_\eta}\|  (\cyl{0}{4-\eta}) \geq ({Q_\bp} - \frac{2Q^\star + 1}{4}) \omega_m (4-\eta)^m .$$

In addition, by (A), a sufficiently small $\varepsilon_0$ imply

$$
\begin{aligned}
({Q_\bp} - \frac{2Q^\star + 1}{4}) \omega_m (4-\eta)^m \overset{\eqref{eq:bound-restricted-current-S}}&{\leq} \|{S_\eta}\|  (\cyl{0}{4-\eta}) 
\overset{\eqref{eq:bound-restricted-current-S-0}}{\leq} \|T\| (\ball{0}{4-\eta/2}) \\
\overset{\text{(A)}}&{\leq} (Q-\frac{2Q^\star-1}{4}) \omega_m (4-\frac{\eta}{2})^m,
\end{aligned} $$

we achieve that ${Q_\bp} \leq Q$ provided $\eta$ is chosen smaller than a geometric constant. On the other hand, 

\[
\|{S_\eta}\| (\cyl{0}{4-\eta}) \leq {Q_\bp} |\Omega^+| + ({Q_\bp}-Q^\star) |\Omega^-| + \bE ({S_\eta},\cyl{0}{4-\eta})\, .
\] 

Using \eqref{e:nonsisa} and the argument above, if $\varepsilon_0$ is sufficiently small we get $\|{S_\eta}\| (\cyl{0}{4-\eta}) \leq ({Q_\bp}-\frac{2Q^\star-1}{4}) \omega_m (4-\eta)^m$. 
Recall that (C) ensures that $\|T\|(\ball{0}{4-\eta})\le\|{S_\eta}\|(\cyl{0}{4-\eta})$, and, using (A), we also have $\|T\| (\ball{0}{4-\eta}) \geq (Q-\frac{2Q^\star+1}{4}) (4-\eta)^m$. Thus necessarily ${Q_\bp}\geq Q$, consequently, we have $$Q_\bp = Q.$$
Next, since $T\res\ball{0}{2}={S_\eta}\res\ball{0}{2}$, then
$$
\begin{aligned}
\bA^2\overset{\eqref{e:Apiccolo}}&{\le}  2^{m}M_0^{-1}\bE^\flat(T,\ball{0}{2}) \le 2^{m}\left(\frac{4-\eta}2\right)^mM_0^{-1}\bE({S_\eta},\cyl{0}{4-\eta}) \\
\overset{\eqref{e:nonsisa}}&{\le} C M_0^{-1}\bE^\flat(T,\ball{0}{4})\, .
\end{aligned}
$$

By the last inequality and $Q_\bp = Q$, we finally proved that we are in position to apply Theorem \ref{thm:first-lip-approx-by-harmonic-sheet} with $\beta=\frac{1}{5m}$ and a sufficiently small parameter $\eta_*$ to be chosen later, provided $\varepsilon_0$ is sufficiently small and $M_0$ is sufficiently large. 

\textbf{Reduction to excess decay for graphs}

From now on we let $(u^+,u^-)$ and $h$ be as in Theorem \ref{thm:first-lip-approx-by-harmonic-sheet}. In particular, recall that $(u^+,u^-)$ is the $E^\beta$-approximation of Theorem \ref{thm:first-lip-approx-and-good-set} and $h$ is a single-valued harmonic function. Moreover, denote by $E$ the cylindrical excess $\bE({S_\eta},\cyl{0}{4-\eta})$ and record the estimates:

\begin{equation}\label{e:bibo}
\bA^2  \le C_0M_0^{-1}E\quad \text{and}\quad
E  \le C_0\bE^\flat(T,\ball{0}{2}),
\end{equation}

where $C_0$ is a geometric constant and the second inequality follows by combining \eqref{e:nonsisa} and \eqref{e:Apiccolo}. Next, define $\pi$ to be the plane given by the graph of the linear function $x\mapsto (Dh(0)x,0)$. Since, by the Schwarz reflection principle and the unique continuation for harmonic functions, we obtain that $h$ is odd, and \(h(x',0)=0\), so, we have  that 

\[
\pi\supset T_{0}\Gamma=\R^{m-1}\times\{0\}.
\]

Moreover, by elliptic estimates, 

\begin{equation}\label{e:grst}
 |\pi|\le C|Dh(0)|\le ( C\D(h,\baseball{0}{\frac 52(4-\eta)}))^\frac 12\overset{\textup{Thm. }\ref{thm:first-lip-approx-by-harmonic-sheet}}{\le} CE^\frac 12.
\end{equation}

Fix $\overline\eta$ to be chosen later. The following inequality is a consequence of the reduction argument given in \cite[Theorem 6.8]{DDHM} where the authors reduce the whole discussion to the analysis of a decay for classical harmonic functions using Theorem \ref{thm:first-lip-approx-by-harmonic-sheet}

\begin{equation}\label{e:gollonzo2}
\bE(\bG_{u^+}+\bG_{u^-},\cyl{0}{1},\pi)\le(2-\overline\eta)^{ - (2-\varepsilon)}\bE(\bG_{u^+}+\bG_{u^-},\cyl{0}{2-\overline\eta})+C\overline\eta E\,.
\end{equation}

Now, we claim that this inequality allows us to conclude \eqref{e:gollonzo}. First of all, by the Taylor expansion of the mass of a Lipschitz graph, \cite[Corollary 3.3]{DS2}, and the bound on Dirichlet energy of $u^\pm$ on the bad set, we conclude
\begin{equation}\label{eq:excess-decay:reduction to graphs}
\begin{aligned}
\bE(\bG_{u^+}+\bG_{u^-},\cyl{0}{2-\overline\eta}) &\leq \bE({S_\eta},\cyl{0}{2-\overline\eta})+ C\int_{\Omega^+\setminus K}|Du^+|^2+C\int_{\Omega^-\setminus K}|Du^-|^2\\
\overset{\eqref{eq:thm-first-lip-approx-by-harm-sheet:Dir-energy-on-bad-set}}&{\leq}  \bE({S_\eta},\cyl{0}{2-\overline\eta}) + C\eta_* E.
\end{aligned}
\end{equation}

In second place, we have

\begin{equation}\label{e:quasi_gollonzo1}
\begin{aligned}
\bE(T,\ball{0}{1},\pi) &\le\bE({S_\eta},\cyl{0}{1},\pi) \\ &\le\bE(\bG_{u^+}+\bG_{u^-},\cyl{0}{1},\pi)+2\be_T(\baseball{0}{1}\setminus K)+2|\pi|^2|\baseball{0}{1}\setminus K|\\
\overset{\eqref{eq:thm-first-lip-approx-by-harm-sheet:excees-measure-of-bad-set}}&{\le} \bE(\bG_{u^+}+\bG_{u^-},\cyl{0}{1},\pi) + C\eta_* E + 2|\pi|^2|\baseball{0}{1}\setminus K|\\
\overset{\eqref{e:grst}, \eqref{e:gollonzo2},\eqref{eq:excess-decay:reduction to graphs}}&{\le}  (2-\overline{\eta})^{2-\varepsilon}\bE({S_\eta},\cyl{0}{2-\overline\eta}) + C\eta_* E +C\overline\eta E.
\end{aligned}
\end{equation}

Using the height bound in Theorem \ref{lemma:height-bound}, for $\varepsilon<\varepsilon_0$ sufficiently small, we have
$$\bh (T, \cyl{0}{2-\overline{\eta}}, \pi_0) \leq C_h \left(\frac{\bE (T, \cyl{0}{4 - 2\overline{\eta}})}{1-\frac{\overline{\eta}}{2}}+ \bA\right)^{\frac{1}{2}} (2-\overline{\eta})^\frac{3}{2},$$
and thus
\begin{equation}\label{eq:milder-decay:height-bound}
  \spt(T)\cap\cyl{0}{2-\overline\eta}\subset\ball{0}{2}\,.  
\end{equation}

Since ${S_\eta}\res\ball{0}{2}=T\res\ball{0}{2}$, we obtain that

\begin{equation*}
\begin{aligned}
\bE^\flat(T,\ball{0}{1}) & \le \bE(T,\ball{0}{1},\pi)\\
\overset{\eqref{e:quasi_gollonzo1},\eqref{eq:milder-decay:height-bound}}&{\le}(2-\overline\eta)^{- (2-\varepsilon)}\left(\frac{2}{2-\overline\eta}\right)^m\bE(T,\ball{0}{2},\pi)+C\eta_*E+\overline\eta E\\
& = (2-\overline\eta)^{- (2-\varepsilon)}\left(\frac{2}{2-\overline\eta}\right)^m\bE^\flat(T,\ball{0}{2})+C\eta_*E+\overline\eta E \\
\overset{\eqref{e:bibo}}&{\le} \biggl[(2-\overline\eta)^{- (2-\varepsilon)}\left(\frac{2}{2-\overline\eta}\right)^m\ +C(\eta_* + \overline\eta)\biggr] \bE^\flat(T,\ball{0}{2}) .
\end{aligned}
\end{equation*}

Hence, since the constant $C$ in the last inequality is independent of the parameters $\eta_*,\overline\eta$, choosing the latter sufficiently small, we conclude \eqref{e:gollonzo}.
\end{proof}

\begin{proof}[Proof of Theorem \ref{thm:improved-excess-decay}]

Firstly, we want to prove that the assumptions \eqref{eq:lemma-milder-decay:small-excess-assump} and \eqref{eq:lemma:milder-decay:densities-p-q-close} of Lemma \ref{lemma:milder-excess-decay} are satisfied for every boundary point $q$ in a neighbourhood of $p$ (see (A) and (B) below). To this end, we notice that, since $p$ is a two-sided collapsed point, by Definition \ref{defi:collapsed}, for every \(\delta>0\), there exists \(\bar\rho =\bar\rho (\delta)\) small such that 

\begin{enumerate}[\upshape (i)]
\item \(\bE^\flat(T,\ball{p}{2\sigma})+ 4 \bA \sigma^2\le \delta\) for every $\sigma\leq \bar\rho$;

\item \(\Theta(T,q)\ge \Theta(T,p)=Q-\frac{Q^\star}{2}\) for all \(q\in \Gamma \cap \ball{p}{2\bar\rho}\).
\end{enumerate}

Next, since $\Theta (T, p)= Q-\frac{
Q^\star}{2}$, if the radius $\bar\rho$ is chosen small enough, we can assure that

\[
\|T\| (\ball{p}{4\bar\rho}) \leq \omega_m \left(Q-\frac{Q^\star}{2}+\frac{1}{8}\right) (4\bar\rho)^m\, .
\]

By a simple comparison, for $\eta$ sufficiently small, if $q\in \ball{p}{\eta}\cap \Gamma$ and $\bar\rho' =\bar \rho-\eta$, then

\begin{equation*}
\begin{aligned}
\|T\| (\ball{q}{4\bar\rho '}) &\leq \|T\| (\ball{p}{4\bar\rho}) \leq \omega_m \left(Q-\frac{Q^\star}{2} + \frac{1}{8}\right) (4\bar\rho)^m \leq \omega_m \left(Q-\frac{Q^\star}{2} + \frac{3}{16}\right) (4\bar\rho')^m\, .
\end{aligned}
\end{equation*}

Next, by the latter inequality and by the monotonicity formula, it follows that

\begin{equation*}
\begin{aligned}
\sigma^{-m} \|T\| (\ball{q}{\sigma}) &\leq e^{\bA (4\bar\rho'-\sigma)} (4\bar \rho')^{-m} \|T\| (\ball{q}{4\bar \rho'})\\
&\leq e^{\bA (4\bar\rho'-\sigma)}\omega_m \left(Q-\frac{Q^\star}{2}+\frac{3}{16}\right)\leq e^{4 \bA\bar\rho} \omega_m \left(Q-\frac{Q^\star}{2}+\frac{3}{16}\right),
\end{aligned}
\end{equation*}

for all $\sigma \leq 4\bar\rho' $. In particular, if $\bar\rho$ is chosen sufficiently small, by the last inequality we then conclude

\begin{equation}\label{eq:thm:improved-excess:bound-density}
\|T\| (\ball{q}{\sigma}) \leq \omega_m \left(Q-\frac{Q^\star}{2}+\frac{1}{4}\right) \sigma^m, \qquad \forall q\in \ball{p}{\eta}\cap \Gamma 
\;\mbox{and}\; \forall \sigma \leq 4\bar \rho' \, .
\end{equation}

So, the density of $T$ at $q$ is bounded above by \eqref{eq:thm:improved-excess:bound-density} and below by (ii). Set now $r:= \min \{\eta, \bar\rho'\}$. For all  points \(q\) in \(\ball{p}{r}\cap \Gamma\) we claim that
\begin{equation}\label{e:start}
\bE^\flat(T,\ball{q}{r}) \le 2^m \bE^\flat (T, \ball{p}{2r}) + C \bA^2 r^2 \overset{\textup{(i)}}{\leq} C\delta.
\end{equation}

Indeed let $\pi$ be a plane for which $\bE^\flat (T, \ball{p}{2r}) = \bE (T, \ball{p}{2r}, \pi)$. By the regularity of
$\Gamma$, we find a plane $\pi_q$ such that $|\pi - \pi_q|\leq C r \bA$ and 
$T_q \Gamma \subset \pi_q$. Then we can estimate
\begin{equation}\label{eq:thm:improved-excess:A}
\begin{aligned}
\bE^\flat(T,\ball{q}{r}) &\le\bE(T,\ball{q}{r},\pi_q) \le C \bE(T,\ball{p}{2r},\pi_q)\\
\overset{\textup{Triangular}}&{\le} C \bE^\flat (T,\ball{p}{2r}) + C r^2 \bA^2 \overset{\textup{(i)}}{\le} C\delta\, .
\end{aligned}
\end{equation}

We will show that the conclusions of the theorem hold for this particular radius $r$ which, without loss of generality, we assume to be $r=1$ and we also assume $p=0$. So, we have proved that we are under the assumptions of Lemma \ref{lemma:milder-excess-decay}, in fact, \eqref{eq:thm:improved-excess:A} and \eqref{eq:thm:improved-excess:bound-density} ensure the following properties for every $q\in \ball{0}{1} \cap \Gamma$

\begin{enumerate}[\upshape (A)]
\item $\bE^\flat (T, \ball{q}{1}) + \bA^2 \leq C \bE^\flat (T, \ball{0}{2}) + C \bA^2 \leq C\delta$,

\item $\|T\| (\ball{q}{s}) \leq (Q-\frac{2Q^\star - 1}{4}) \omega_m s^m$ for every $s\leq 1$.  
\end{enumerate}

We now fix any point \(q\in \Gamma\cap \ball{0}{1}\) and define $\mathfrak{m}(s) := \bE^\flat (T, \ball{q}{s})$. We claim that

\begin{equation}\label{e:dec_et}
\mathfrak{m}(s)\le Cs^{2-2\varepsilon}\max\{\mathfrak{m}({\textstyle{\frac{1}{4}}}),\mathfrak{m}({\textstyle{\frac{1}{2}}})\} + Cs^2\bA^2, \qquad \forall s\in (0,\frac12).
\end{equation}

In order to prove \eqref{e:dec_et}, we firstly prove for $s=2^{-k-1}$ and for all $k\in\N$ that 

\begin{equation}\label{eq:lemma:milder-decay:induction}
\mathfrak{m}(2^{-k-1})\leq C \max\{2^{(2\epsilon -2)k}\mathfrak{m}(\frac{1}{4}),2^{(2\epsilon -2)k+2}\mathfrak{m}(\frac{1}{2})\} + C 2^{-2k-4}\bA
\end{equation}

is valid and then we will show how to derive \eqref{e:dec_et} from \eqref{eq:lemma:milder-decay:induction}. The proof of \eqref{eq:lemma:milder-decay:induction} will be done by induction. Notice that inequality \eqref{eq:lemma:milder-decay:induction} is trivially true for $k=0$, indeed,

$$ \mathfrak{m}(\frac{1}{2}) \leq 2^2\mathfrak{m}(\frac12) \leq \max\{ \mathfrak{m}(\frac14), 2^2\mathfrak{m}(\frac12) \}. $$

If the inequality is true for $k_0\geq 0$, we want to show it for $k=k_0+1$. We set $\sigma=2^{-k-2}$ and notice that, by inductive assumption, we conclude that

\begin{equation}\label{eq:excess-decay-induction-step0}
\begin{aligned}
\mathfrak{m}(2^{-k-1}) \leq \mathfrak{m}(4\sigma) = \mathfrak{m}(2^{-k_0-1}) &\leq \max\{ 2^{(2\varepsilon -2)k_0}\mathfrak{m}(\frac14),2^{(2\varepsilon -2)k_0+2}\mathfrak{m}(\frac12) \}\\
&\leq \max\{\mathfrak{m}({\textstyle{\frac{1}{4}}}),\mathfrak{m}({\textstyle{\frac{1}{2}}})\} \le \mathfrak{m}(1) \overset{(\textup{A})}{\le} C\delta.
\end{aligned}
\end{equation}

Hence, provided we choose \(\delta=\delta(m,n, Q^\star,Q)>0\) small in (A) and consequently \(r\) is sufficiently small too, we are in position to apply Lemma \ref{lemma:milder-excess-decay} which assures that  

\begin{equation*}
\begin{aligned}
\mathfrak{m}(2^{-(k+1)-1}) = \mathfrak{m}(\sigma) \overset{\eqref{eq:thm:milder-excess-decay}}&{\le} \max\{2^{-2+2\varepsilon}\theta(2\sigma),2^{-4+4\varepsilon}\theta(4\sigma)\} \leq C 2^{-4+4\varepsilon}\theta(4\sigma)\\
&= C 2^{-4+4\varepsilon} \max\{\mathfrak{m}(4\sigma), M_0\bA^2(4\sigma)^2 \}\\
\overset{\eqref{eq:excess-decay-induction-step0}}&{\le} C\max\{ 2^{(2\varepsilon -2)(k_0+1)}\mathfrak{m}(\frac14),2^{(2\varepsilon -2)(k_0+1)+2}\mathfrak{m}(\frac12) \} + M_0\bA^2\sigma^2 ,
\end{aligned}
\end{equation*}

where we recall that $C$ and $M_0$ are both constants that depends on $m,n,Q^\star, Q$ and $\varepsilon$, which finishes our induction steps and proves \eqref{eq:lemma:milder-decay:induction}. To prove \eqref{e:dec_et}, we take $s\in(0,\frac12)$ and $k_s\in\N$ such that $s\in(2^{-k_s-2}, 2^{-k_s-1})$, hence, by \eqref{eq:lemma:milder-decay:induction},

$$ \mathfrak{m}(s) \leq \mathfrak{m}(2^{-k_s-1}) \leq C\max\{ 2^{(2\varepsilon -2)k_s}\mathfrak{m}(\frac{1}{4}), 2^{(2\varepsilon -2)k_s+2}\mathfrak{m}(\frac{1}{2})\} + C2^{-2k_s-4}\bA^2 , $$

taking into account in the last inequality that $s^{2-2\varepsilon}>4^{2-2\varepsilon}\cdot 2^{(2\varepsilon -2)k_s}$, we finish the proof of \eqref{e:dec_et}. We then conclude, for \(\rho\in (0,\frac{1}{2})\), that the following equation holds

\begin{equation}\label{e:dec_rrho}
\begin{aligned}
\bE(T,\ball{q}{\rho})\le \bE^\flat(T,\ball{q}{\rho}) \overset{\eqref{e:dec_et}}&{\le} C\rho^{2-2\varepsilon}\max\{\mathfrak{m}(\frac{1}{2}), \mathfrak{m}(\frac{1}{4})\} +C\rho^2\bA^2\\
&\le C\rho^{2-2\varepsilon}\bE^\flat(T,\ball{q}{1}) + C\rho^2\bA^2\\
\overset{(\mathrm{A})}&{\le} C \rho^{2-2\varepsilon} \bE^\flat (T, \ball{0}{2}) + C \rho^{2-2\varepsilon}\bA^2 . 
\end{aligned}
\end{equation}

Furthermore, by (A), the estimate is trivial for $\frac{1}{2} \leq \rho < 1$. For \(0<t<s<1\), define \(\pi(q,s)\) and \(\pi(q,t)\) to be the optimal planes for \(\bE^\flat(T,\ball{q}{t})\) and \(\bE^\flat(T,\ball{q}{s})\), respectively. So, \eqref{e:dec_rrho} implies 
$$
\begin{aligned}
|\pi(q,s)-\pi(q,t)|^2 &= \frac{1}{\|T\|(\ball{q}{s})}\int_{\ball{q}{s}}|\pi(q,t)-\pi(q,s)|^2d\|T\| \\
\overset{\eqref{eq:thm:improved-excess:bound-density}}&{\le} C\bE(T,\ball{q}{s},\pi(q,s))+C\bE(T,\ball{q}{t}),\pi(q,t))\\
\overset{\eqref{e:dec_rrho}}&{\le} C s^{2-2\varepsilon}\bE^\flat(T,\ball{0}{2})+Cs^{2-2\varepsilon}\bA^2 .
\end{aligned}
$$
Letting $t$ goes to $0$ in the last equations and thanks to the compactness of $G_m(\mathbb{R}^{m+n})$, we obtain the existence of a limit $\pi(q)$ such that
\begin{equation}\label{e:franco00}
|\pi(q)-\pi(q,\rho)|^2\le C \rho^{2-2\varepsilon}\bE^\flat(T,\ball{0}{2})+C \rho^{2-2\varepsilon}\bA^2\quad,\forall\,\rho< 1\,.
\end{equation}
Hence, for all $\rho\in (0,1)$, we conclude that
\begin{equation}
\begin{aligned}
\bE^\flat (T,\ball{q}{\rho}) &\leq \bE^\flat (T,\ball{q}{\rho}, \pi(q)) \\
\overset{\textup{Triangular}}&{\leq} C\bE^\flat (T,\ball{q}{\rho}, \pi(q,\rho)) + C |\pi(q,\rho)-\pi(q)|^2\\
&= C\bE^\flat(T,\ball{q}{\rho}) + C |\pi(q,\rho)-\pi(q)|^2\\
\overset{\eqref{e:dec_rrho},\eqref{e:franco00}}&{\leq} C s^{2-2\varepsilon}\bE^\flat(T,\ball{0}{2})+Cs^{2-2\varepsilon}\bA^2.
\end{aligned}   
\end{equation} 
which concludes the proof of \eqref{eq:decay-improved-excess}.

We now turn to \eqref{eq:thm-improved-excess-decay:height-bound}, let 
$$
S_\rho = T\res\left(B_\rho(q,\pi(q))\times B^n_\rho(q,\pi(q)^\perp)\right).
$$
Hence, we immediately have $T\res\ball{q}{\rho}=S_\rho\res\ball{q}{\rho}$. Moreover, arguing as in (B), (C) and (D) in the proof of Lemma \ref{lemma:milder-excess-decay}, we are under Assumption \ref{assump:general.asump-first-lip}, thus we can apply the Height bound (Lemma \ref{lemma:height-bound}) to obtain
\begin{equation}\label{eq:thm:height-bound01}
\bh(S_\rho,\cyl{q}{\rho}, \pi(q)) \leq C_h\left(\rho^{-1}\bE(S_\rho,\cyl{q}{2\rho},\pi(q)) + \bA\right)^{\frac12}\rho^{\frac32}, \ \ \forall \rho\in (0,\frac12).
\end{equation}

As in \eqref{e:nonsisa}, we obtain that
\begin{equation}\label{eq:thm:height-bound02}
\bE(S_\rho, \cyl{q}{\rho},\pi(q)) \leq C\bE^\flat(T,\ball{q}{\sqrt{2}\rho},\pi(q)), \ \ \forall \rho\in (0,\frac{1}{\sqrt{2}}).
\end{equation}

We are ready to conclude \eqref{eq:thm-improved-excess-decay:height-bound} as follows, for every $\rho\in (0,\frac1{2\sqrt{2}})$, 

\begin{equation}
\begin{aligned}
\bh(T,\ball{q}{\rho},\pi(q)) & = \bh(S_\rho,\ball{q}{\rho},\pi(q)) \overset{\eqref{eq:thm:height-bound01}}{\le} C_h\left(\rho^{-1}\bE(T,\cyl{q}{2\rho},\pi(q)) + \bA\right)^{\frac12}\rho^{\frac32} \\
\overset{\eqref{eq:thm:height-bound02}}&{\le} C_h\left(\rho^{-1}\bE(T,\cyl{q}{2\sqrt{2}\rho},\pi(q)) + \bA\right)^{\frac12}\rho^{\frac32},
\end{aligned}
\end{equation}

it is sufficient to apply the improved excess decay, \eqref{eq:decay-improved-excess}, to conclude the proof.
\end{proof}

\subsection{Uniqueness of tangent cones at two-sided collapsed points}

In the spirit of Theorem \ref{thm:uniqueness-tang-cones-decay} and Lemma \ref{lemma:Holder-continuity-tang-cones} which state the uniqueness of tangent cones and the H\"older continuity of the map $q\mapsto T_q$ for $(C_0, r_0, \alpha_0)$-almost area minimizing currents of dimension $2$, we prove the uniqueness of tangent and the H\"older continuity of the same map for area minimizing currents of arbitrary dimension $m$. We state below the analog of \cite[Theorem 6.3]{DDHM} when the boundary is taken with multiplicity $Q^\star$.

\begin{theorem}[Uniqueness of tangent cones at two-sided collapsed points]\label{thm:uniqueness-dimension-m}
Let $T, p, U$ and $r$ be as in Theorem \ref{thm:improved-excess-decay}. Then for all $q\in\ball{p}{r}\subset U$, we have that $q$ is a two-sided collapsed point with $\Theta^m(T,q)=\Theta^m(T,p)$ and there is a unique tangent cone $T_q = Q\a{\pi(q)^+}+(Q-Q^\star)\a{\pi(q)^-}$ to $T$ at $q$, where $\pi(q)$ is an $m$-dimensional plane. Moreover, for any $\varepsilon>0$, there is $C=C(\varepsilon)>0$, such that

\begin{equation}\label{eq:thm-uniq-dim-m:Holder}
|\pi (q) - \pi (z)| \leq C \left(r^{\eps-1} \left(\bE^\flat (T, \ball{p}{2r})\right)^{\frac{1}{2}} + \bA r^\varepsilon\right) |z-q|^{1-\varepsilon}, \ \forall z\in\ball{p}{r}.
\end{equation}
\end{theorem}
\begin{remark}
{Note that, as we have proved in Lemma \ref{lemma:collapsed-set-is-open} in dimension $2$ using the characterization of the tangent cones in the $2d$ setting, Theorem \ref{thm:uniqueness-dimension-m} ensures, for arbitrary dimension $m$, that the set of two-sided collapsed points is relatively open in $\Gamma$. Furthermore, it also guarantees that the density is constant in $\ball{p}{r}\cap\Gamma$.}
\end{remark}
The following proof goes along the same lines of \cite[Theorem 6.3]{DDHM} and we report it here for completeness's sake. 
\begin{proof}
We start taking $\pi(q)$ as the plane given by the improved excess decay, see Theorem \ref{thm:improved-excess-decay}. Let us now prove that, if $T_q$ is a tangent cone to $T$ at $q$ w.r.t. the sequence $\rho_k\to 0$, then its support is $\pi(q)$. By rescaling we have that
$$ \bE(T,\ball{q}{\rho_k}, \pi(q)) \leq C\bE(T_{q,\rho_k},\ball{0}{2}, \pi(q)).$$
The latter rescaling and the improved excess decay, i.e., \eqref{eq:decay-improved-excess}, furnish
\begin{equation}\label{eq:thm:uniqueness-dim-m:support-pi(q)}
\bE(T_{q,\rho_k},\ball{0}{2}, \pi(q)) \leq C\left(\frac{\rho_k}{r}\right)^{2-2\varepsilon}\bE^\flat(T,\ball{p}{2r})+C\rho_k^{2-2\varepsilon}r^{2\varepsilon}\bA^2, \ \ \forall \rho_k < r.
\end{equation} 
We now let $\rho_k\to 0$ in \eqref{eq:thm:uniqueness-dim-m:support-pi(q)} to conclude that $\bE(T_{q},\ball{0}{2}, \pi(q)) = 0$ and hence $T_q$ is supported in $\pi(q)$.
We conclude that the tangent cone is unique and, by a standard argument involving the Constancy Lemma as already used many times above, it takes the form

\[
T_q={Q(q)}\a{\pi(q)^+}+({Q(q)}-Q^\star)\a{\pi(q)^-}\,,
\]

for some \({Q(q)}\in \mathbb N\), since the tangent cone is an integral current. By assumption we have that $p$ is a two-sided collapsed point and $q\in U$, ${Q(q)}-\frac{Q^\star}{2} =\Theta(T,q)\ge \Theta^m(T,p)$. Moreover, by \eqref{eq:lemma:milder-decay:densities-p-q-close}, we obtain $Q(q) -\frac{Q^\star}{2} \leq \Theta^m(T,p) + \frac{1}{4}$ and therefore \(\Theta^m(T,q) = \Theta^m(T,p)\). Finally, in order to finish the proof of the theorem, for \(0<t<s<1\), we let \(\pi(q,s)\) and \(\pi(q,t)\) such that \(\bE^\flat(T,\ball{q}{t}) = \bE^\flat(T,\ball{q}{t}, \pi(q,t))\) and \(\bE^\flat(T,\ball{q}{s}) = \bE^\flat(T,\ball{q}{s}, \pi(q,s))\). We now take $0<t<\rho:=|q-z|<r$ and note that

\begin{equation*}
\begin{aligned}
|\pi(q,t)-\pi(z,t)|^2  &=\mint_{\ball{q}{t}\cap \ball{z}{t}}|\pi(q,t)-\pi(z,t)|^2d\|T\|\\
& \le \frac{C}{\omega_mt^m}\int_{\ball{q}{t}}|\vec T-\pi(q,t)|^2+\frac{C}{\omega_m t^m}\int_{\ball{z}{t}}|\vec T-\pi(z,t)|^2d\|T\|\\
& =C(\bE^\flat(T,\ball{q}{t})+\bE^\flat(T,\ball{z}{t}))\\
& \le C \left(\frac{\rho}{r}\right)^{2-2\varepsilon}  \bE^\flat (T, \ball{p}{2r}) + C \rho^{2-2\varepsilon}r^{2\varepsilon} \bA^2,
\end{aligned}
\end{equation*}

where, in the second line, we have used that \(\|T\|(\ball{p}{t}) \ge c t^m\), a simple consequence of the monotonicity formula. Hence, the latter inequality gives

$$ |\pi(q,t)-\pi(z,t)| \leq  C\left( r^{-2+2\varepsilon}  \bE^\flat (T, \ball{p}{2r}) +  r^{2\varepsilon} \bA^2\right)^{\frac{1}{2}}|q-z|^{1-\varepsilon},$$

we let $t$ goes to $0$ to conclude \eqref{eq:thm-uniq-dim-m:Holder}.
\end{proof}

We state an important corollary of Theorems \ref{thm:uniqueness-dimension-m} and \ref{thm:improved-excess-decay} which will be used often in the remaining chapters and relates the height when we change the reference plane to an optimal plane to the excess in $p$ instead of consider the tangent plane at $q$. 

\begin{corollary}[Height bound relative to tilted optimal planes]\label{cor:height-tilted-planes}
Let $\Gamma, T, p, q, \pi(q)$ and $r$ be as in Theorem \ref{thm:uniqueness-dimension-m} and let $\pi$ be an optimal plane for $\bE^\flat(T,\ball{p}{2r})$. If we denote by $\bp_\pi, \bp_\pi^\perp, \bp_{\pi(q)}$ and $\bp_{\pi(q)}^\perp$ respectively the
orthogonal projections onto $\pi, \pi^\perp, \pi (q)$ and $\pi (q)^\perp$, then, for all $q \in \Gamma\cap \ball{p}{r}$, we have
\begin{equation}\label{eq:thm-height-tilted:dist-optimal-planes}
|\pi (q) - \pi| \leq C (\bE^\flat(T,\ball{p}{2r})^{\frac12} + \bA r),
\end{equation}
\begin{equation}\label{eq:thm-height-tilted:subset-q-set} 
	\spt(T) \cap \ball{q}{\frac{r}{2}}  \subset \left\{ x\in\R^{m+n} : \abs{\bp^\perp_{\pi(q)}(x-q)} \le C (r^{-1} \bE^\flat(T,\ball{p}{2r})^\frac{1}{2} + \bA)^{\frac{1}{2}} \abs{x-q}^{\frac32}\right\}, 
\end{equation}
\begin{equation}\label{eq:thm-height-tilted:subset-p-set}
	\spt(T) \cap \ball{q}{\frac{r}{2}}  \subset  \left\{x\in\R^{m+n}  \colon \abs{\bp_\pi^\perp(x-q)} \le C (\bE^\flat(T,\ball{p}{2r})^{\frac12} + \bA r)^{\frac{1}{2}} \abs{x-q}\right\}.
\end{equation}
\end{corollary}
\begin{proof}

To prove \eqref{eq:thm-height-tilted:dist-optimal-planes}, we proceed as follows
\begin{equation}\label{e:proof:tilt_pi(q)1}
\begin{aligned}
|\pi-\pi(q)|^2 &\leq 2|\pi - \pi(p)|^2 + 2|\pi(p) - \pi(q)|^2 \\
\overset{\eqref{eq:thm-uniq-dim-m:Holder}}&{\leq} 2|\pi - \pi(p)|^2 + C(\bE^\flat(T, \ball{p}{2r})^\frac{1}{2}+\bA r)^2,
\end{aligned}
\end{equation}

and 

\begin{equation}\label{e:proof:tilt_pi(q)2}
\begin{aligned}
|\pi - \pi(p)|^2 &\leq C\frac{1}{\|T\|(\ball{p}{2r})}\int_{\ball{p}{2r}}\left(|\pi - \vec{T}|^2+|\vec{T} - \pi(p)|^2\right)d\|T\| \\
\overset{(\ast)}&{\leq} C\bE^\flat(T,\ball{p}{2r}) + C\frac{1}{\|T\|(\ball{p}{2r})}\int_{\ball{p}{2r}}|\vec{T} - \pi(p)|^2d\|T\|\\
\overset{\eqref{eq:decay-improved-excess}, (\ast)}&{\leq} C(\bE^\flat(T, \ball{p}{2r})^\frac{1}{2}+\bA r)^2,
\end{aligned}
\end{equation}

where in $(\ast)$ we have used the standard argument with the monotonicity formula to obtain a bound $\|T\|(\ball{p}{2r}) \geq cr^m$. Therefore \eqref{e:proof:tilt_pi(q)1} and \eqref{e:proof:tilt_pi(q)2} prove \eqref{eq:thm-height-tilted:dist-optimal-planes}. Note that \eqref{eq:thm-height-tilted:subset-q-set} follows immediately from \eqref{eq:thm-improved-excess-decay:height-bound}. We next observe that 

\begin{equation}\label{eq:tilting-cor:distance-between-projections}
    \abs{\bp_\pi^\perp -{\bp_{\pi(q)}}^\perp}^2= \abs{\bp_\pi -{\bp_{\pi(q)}}}^2\le C \abs{\pi- \pi(q)}^2.
\end{equation}

Furthermore, for $x\in \ball{q}{r}\cap \spt (T)$, it follows

\begin{equation*}
\begin{aligned}
|\bp_\pi^\perp (x-q)|^2 &\leq C|x-q|^2 |\bp_\pi^\perp - {\bp_{\pi(q)}}^\perp|^2+ C|{\bp_{\pi(q)}}^\perp (x-q)|^2\\
&\leq C|\pi - \pi(q)|^2|x-q|^2 + C|{\bp_{\pi(q)}}^\perp (x-q)|^2\\
\overset{\eqref{eq:thm-height-tilted:dist-optimal-planes}}&{\leq} C (\bE(T,\ball{p}{2r})^{\frac{1}{2}} + \bA r)^2 |x-q|^2 + C|{\bp_{\pi(q)}}^\perp (x-q)|^2\\
\overset{\eqref{eq:thm-height-tilted:subset-q-set}}&{\leq} C (\bE(T,\ball{p}{2r})^{\frac{1}{2}} + \bA r)^2 |x-q|^2 + C (r^{-1} \bE(T,\ball{p}{2r})^{\frac12} + \bA)\abs{x-q}^{\frac32}\\
&\leq C (\bE(T,\ball{p}{2r})^{\frac12} + \bA r)^{\frac{1}{2}} |x-q| ,
\end{aligned}
\end{equation*}

where in the last inequality the fact $|x-q|<r$ took place, thus the latter inequality proves \eqref{eq:thm-height-tilted:subset-p-set}.

\end{proof}

\subsection{Lusin type strong Lipschitz approximations}

As we remarked at the beginning of this section, Theorem \ref{thm:first-lip-approx-and-good-set} provides an approximation which is not enough for our purposes. In the last subsection, we used the harmonic approximation, Theorem \ref{thm:first-lip-approx-by-harmonic-sheet}, to obtain the superlinear excess decay, Theorem \ref{thm:improved-excess-decay}, which will now be used to provide our desired approximation with faster decays and stronger estimates as it is precisely stated in Theorem \ref{thm:lip-approx-superlinear-decay}.

\begin{assumption}\label{assump:superlinear-lip-approx}
Let $T$ and $\Gamma$ be as in Assumption \ref{assumptions} with $C_0=0$, $\pi$ be a $m$-dimensional subspace, $\{e_i\}_{i=1}^m$ a basis of $\R^m$ and $q\in\Gamma$. We use the following notations $\pi^\prime = \operatorname{span}(\bp_\pi(e_1),\cdots, \bp_\pi(e_{m-1})), \psi_1:(q+\pi^\prime)\to(q+\operatorname{span}(\bp_\pi(e_m))), \psi:\gamma\subset(q+\pi)\rightarrow (q+\pi)^\perp,\psi_2:(q+\pi^\prime)\rightarrow(q+\operatorname{span}(\bp_\pi(e_m)))\times (q+\pi)^\perp , \psi_2(x) = (\psi_1(x), \psi(x, \psi_1(x)))$ with $\gr{\psi_1} = \gamma, \Gamma=\gr{\psi_2}$, and $\psi$ is of class $C^{3,\alpha}$. We assume that
\begin{enumerate}[\upshape (i)]
    \item In the excess decay, Theorem \ref{thm:improved-excess-decay}, $p = 0 \in\Gamma$ and $r=1$,

    \item $\bE^\flat (T, \ball{0}{2}) + \bA < \varepsilon_1,$ where $\varepsilon_1 = \varepsilon_1(m,n,Q^\star,Q)>0$ is a small constant.
\end{enumerate}
\end{assumption}

We would like to point out that the approximation and the estimates in the following theorem hold for any point in a suitable ball of the fixed two-sided collapsed point $p=0$, compare with Theorem \ref{thm:first-lip-approx-and-good-set} where the approximation and estimates are build in balls centered in the fixed two-sided collapsed point $p$.

\begin{theorem}[Strong Lipschitz approximation]\label{thm:lip-approx-superlinear-decay}
Let $T, \Gamma, \psi$ and $\gamma=\bp_\pi(\Gamma)$ be as in Assumption \ref{assump:superlinear-lip-approx}, $q\in\Gamma\cap\ball{0}{1}$, $r<\frac{1}{8}$ and $\pi$ be a plane such that $T_q\Gamma\subset\pi$ and $\bE(T,\cyltilted{q}{4r}{\pi})<\varepsilon_1$. Then there are a closed set $K\subset\mathrm{B}_r (q, \pi)$ and a $\left(Q-{\textstyle \frac{Q^\star}{2}}\right)$-valued map $(u^+, u^-)$ on $\mathrm{B}_r (q, \pi)$ which collapses at the interface $(\gamma, Q^\star\a{\psi})$ satisfying the following estimates:
\begin{equation}
\Lip (u^\pm)\leq C  (\bE(T, \cyltilted{q}{4r}{\pi}) + \bA^2r^2)^{\sigmaexp}\label{eq:thm:superlinear-lip-approx:lip-bound}
\end{equation}
\begin{equation}
{\operatorname{ osc}} (u^\pm) \leq C (\bE(T, \cyltilted{q}{4r}{\pi}) + \bA^2r^2)^{\frac{1}{2}} r\label{eq:thm:superlinear-lip-approx:osc-bound}
\end{equation}
\begin{equation}
\bG_{u^\pm} \res[(K \cap \Omega^\pm) \times \pi^\perp] = T \res [(K\cap \Omega^\pm)\times \R^n]\label{eq:thm:superlinear-lip-approx:gr=T}
\end{equation}
\begin{equation}
|\mathrm{B}_r (q, \pi)\setminus K|\leq C (\bE(T, \cyltilted{q}{4r}{\pi}) + \bA^2r^2)^{1+\sigmaexp}r^m\label{eq:thm:superlinear-lip-approx:bad-set-estimate}
\end{equation}
\begin{equation}
\be_T(\mathrm{B}_r (q, \pi)\setminus K) \le  C (\bE(T, \cyltilted{q}{4r}{\pi}) + \bA^2r^2)^{1+\sigmaexp}r^m\label{eq:thm:superlinear-lip-approx:excess-measure-estimate}
\end{equation}
\begin{equation}
\int_{\mathrm{B}_r (q, \pi)\setminus K} \abs{Du}^2 \le  C (\bE(T, \cyltilted{q}{4r}{\pi}) + \bA^2r^2)^{1+\sigmaexp}r^m \label{eq:thm:superlinear-lip-approx:Dir-energy-bound}
\end{equation}
\begin{equation}
\left| \be_T(F) - \frac12 \int_F \abs{Du^\pm}^2\right| \le C (\bE(T, \cyltilted{q}{4r}{\pi}) + \bA^2r^2)^{1+\sigmaexp} r^m, \quad \forall F \subset \Omega^\pm\text{ measurable,}\label{eq:thm:superlinear-lip-approx:Dir-energy-excess-measure}
\end{equation}

where $\Omega^\pm$ are the two regions in which $\mathrm{B}_r (q, \pi)$ is divided by $\gamma$, $C\ge 1$ and $\sigmaexp\in ]0,\frac{1}{4}[$ are two positive constants which depend on $m, n, Q^\star$ and $Q$. 
\end{theorem}
\begin{proof}
Our strategy to prove this theorem is to go back to the interior estimates done in \cite{DS3}. So we will divide the proof into two steps, in Step 1 we will prove further estimates provided by the interior case which are needed to conclude our estimates \eqref{eq:thm:superlinear-lip-approx:lip-bound}-\eqref{eq:thm:superlinear-lip-approx:Dir-energy-excess-measure}, and in Step 2 we will exhibit how to obtain the theorem from the interior case.

\textbf{Step 1:} If we assume that $\varepsilon_1$ is smaller than the constant of \cite[Theorem 2.4]{DS3} (also denoted by $\varepsilon_1$), $x\in\spt(T)$, and the cylinder $\cyltilted{x}{4\rho}{\pi}$ does not intersect $\Gamma$ and is contained in $\cyltilted{q}{4r}{\pi}$. Then, \cite[Theorem 2.4]{DS3} provides a map $f: \mathrm{B}_\rho (x, \pi) \to \mathcal{A}_Q (\pi^\perp)$ (or $\mathcal{A}_{Q-Q^\star} (\pi^\perp)$) and a closed set $\bar K\subset \mathrm{B}_\rho (x, \pi) $ such that

\begin{equation}\label{eq:interior-lip-approx:lip-bound}
\Lip (f) \leq C_{21} \bE(T,\cyl{x}{4\rho})^{\sigmaexp},
\end{equation}
\begin{equation}
\bG_f\res (\bar K\times \R^n)=T\res (\bar K\times\R^{n}),\label{eq:interior-lip-approx:gr=T}
\end{equation}
\begin{equation}
|\mathrm{B}_\rho (x, \pi)\setminus \bar K| \leq C\bE(T,\cyl{x}{4\rho})^{1+\sigma}\rho^m, \label{eq:interior-lip-approx:bad-set-estimate}
\end{equation}
\begin{equation}
\left| \|T\| (\cyl{x}{\rho}) - Q\omega_m\rho^m - \frac{1}{2}\int_{\mathrm{B}_{\rho} (x, \pi)} |Df|^2\right| \leq C\bE(T,\cyl{x}{4\rho})^{1+\sigma}\rho^m, \label{eq:interior-lip-approx:mass-dir-energy}
\end{equation}

In order to simplify our notation, we assume that $\pi= \mathbb R^m \times \{0\}$ and use the shorthand notation $\mathrm{B}_t (x)$ for $\mathrm{B}_t (x, \pi)$. It remains to prove the analog of \eqref{eq:thm:superlinear-lip-approx:excess-measure-estimate}, \eqref{eq:thm:superlinear-lip-approx:Dir-energy-bound} and \eqref{eq:thm:superlinear-lip-approx:Dir-energy-excess-measure} when $q$ is replaced by the interior point $x$. Notice that \eqref{eq:interior-lip-approx:lip-bound} and \eqref{eq:interior-lip-approx:bad-set-estimate} give

\[ \int_{F\setminus \bar K} \abs{Df}^2 \le C \bE(T,\cyl{x}{4\rho})^{2\sigmaexp} \abs{\baseball{x}{\rho}\setminus \bar K} \le C \bE(T,\cyl{x}{4\rho})^{1+3\sigmaexp}\rho^m, \]

for every $F \subset \mathrm{B}_\rho (x)$ measurable, hence we achieve \eqref{eq:thm:superlinear-lip-approx:Dir-energy-bound} taking $F= \baseball{x}{\rho}$. Next recall that either $\|T\| (\cyl{x}{\rho}) - Q\omega _m\rho^m = \be_T (\baseball{x}{\rho})$ or $\|T\| (\cyl{x}{\rho}) - (Q-Q^\star)\omega _m\rho^m = \be_T (\baseball{x}{\rho})$, hence \eqref{eq:interior-lip-approx:mass-dir-energy} can be reformulated as

\[ \left|\be_T (\baseball{x}{\rho}) - \frac12\int_{\baseball{x}{\rho}} \abs{Df}^2\right| \le  C \bE(T,\cyl{x}{4\rho})^{1+\sigmaexp}\rho^m\, . \]
We also have that

\begin{equation}\label{eq:lip-approx:lipschitz-bound}
\begin{aligned}
\frac{1}{2} \int_{\baseball{x}{\rho}} \abs{Df}^2 \leq\left(\bE(T, \cyl{x}{4\rho}) +  C\bE(T, \cyl{x}{4\rho})^{1+\sigmaexp}\right)\rho^m\leq C\bE(T, \cyl{x}{4\rho}) \rho^m .
\end{aligned}
\end{equation}

The Taylor expansion of the area functional, \cite[Corollary 3.3]{DS2}, and \eqref{eq:interior-lip-approx:lip-bound} give

\begin{equation}\label{eq:lip-approx:excess-measure-dir-energy}
\left|\be_{\bG_f}(F)-\frac12\int_F \abs{Df}^2\right| \le C\Lip(f)^2 \int_{F}\abs{Df}^2 \le C \bE(T, \cyl{x}{4\rho})^{1+2\sigmaexp}\rho^m,
\end{equation}

for every $F \subset \baseball{x}{\rho}$ measurable. Therefore, by \eqref{eq:lip-approx:lipschitz-bound} and \eqref{eq:lip-approx:excess-measure-dir-energy}, we obtain 

$$
\begin{aligned}
\be_T(\baseball{x}{\rho}\setminus \bar K) &= \be_T(\baseball{x}{\rho}) - \be_{\bG_f}(\baseball{x}{\rho}\cap \bar K) \\
&\le \left|\be_T(\baseball{x}{\rho})-{\frac12} \int_{\baseball{x}{\rho}} \abs{Df}^2 \right|+ \left|{\frac12} \int_{\baseball{x}{\rho}\cap \bar K} |Df|^2- \be_{\bG_f}(\baseball{x}{\rho}\cap \bar K)\right|\\
&\quad + \int_{\baseball{x}{\rho}\setminus \bar K} \abs{Df}^2\\
&\le C\bE(T,\cyl{x}{4\rho})^{1+\sigmaexp}\rho^m,
\end{aligned}
$$

which is \eqref{eq:thm:superlinear-lip-approx:excess-measure-estimate}. Finally, \eqref{eq:thm:superlinear-lip-approx:Dir-energy-excess-measure} is implied by the last inequality, \eqref{eq:thm:superlinear-lip-approx:Dir-energy-bound} and \eqref{eq:lip-approx:excess-measure-dir-energy} as follows, for every $F\subset \baseball{x}{\rho}$ measurable we have

$$
\begin{aligned}
\left| \be_T(F) - \frac12 \int_F \abs{Df}^2\right| &\le \left|\be_{\bG_{
f}}(F\cap \bar{K}) - \frac12 \int_{F\cap \bar{K}} \abs{Df}^2\right|+ \be_T(F\setminus \bar{K}) + \frac12 \int_{F\setminus \bar{K}} \abs{Df}^2\\
&\le C \bE(T,\cyl{x}{4\rho})^{1+\sigmaexp}\rho^m .
\end{aligned}
$$

\textbf{Step 2:} Without loss of generality we assume that $T_q \Gamma = \R^{m-1}\times \{0\}$, $\pi = \R^m \times \{0\}$. We then use $\cyl{q}{s}$ in place of $\cyltilted{q}{s}{\pi}$, $\baseball{q}{s}$ in place of
$\mathrm{B}_s (q, \pi)$, and $\bp$ to be the orthogonal projection onto $\pi$. 
By Assumption \ref{assump:superlinear-lip-approx}, we have

\begin{equation}\label{eq:lip-approx:assumptions1}
\partial T \res \cyl{q}{4r} = Q^\star \a{\Gamma\cap \cyl{q}{4r}}\nonumber\\
\qquad\mbox{and}\qquad  \bp_\sharp (\partial T \res \cyl{q}{4r}) = Q^\star\a{\gamma\cap \baseball{\bp(q)}{4r}}.
\end{equation}

As in the previous sections, denote by $\Omega^+$ and $\Omega^-$ the two connected components of $\baseball{q}{4r}\setminus \gamma$, we have
\begin{equation}\label{Eq:SecondLip:PushForwardOnTheBase}
\bp_\sharp T\res \cyl{q}{4r} = Q \a{\Omega^+} + (Q-Q^\star) \a{\Omega^-}\, .
\end{equation} 
We denote $L_0$ be the $m$-cube $q+ [-r,r]^m$ and, for any natural number $k$, let $\mathcal{C}_k$ be a collection of $m$-cubes given by

$$ \mathcal{C}_k := \left \{L=q + r 2^{-k} x + [-2^{-k}r, 2^{-k} r]^m: x\in \mathbb Z^m, k\in\N, L\subset L_0, L\cap\baseball{q}{r}\neq \emptyset  \right \}.$$

Analogously to \cite[Section 7.2]{DDHM}, we take a number $N\in \N$ such that the $2^{4-N}\sqrt{m}r$-neighborhood of $\cup_{L\in \mathcal{C}_N} L$ is
contained in $\cyl{q}{4r}$ and we will proceed with the construction of a Whitney decomposition of 
\[
\tilde\Omega := \bigcup_{L\in \mathcal{C}_N} L \setminus \gamma.
\]
Here and in what follows we set
\[
{\operatorname{ sep}}\, (L, \gamma) := \min \{|x-y|:x\in \gamma, y\in L \}\, .
\] 
We firstly define the following sets of $m$-cubes $\mathcal{R}_N = \mathcal{C}_N$,

\begin{equation*}
\begin{aligned}
\mathcal{W}_N &:= \left\{ L\in\mathcal{R}_N: \operatorname{diam}(L) \leq \frac{1}{16} \operatorname{sep}(L,\gamma) \right\}.
\end{aligned}
\end{equation*}

If $L\in\mathcal{R}_N\setminus\mathcal{W}_N$, we subdivide $L$ in $2^m$ $m$-subcubes of side $2^{-(N+1)} r$ and assign them to $\mathcal{R}_{N+1}$. We proceed inductively to define $\mathcal{W}_k$ and $\mathcal{R}_{k+1}$ for every $k\geq N$. Therefore, we obtain a Whitney decomposition $\mathcal{W} = \cup_{k\geq N} \mathcal{W}_k$ which is a collection of closed dyadic $m$-cubes such that 
\begin{align}
\operatorname{int}(L)\cap\operatorname{int}(H) &= \emptyset, \ \text{for all} \ L,H\in\mathcal{W}, \label{eq:second-lip:empty-intersec}\\
\Omega^+\cup\Omega^- &\subset \cup_{L\in\mathcal{W}}L, \label{eq:second-lip:covers-omega}\\
\frac{15}{16}\frac{1}{32} {\operatorname{ sep}}\, (L, \gamma) &<\; \diam(L) \le \frac{1}{16}{\operatorname{ sep}}\, (L, \gamma), \ L\in\mathcal{W}.\label{eq:second-lip:diam-sep}
\end{align}

Note that \eqref{eq:second-lip:empty-intersec} readily follows from the construction. In regard to \eqref{eq:second-lip:covers-omega}, we take $z\in\Omega^\pm$ to obtain two mutually exclusive cases that could happen, namely, either there exists $L\in\mathcal{W}_N$ such that $z\in L$, or for every $L\in\mathcal{W}_N$ results $z\notin L$. In the first case we finish readily the proof. In the second case take $L$ such that $L\in\mathcal{R}_N\setminus\mathcal{W}_N$ and $z\in L$. Then we may subdivide it and pass to the next generation and find a new cube $L'\in\mathcal{R}_{N+1}$ such that $z\in L'$. Now we may apply the same reasoning inductively and construct a sequence $L_{k,j_{k,z}}$, $k\ge N$ such that $z\in L_{k,j_{k,z}}$, $L_{k,j_{k,z}}\subseteq L_{k+1,j_{k+1,z}}$, and $L_{k,j_{k,z}}\notin\mathcal{W}$ for every $k\ge N$. If the sequence $(L_{k,j_{k,z}})_k$ is not constant for large values of $k$, then the diameters of $L_{k,j_{k,z}}$ goes to zero as $k$ goes to infinity, and thus we obtain for sufficiently large $k$ that 

$$ \diam(L_{k,j_{k,z}}) \leq \frac{1}{16}\operatorname{sep}(L_{N,j_{N,z}}, \gamma) \leq \frac{1}{16}\operatorname{sep}(L_{k,j_{k,z}}, \gamma), $$

since $L_{k,j_{k,z}}\subset L_{N,j_{N,z}}$ which ensures that $L_{k,j_{k,z}}\in\mathcal{W}$ and therefore \eqref{eq:second-lip:covers-omega}. To prove \eqref{eq:second-lip:diam-sep}, observe that $\operatorname{sep}(L,\gamma)\le\operatorname{sep}(\tilde{L},\gamma)+\diam(L)$ for every $L\in\mathcal{C}_k, \tilde{L}\in\mathcal{C}_{k-1}$ and $L\subseteq\tilde{L}$. 
By construction for each $L\in\mathcal{W}_k$ there exists $\tilde{L}\in\mathcal{R}_{k-1}\setminus\mathcal{W}_{k-1}$ such that $L\subseteq\tilde{L}$. Thus $\frac{15}{16}\operatorname{sep}(L, \gamma)\le\operatorname{sep}(\tilde{L},\gamma)$, $2\diam(L)=\diam(\tilde{L})>\frac1{16}\operatorname{sep}(\tilde{L},\gamma)$. So $\diam(L)>\frac1{32}\operatorname{sep}(\tilde{L},\gamma)\ge\frac{15}{16}\frac1{32}\operatorname{sep}(L,\gamma)$. 

Another important property of this family is that if a $m$-cube stops then its neighbours of next generation must also stop. Precisely, let $H\in\mathcal{W}_j,L\in\mathcal{C}_{j+1}$ and $H\cap L \neq\emptyset$, then
\begin{equation}\label{eq:second-lip:ngbh-also-stops}
\operatorname{sep}(L,\gamma) \geq \operatorname{sep}(H,\gamma) - \operatorname{diam}(L) \overset{\eqref{eq:second-lip:diam-sep}}{\geq} 16\operatorname{diam}(H) - \operatorname{diam}(L)\ge31\operatorname{diam}(L) \geq 16\operatorname{diam}(L).
\end{equation}
The chain of inequalities above guarantees that $\operatorname{sep}(L,\gamma)\geq 16\operatorname{diam}(L)$ which is the very definition of the family $\mathcal{W}_{j+1}$, i.e., $L\in\mathcal{W}$.

We denote with $c_L$ the center of the $m$-cube $L\in \mathcal{W}$ and set $r_L:= 3\, \diam(L)$ so that 

\begin{equation}\label{eq:strong-lip:whitney:L-subset-ball-one-quarter}
    L \subset \baseball{c_L}{\frac14 r_L}.
\end{equation}

We claim that, for each cube $L$, the current $T$ restricted to the cylinder $\cyl{c_L}{4r_L}$ satisfies the assumptions of \cite[Theorem 2.4]{DS3}. Firstly note that, by \eqref{eq:second-lip:diam-sep}, we have $\cyl{c_L}{4r_L}\cap\Gamma = \emptyset$ and thus $\partial T \res \cyl{c_L}{4r_L}=0$, we also obtain by the choice of $N$ that $\ball{c_L}{6r_L}\subset \ball{q}{4r}$. Moreover, either $\baseball{c_L}{4r_L} \subset \Omega^+$ or $\baseball{c_L}{4r_L}\subset \Omega^-$ and thus by \eqref{Eq:SecondLip:PushForwardOnTheBase} we have

$$\bp_\sharp T \res \cyl{c_L}{4r_L} = \left\{\begin{array}{ll} 
Q \a{\baseball{c_L}{4r_L}},&\text{if } c_L\in\Omega^+, \\
(Q-Q^\star) \a{\baseball{c_L}{4r_L}},&\text{if } c_L\in\Omega^- .
\end{array}\right.  $$ 

It remains to prove that the excess is small enough. Towards this goal we will make use of the excess decay, Theorem \ref{thm:improved-excess-decay}. In fact, we start distinguishing the two cases $r_L = 2^{-N}r$ and $r_L< 2^{-N}r$.
If $r_L = 2^{-N}r$, rescaling the excess and by the assumptions of the theorem, we easily obtain that
\[ 
\bE(T,\cyl{c_L}{4r_L}) \le 2^{mN} \bE(T,\cyl{q}{4r}) < \varepsilon_1. 
\]
Now, for each $L \in \mathcal{W}$ with $r_L < 2^{-N}r$, we let $x_L$ be the orthogonal projection of $c_L$ on $\gamma$ and $q_L\in \Gamma$ be the point $(x_L, \psi(x_L))$. The first inequality of \eqref{eq:second-lip:diam-sep} implies that 
\begin{equation*}
\cyl{c_L}{4r_L}\subset \cyl{q_L}{15r_L}.
\end{equation*}
From our choice of $N$, taking $\varepsilon_1$ smaller if necessary, by \eqref{eq:thm-height-tilted:subset-p-set}, we have

\begin{equation}\label{eq:lip-approx:cyl-and-balls-subset-of-cyl-4r}
\spt(T)\cap\cyl{q_L}{16r_L}\subset \ball{q_L}{17r_L}\subset \cyl{q}{4r}.
\end{equation}

So, we deduce that

\begin{equation}
\begin{aligned}
\bE(T, \cyl{c_L}{4r_L}) \overset{\eqref{eq:lip-approx:cyl-and-balls-subset-of-cyl-4r}}&{\leq} \bE(T, \ball{q_L}{17r_L},\pi)\\
\overset{\textup{Triangular}}&{\le} C\bE(T, \ball{q_L}{17r_L},\pi(q_L)) + C|\pi - \pi(q_L)|\\
&\leq C\bE (T, \cyl{q}{4r}) + \bA^2 r^2< \varepsilon_1,
\end{aligned}
\end{equation}

where in the last inequality we have used the excess decay, Theorem \ref{thm:improved-excess-decay}, and \eqref{eq:thm-height-tilted:dist-optimal-planes}. So, provided $\varepsilon_1$ is chosen sufficiently small, we can apply Step 1 in every cylinder $\cyl{c_L}{4r_L}$ and obtain either a $Q$-valued or a $(Q-Q^\star)$-valued map $f_L$ on each half ball $\basehalfball{+}{c_L}{r_L}$ or $\basehalfball{-}{c_L}{r_L}$ and a closed set $K_L \subset \basehalfball{\pm}{c_L}{r_L}$ such that 

\begin{equation}
\Lip(f_L) \le  C \bE(T,\cyl{c_L}{4r_L})^{\sigmaexp}, \label{eq:strong-approx-in_L_1}
\end{equation}
\begin{equation}
\bG_{f_L}\res(K_L \times \R^n)=  T\res(K_L \times \R^n), \label{eq:strong-lip:graph-equal-current-on-L}
\end{equation}
\begin{equation}
\abs{ \baseball{c_L}{r_L}\setminus K_L} \le  C \bE(T,\cyl{c_L}{4r_L}) ^{1+\sigmaexp} r_L^m, \label{eq:strong-approx-in_L_2}
\end{equation}
\begin{equation}
\be_T(\baseball{c_L}{r_L}\setminus K_L) \le  C \bE(T,\cyl{c_L}{4r_L}) ^{1+\sigmaexp}r_L^m,
\end{equation}
\begin{equation}
\int_{\baseball{c_L}{r_L}\setminus K_L} \abs{Df_L}^2 \le  C \bE(T,\cyl{c_L}{4r_L})^{1+\sigmaexp}r_L^m, \label{eq:strong-approx-in_L_3}
\end{equation}
\begin{equation}
\left| \be_T(F) - \frac12 \int_F \abs{Df_L}^2\right| \le  C \bE(T,\cyl{c_L}{4r_L})^{1+\sigmaexp} r_L^m, \quad \forall F \subset \baseball{c_L}{r_L} \text{ measurable}.	 \label{eq:strong-approx-in_L_4}
\end{equation}

Next, for each $L$ we let $\mathcal{N}_{\geq} (L)$ be what we call the neighboring $m$-cubes in $\mathcal{W}$ with larger radius, i.e. 
$$
\mathcal{N}_{\geq}(L)=\{ H \in \mathcal{W} \colon H \cap L \neq \emptyset, r_H \ge r_L \}.
$$

We use the good sets provided by the interior approximation to define 

$$
\begin{aligned}
K'_L = K_L \cap \bigcap_{H\in \mathcal{N}_{\geq} (L)} K_{H}, \quad K^+ = \bigcup_{L \in \mathcal{W}, L\subset \Omega^+} K'_L\cap L, \quad K^- =\bigcup_{L \in \mathcal{W}, L\subset \Omega^-} K'_L\cap L.
\end{aligned}
$$

Note that, if $H\in\mathcal{N}_{\geq}(L)$, we already know by \eqref{eq:strong-lip:whitney:L-subset-ball-one-quarter} that $L\subset\baseball{c_L}{\frac{1}{4}r_L}$ and then, since $H\cap L\neq \emptyset$, we deduce that $L \subset B_{r_{H}} (c_{H})$. Such a fact guarantees that $K^{\prime}_L\neq\emptyset$ and the following:
\begin{equation}\label{eq:strong-lip:L-minus-K'_L}
\begin{aligned}
|L\setminus K'_L|&\leq \abs{L\setminus K_L}+\sum_{H\in\mathcal{N}_{\geq}(L)}\abs{L\setminus K_H}\leq \abs{L\setminus K_L}+\sum_{H\in\mathcal{N}_{\geq}(L)}\abs{\baseball{c_H}{r_H}\setminus K_H} \\
\overset{\eqref{eq:strong-approx-in_L_2}}&{\leq} C \bE(T, \cyl{q}{4r})^{1+\sigmaexp}r_L^m <C \varepsilon_1 r_L^m,
\end{aligned}
\end{equation}

where we also use that the cardinality of $\mathcal{N}_{\geq} (L)$ is bounded by a geometric constant $C^\prime$ and $r_H$ and $r_L$ are comparable by construction of the Whitney decomposition, since $H\cap L\neq \emptyset$, $H,L\in\mathcal{W}$. 

Furthermore, we set up two functions defined on $K^+$ and $K^-$, respectively, given by $\tilde{u}^+(x) := f_L(x)$, for any $x\in L\cap K^+, L\in\mathcal{W}$, and $\tilde{u}^-(x):= f_L(x)$, for all $x\in L\cap K^-, L\in\mathcal{W}$, thanks to \eqref{eq:strong-lip:graph-equal-current-on-L} these functions are well defined. Note that these functions are defined on each square in $\mathcal{W}$ which are away from the boundary $\gamma$, it means that we have to properly extend these functions in order to have a $(Q-\frac{Q^\star}{2})$-valued map which collapses at the interface. Indeed, we will refine this idea in the sequel.

We now claim the validity of the following: 

\begin{equation}
\Lip (\tilde{u}^\pm) \leq  C(\bE(T, \cyl{q}{4r}) + \bA^2r^2)^{\sigmaexp},\label{e:muta1}
\end{equation}
\begin{equation}
\bG_{\tilde{u}^\pm}\res (K^\pm\times \R^n) =  T \res (K^\pm\times \R^n),\label{e:muta2}
\end{equation}
\begin{equation}
\be_T (L\setminus K'_L) \leq   C (\bE(T, \cyl{q}{4r}) + \bA^2r^2)^{1+\sigmaexp}r_L^m,\label{e:mutanda4.5}
\end{equation}
\begin{equation}
\int_{L\setminus K'_L} \abs{D\tilde{u}^\pm}^2 \le  C(\bE(T, \cyl{q}{4r}) + \bA^2r^2)^{1+\sigmaexp}r_L^m.\label{e:mutanda5}
\end{equation}
Before the proof of this claim, we show how to prove \eqref{eq:thm:superlinear-lip-approx:lip-bound}-\eqref{eq:thm:superlinear-lip-approx:Dir-energy-excess-measure} from it. Define the good set as $K = K^+\cup K^-$ and notice that in view of these last inequalities, up to now, we have finished the proof of \eqref{eq:thm:superlinear-lip-approx:lip-bound}, \eqref{eq:thm:superlinear-lip-approx:osc-bound}, and \eqref{eq:thm:superlinear-lip-approx:gr=T}. Observe that 

\begin{equation}\label{eq:strong-lip:sum-r_L}
\sum_{L\in \mathcal{W}} r_L^m \leq C r^m\, ,
\end{equation}

which furnishes \eqref{eq:thm:superlinear-lip-approx:bad-set-estimate}, \eqref{eq:thm:superlinear-lip-approx:excess-measure-estimate} and \eqref{eq:thm:superlinear-lip-approx:Dir-energy-bound} by summing over $L\in\mathcal{W}$, respectively, \eqref{eq:strong-lip:L-minus-K'_L}, \eqref{e:mutanda4.5} and \eqref{e:mutanda5}.  Regarding \eqref{eq:thm:superlinear-lip-approx:Dir-energy-excess-measure}, we proceed as follows, fix a measurable set $F\subset \Omega^\pm$ and observe that, for any $m$-cube $L$ in the Whitney decomposition $\mathcal{W}$ of $\Omega^\pm$ we have that
$$
\begin{aligned}
\left| \be_T (F\cap L) - \frac{1}{2} \int_{F\cap L} |D\tilde{u}^+|^2\right| \overset{\textup{Triangular}}&{\leq} \left| \be_T (F\cap L\cap K^\pm) - \frac{1}{2} \int_{F\cap L\cap K^\pm} |D\tilde{u}^+|^2\right|\\
&\quad+ \be_T (L\setminus K^\pm) + \Lip (\tilde{u}^+)^2 |L\setminus K^\pm|\\
\overset{\eqref{eq:strong-lip:L-minus-K'_L},\eqref{e:muta1},\eqref{e:mutanda4.5}}&{\leq} \left| \be_T(F\cap L\cap K^\pm) - \frac12 \int_{F \cap L\cap K^\pm} \abs{Df_L}^2\right| \\
&\quad + C (\bE(T, \cyl{q}{4r}) + \bA^2r^2)^{1+\sigmaexp} r_L^m\\
\overset{\eqref{eq:strong-approx-in_L_4}}&{\leq} C (\bE(T, \cyl{q}{4r}) + \bA^2r^2)^{1+\sigmaexp} r_L^m\, .
\end{aligned}
$$
By \eqref{eq:strong-lip:sum-r_L}, summing over $L\in\mathcal{W}$, we obtain \eqref{eq:thm:superlinear-lip-approx:Dir-energy-excess-measure}. Now, we turn our attention to the proof of the claim.

We start with the proof of the Lipschitz bound in \eqref{e:muta1}, we let $H, L \in \mathcal{W}$ with $\diam(H) \ge \diam(L)$ and $x \in H\cap K_H, y \in L\cap K_L$, hence

\begin{itemize}
    \item If $H \cap L \neq \emptyset$, by the very definitions of $\tilde{u}$ and $K^\pm$, we know that $\tilde{u}^\pm=f_H$ on $K^\pm\cap H$. Since $r_H\geq r_L$, $H\in\mathcal{N}_{\geq}(L)$ and then $K_L^{\prime}\neq\emptyset$ as mentioned above. So, we can take any $z\in K_L^{\prime}\subset K_L\cap K_H$, to have that the Lipschitz bound on each $m$-cube, i.e., \eqref{eq:strong-approx-in_L_1}, ensures
    $$\G (\tilde{u}^\pm (x), \tilde{u}^\pm (y))\leq \G (\tilde{u}^\pm (x), \tilde{u}^\pm (z))+\G (\tilde{u}^\pm (z), \tilde{u}^\pm (y))\leq  C\bE(T,\cyl{q}{4r})^{\sigmaexp} |x-y|.$$
    
    \item If $H \cap L = \emptyset$, let $x_\gamma, y_\gamma\in\gamma$ such that $\dist(x,\gamma)=\dist(x,x_\gamma)$ and $\dist(y,\gamma)=\dist(y,y_\gamma)$. We directly obtain that \begin{equation}\label{eq:lip-approx:01}
    \begin{aligned}
    \G (\tilde{u}^\pm (x), Q^\pm \a{\psi (x_\gamma)}) &\leq C\bE(T,\cyl{q}{4r})^{\frac{1}{2}}|x-x_\gamma|,\\
    \G (\tilde{u}^\pm (y), Q^\pm \a{\psi (y_\gamma)}) &\leq C\bE(T,\cyl{q}{4r})^{\frac{1}{2}}|y-y_\gamma|,
    \end{aligned}
    \end{equation}
    where $Q^+ = Q$ and $Q^- = Q - Q^\star$. Indeed, both inequalities are due to the following facts: $\dist(x, \gamma)$ and $r_L$ are comparable, see \eqref{eq:second-lip:diam-sep}, and $\spt(\tilde{u}^{\pm}(x)),(x_{\gamma},\psi(x_{\gamma}))\in\spt(T)$ then we can readily apply the height bound, in the cylinder $\cyl{x_\gamma}{16r_L}$, given in \cite[Proposition 7.3]{DDHM} (this proposition is just a corollary of \cite[Thm 6.3 and Assump 7.1 (iv)]{DDHM} which are the counterparts of Theorem \ref{thm:improved-excess-decay} and Corollary \ref{cor:height-tilted-planes}, so it can be derived with the exact same proof). Note also that, by the regularity of \(\Gamma\), we obtain $|\psi (x_\gamma)-\psi (y_\gamma)|\leq C\bA^{\sfrac{1}{2}} |x_\gamma-y_\gamma|$. As a consequence, by \eqref{eq:lip-approx:01}, we can estimate
    $$
    \begin{aligned}
    \G (\tilde{u}^\pm (x), \tilde{u}^\pm (y)) &\leq \G (\tilde{u}^\pm (x), Q^\pm \a{\psi (x_\gamma)}) + ({Q^\pm})^{\frac{1}{2}} |\psi (x_\gamma) - \psi (y_\gamma)| \\
    &\quad + \G (\tilde{u}^\pm (y), {Q^\pm} \a{\psi (y_\gamma)}) \leq C\left( \bE(T,\cyl{q}{4r}) + \bA^2 r^2 \right)^\sigmaexp |x-y|
    \end{aligned}
    $$
    where we have used that $\sigmaexp\leq \frac{1}{4}$ and that
    \[
    |x-x_\gamma|+|y_\gamma-y|=\dist(x,\gamma)+\dist(y,\gamma)\le C(r_L+r_H)\le Cr_H\le C\frac{3\sqrt{m}}2|x-y|.
    \] 
    To see the last inequality observe that when $H \cap L = \emptyset$ by \eqref{eq:second-lip:ngbh-also-stops} we have that $|x-y|\ge\frac{2r_H}{3\sqrt{m}}$.
\end{itemize}
In particular, we have also proved that $(\tilde{u}^+, \tilde{u}^-)$ has a Lipschitz extension to $(K^\pm\cup \gamma)\cap \baseball{q}{r}$ which, by \eqref{eq:lip-approx:01}, collapses at the interface $(\gamma \cap \baseball{q}{r}, Q^\star\a{\psi})$. We next extend $\tilde{u}^\pm$ to the whole $\Omega^\pm$, and denote by $u^\pm$, keeping the Lipschitz estimate \eqref{e:muta1} up to a multiplicative geometric constant, c.f., \cite[Theorem 1.7]{DS1}. Finally, inequality \eqref{e:mutanda5} follows directly by the estimates on the bad set and the Lipschitz bound, i.e., \eqref{eq:strong-lip:L-minus-K'_L} and \eqref{e:muta1}. Concerning inequality \eqref{e:mutanda4.5}, we obtain it by \eqref{eq:strong-approx-in_L_4} and \eqref{e:muta1}. To conclude the proof of the theorem, we notice that equation \eqref{e:muta2} follows from \eqref{eq:strong-lip:graph-equal-current-on-L}.

\end{proof}

\section{Center manifolds \texorpdfstring{$\mathcal{M}^\pm$}{Lg} with boundary \texorpdfstring{$\Gamma$}{Lg}}

In this section we work under the assumption that $0\in\R^{m+n}$ is a two-sided collapsed point and that $T_0\Gamma =\R^{m-1}\times\{0\}$ and therefore, from Theorem \ref{thm:uniqueness-dimension-m}, the tangent cone of $T$ at $0$ is $Q\a{\pi_0^+}+(Q-Q^\star)\a{\pi_0^-}$, where 
\[
\pi_0^\pm = \{x\in \R^{m+n} :\pm x_m>0, x_{m+1} = \ldots = x_{n+m}=0\}\, .
\]

Following the notation that we have used up to now, we denote by $\gamma$ the projection onto $\pi_0$ of $\Gamma$ and, given any sufficiently small open set $\Omega\subset \pi_0$ in $\R^m$ which is contractible and contains $0$,
we denote by $\Omega^\pm$ the two regions in which $\Omega$ is divided by $\gamma$, i.e., the portions on the right and left of $\gamma$. In this section, we build two distinct  $m$-dimensional submanifolds $\mathcal{M}^\pm$ of class $C^3$ which will be called, respectively, \textbf{left and right center manifolds}. Both center manifolds will have $\Gamma\cap \cyltilted{0}{\sfrac{3}{2}}{\pi_0}$ as their boundary, when considered as submanifolds in the cylinder $\cyltilted{0}{\sfrac{3}{2}}{\pi_0}$ and they will be $C^{3, \kappa}$ for a suitable positive $\kappa$ up to the boundary. Additionally, at each point $p\in \Gamma \cap \cyltilted{0}{\sfrac{3}{2}}{\pi_0}$ the tangent space to both manifolds will be the same which is the tangent cone to $T$ at $p$ denoted by $\pi (p)$ as in Theorem \ref{thm:uniqueness-dimension-m}. In particular $\mathcal{M} = \mathcal{M}^+\cup \mathcal{M}^-$ will be
a $C^{1,1}$ submanifold in $\cyltilted{0}{\sfrac{3}{2}}{\pi_0}$ without boundary. 

Our aim in this section is to provide a new approximation of the current $T$, the way we will do this is building the center manifolds $\mathcal{M}^\pm$ which can be understood as an average of the sheets of the current $T$ on each side of $\Gamma$, in the construction of the center manifold we will fabricate maps $\mathcal{N}^\pm$ which are defined in $\mathcal{M}^\pm$ and show that these maps $\mathcal{N}^\pm$ approximate the current in the sense of the Lipschitz approximations furnished in the previous chapters, e.g., Theorem \ref{thm:lip-approx-superlinear-decay}. With respect to the final argument of this work to conclude the proof of Theorem \ref{thm:main-theorem} using Theorem \ref{thm:collapsed-is-regular}, we desire to prove that $\mathcal{N}^\pm$ is identically zero and thus the current has to satisfies $T\res\cyltilted{0}{\sfrac{3}{2}}{\pi_0} = Q\a{\mathcal{M}^+}+(Q-Q^\star)\a{\mathcal{M}^-}$ which assures that $0$ is a two-sided regular point of $T$. This strategy will be developed in the remaining part of this work where we begin with the construction of the center manifolds and the approximating maps, after that we use the theory of $(Q-\frac{Q^\star}{2})$-Dir minimizing maps, see Section \ref{sec:linear}, to obtain that $\mathcal{N}^\pm|_{\mathcal{M}^\pm} \equiv 0$. 

\subsection{Construction of the Whitney decomposition}

Since the algorithm is the same for both sides of $\gamma$, it means that we can repeat the same frame to build both center manifolds. We can focus without loss of generality on the construction of $\mathcal{M}^+$. We start by describing a procedure which furnishes a suitable Whitney-type decomposition of $\basehalfball{+}{0}{3/2}$ with cubes whose sides are parallel to the coordinate axes and have sidelength $2 \ell(L)$. The center of any such cube $L$ will be denoted by $c_L$ and its sidelength will be denoted by $2 \ell (L)$. 
We start by introducing a family of dyadic cubes $L\subset \pi_0$ in the following way: for $j\geq N_0$, where $N_0$ is an integer whose choice will be specified below,
we introduce the families 
\[
\sC_j:=\{L:\,L\text{ is a dyadic cube of side }\ell(L)=2^{-j}\text{ and } \basehalfball{+}{0}{3/2}\cap L\neq\emptyset\}\, ,
\]
where we recall that, for $s>0$, $\basehalfball{\pm}{0}{s}$ are the connected components of $\baseball{0}{s}\setminus \gamma$. For each $L$ define a radius
\[
r_L:=M_0\sqrt m\ell(L)\,,
\]
with $M_0\geq 1$ to be chosen later.  
We then subdivide $\sC := \cup_j \mathscr{C}_j$ into, respectively, \textbf{boundary cubes} and \textbf{non-boundary cubes}
\begin{align*}
\sC^\flat & :=\{L\in \sC:\,\dist(c_L,\gamma)<64 r_L\}, \quad \sC_j^\flat = \mathscr{C}^\flat \cap \mathscr{C}_j,\\
\sC^\natural & :=\{L\in \sC:\,\dist(c_L,\gamma)\ge 64 r_L\}, \quad \mathscr{C}^\natural_j = \mathscr{C}^\natural \cap \mathscr{C}_j.
\end{align*}

Observe that some boundary cubes can be completely contained in $\basehalfball{+}{0}{3/2}$. For this reason we prefer to use the term ``non-boundary'' rather than ``interior'' for the cubes in $\sC^\natural$. Indeed in what follows, without mentioning it any further, we will often use the same convention for
several other subfamilies of $\mathscr{C}$.

\begin{definition}
If $H, L \in \mathscr{C}$ we say that:
\begin{itemize}
\item $H$ is a \textbf{descendant} of $L$ and $L$ is an \textbf{ancestor} of $H$, if $H\subset L$;
\item $H$ is a \textbf{child} of $L$ and $L$ is the \textbf{parent} of $H$, if $H\subset L$ and $\ell (H) = \frac{1}{2} \ell (L)$;
\item $H$ and $L$ are \textbf{neighbors} if $\frac{1}{2} \ell (L) \leq \ell (H) \leq \ell (L)$ and $H\cap L \neq \emptyset$. 
\end{itemize}
\end{definition}

Note, in particular, the following elementary consequence of the subdivision of $\mathscr{C}$:

\begin{lemma}\label{l:child_parent}
Let $H$ be a boundary cube. Then any ancestor $L$ and any neighbor $L$ with $\ell (L) = 2 \ell (H)$ is necessarily a boundary cube. In particular: the descendant of a non-boundary cube is a non-boundary cube.
\end{lemma}
\begin{proof} For the case of ancestors it suffices to prove that if $L$ is the parent of a boundary cube $H$, then $L$ is a boundary cube. Since the parent of $H$ is a neighbor of $H$ with $\ell (L) = 2\ell (H)$, we only need to show the second part of the statement of the lemma. The latter is a simple consequence of the following chain of inequalities:
$$\begin{aligned}
\dist (c_L, \gamma) &\leq \dist (c_H, \gamma) + |c_H-c_L| = \dist (c_H, \gamma) + 3 \sqrt{m} \ell (H) \\
&< 64 r_H + 3 \frac{r_H}{M_0} \leq \left(64 + 3 M_0^{-1}\right) \frac{r_L}{2} \leq \frac{67}{2} r_L < 64 r_L .
\end{aligned}$$
\end{proof}

\begin{definition}[Satellite balls]\label{defi:satellite-balls}
Following the notations above, we set:
\begin{enumerate}[\upshape (i)]
\item If $L\in\sC^\natural$, then we define the \textbf{non-boundary satellite ball $\bB_L := \ball{p_L}{64r_L}$} where $p_L\in\spt(T)$ such that $\bp_{\pi_0}(p_L)=c_L$, such \(p_L\) is a priori not unique, and $\pi_L$ is a plane which minimizes the excess in $\bB_L$, namely $\bE (T, \bB_L)= \bE (T, \bB_L, \pi_L)$,

\item If $L\in \mathscr{C}^\flat$, then we define the \textbf{boundary satellite ball $\bB^\flat_L := \ball{p_L^\flat}{2^7 64r_L}$} where $p_L^\flat\in\Gamma$ is such that $|\bp_{\pi_0}(p_L^\flat)-c_L|=\dist(c_L,\gamma)$. Note that in this case the point $p^\flat_L$ is uniquely determined because $\Gamma$ is regular and $\bA$ is assumed to be sufficiently small. Likewise $\pi_L$ is a plane which minimizes the excess $\bE^\flat$, namely such that $\bE^\flat (T,\bB^\flat_L) = \bE (T, \bB^\flat_L, \pi_L)$ and $T_{p^\flat_L} \Gamma \subset \pi_L$.
\end{enumerate}
\end{definition}

A simple corollary of Theorem \ref{thm:improved-excess-decay} and Corollary \ref{cor:height-tilted-planes} is the following lemma.

\begin{lemma}\label{lem:raffinamento}
Let $T$ and $\Gamma$ be as in Assumption \ref{assumptions}. Then there
is a positive dimensional constant $C (m,n)$ such that, if  the starting size of the Whitney decomposition satisfies $2^{N_0} \geq C (m,n) M_0$, then the satellite balls $\bB^\flat_L$ and $\bB_L$ are all contained in $\bB_2$. Moreover,  there exists $\varepsilon_1$ such that, for any choice of $M_0 , \alpha_\be>0$ and $\alpha_\bh<\frac 12$, if

\begin{equation}\label{e:eps_condition}
\bE^\flat(T,\bB_2)+ \|\psi\|_{C^{3,a_0}}^2<\varepsilon_1\,,
\end{equation}

then for every cube $L\in\sC^\flat$ we have 
\begin{equation}
    \bE^\flat(T, \bB^\flat_L) \le C_0 \varepsilon_1 r_L^{2-2\alpha_\be} , \label{e:fa}\\
\end{equation}
\begin{equation}
    \bh(T, \bB_L^\flat,\pi_L)\le C_0 \varepsilon_1^{\sfrac{1}{4}} r_L^{1+\alpha_\bh} , \label{e:fb}\\
\end{equation}
\begin{equation}
    |\pi_L - \pi_0| \le C_0 \varepsilon_1^{\sfrac{1}{2}} ,\label{e:fc}\\
\end{equation}
\begin{equation}
    \abs{\pi_L - \pi(p_L^\flat)} \le C_0  \varepsilon_1^{\sfrac{1}{2}} r_L^{1-\alpha_\be}, \label{e:fd}
\end{equation}
where, \(\pi(p_L^\flat)\) is the $m$-dimensional tangent plane supporting the tangent cone to $T$ at $p_L^\flat$ and $C_0$ depends only upon $\alpha_\be$, $\alpha_\bh$, $m$ and $n$. 
\end{lemma}
\begin{proof} 
Take $x\in\mathbf{B}_L$, using the height bound in (iii) of Assumption \ref{assump:superlinear-lip-approx}  we obtain that 

$$\abs{x}\leq 64r_L+\abs{p_L} \leq 64\sqrt{m}M_02^{-N_0} + \abs{c_L} + C\varepsilon_0^{\sfrac{1}{2}} \abs{p_L}, $$

recalling that $c_L\in\ball{0}{\sfrac{3}{2}}$, possibly choosing $\varepsilon_0$ small enough and taking the constants $C(m,n), N_0$ big enough, we certainly obtain that $\mathbf{B}_L\subset\ball{0}{2}$. The proof that $\mathbf{B}^\flat_L\subset\ball{0}{2}$ is analogous with the exception that we might multiply the dimensional constant by $2^7$. Inequality \eqref{e:fc} is a direct consequence of \eqref{eq:thm-height-tilted:dist-optimal-planes}. In this proof we will use the improved excess decay, i.e., Theorem \ref{thm:improved-excess-decay}, with $q=p^\flat_L, p=0, r=1, \rho = 2^764r_L$. To prove estimate \eqref{e:fa}, we do as follows

$$ \bE^\flat (T, \bB_L^\flat) = \bE^\flat (T, \ball{p_L^\flat}{2^7 64r_L}) \overset{\eqref{eq:decay-improved-excess}}{\leq} C (2^7 64r_L)^{2-2\alpha_\be}  \bE^\flat (T, \ball{0}{2}) + C (2^7 64r_L)^{2-2\alpha_\be}\bA^2  ,$$

and thus \eqref{e:eps_condition} concludes the proof of \eqref{e:fa}. With the same argument we prove \eqref{e:fb} using in this time the height bound given by the excess decay, i.e., \eqref{eq:thm-improved-excess-decay:height-bound}. It remains to prove \eqref{e:fd}, by the monotonicity formula, and recalling that $\Theta (T, p^\flat_L) = Q-\frac{Q^\star}{2} \geq \frac{3}{2}$, we have 
\[
\|T\| (\bB^\flat_L) \geq \omega_m (2^7 64 r_L)^m. 
\]
Therefore, by the improved excess decay, we obtain
\[
\bE (T, \bB^\flat_L, \pi_L) \leq \bE (T, \bB^\flat_L, \pi (p^\flat_L)) \overset{\eqref{eq:decay-improved-excess},\eqref{e:eps_condition}}{\leq} C_0 \varepsilon_1 r_L^{2-2\alpha_\be}\, .
\]
Hence
\[
|\pi (p^\flat_L) - \pi_L|^2 \leq C_0 \big (\bE (T, \bB^\flat_L, \pi_L) + \bE (T, \bB^\flat_L, \pi (p^\flat_L))\big) \leq C_0 \varepsilon_1 r_L^{2-2\alpha_\be}\, .
\]
\end{proof}

\subsection{Stopping conditions of the Whitney decomposition}

We will now start to refine our Whitney decomposition putting into account the properties of small excess and height bound of the current, in the sense of Lemma \ref{lem:raffinamento}, to then obtain further stronger information about the current on each cube of the decomposition. To this end let \(C_\be, C_\bh\) be two large positive constants that will be fixed later.  We take a cube $L\in\sC_{N_0}$ and we {\bf do not} subdivide it if either the excess ''is too big'' or the current ''is too high'', precisely if it belongs to one of the following sets:

\begin{enumerate}[\upshape (1)]
\item $\sW_{N_0}^{\be}:=\{L\in\sC_{N_0}^\natural:\, \bE(T, \bB_L)> C_\be \varepsilon_1 \ell(L)^{2-\alpha_\be}\}$;
\item $\sW_{N_0}^{\bh}:=\{L\in\sC_{N_0}^\natural:\, \bh(T, \bB_L,\pi_L)> C_\bh\varepsilon_1^{\sfrac{1}{2m}} \ell(L)^{1+\alpha_\bh}\}$.
\end{enumerate}

We then define 
\[
\sS_{N_0}:=\sC_{N_0}\setminus \left(\sW_{N_0}^{\be}\cup\sW_{N_0}^{\bh}\right)\,.
\]

The cubes in $\sS_{N_0}$ will be subdivided in their children, it means that we are subdividing the cube whenever the current is well behaved in it. In what follows we aim to show that the current is well behaved in the whole ball, it means that we will ensure that $\sW_{N_0} := \sW_{N_0}^{\be}\cup\sW_{N_0}^{\bh} = \emptyset$, and therefore $\mathscr{C}_{N_0} = \sS_{N_0}$, by choosing $C_\be$ and $C_\bh$ large enough, depending only upon $\alpha_\bh, \alpha_\be, M_0$ and $N_0$, see Proposition \ref{pr:tilting_cm} below. 

We next describe the refining procedure assuming inductively that for a certain step $j\geq N_0+1$ we have defined
the families $\sW_{j-1}$ and $\sS_{j-1}$. 
In particular we consider all the cubes $L$ in $\mathscr{C}_j$ which are contained in some element of $\sS_{j-1}$. Among them we select and set aside in the classes $\sW_j:=\sW_j^{\be}\cup\sW_j^{\bh}\cup\sW_j^{\mathbf{n}}$ those cubes where the following stopping criteria are met:

\begin{enumerate}[\upshape (1)]
\item $\sW_j^{\be}:=\{L \mbox{ child of } K \in\sS_{j-1}:\, \bE(T, \bB_L)> C_\be\varepsilon_1 \ell(L)^{2-2 \alpha_\be}\}$,

\item $\sW_j^{\bh}:=\{L \mbox{ child of } K\in \sS_{j-1}:\, L\not\in \sW_j^\be\mbox{ and } \bh(T, \bB_L,\pi_L)> C_\bh\varepsilon_1^{\sfrac{1}{2m}} \ell(L)^{1+\alpha_\bh}\}$,
 
\item $\sW_j^{\mathbf{n}}:=\{L\mbox{ child of } K\in \sS_{j-1}:\, L\not\in \sW_j^\be\cup \sW_j^\bh \mbox{ but } \exists L'\in \sW_{j-1} \text{ with }L\cap L' \neq \emptyset\}$.
\end{enumerate}

We keep refining the decomposition in the set  
\[
\sS_j:=\left\{ L\in \mathscr{C}_j \mbox{ child of } K\in \sS_{j-1}\right\} \setminus \sW_j \, .
\] 
Observe that it might happen that the child of a cube in $\sS_{j-1}$ does not intersect $\basehalfball{+}{0}{3/2}$: in that case, according to our definition, the cube does not belong to $\sS_j$ neither to $\sW_j$: it is simply discarded. As already mentioned, we use the notation $\sS_j^\flat$ and $\sS_j^\natural$ respectively for $\sS_j\cap \sC^\flat$ and $\sS_j\cap \sC^\natural$.
Furthermore we set
\begin{align*}
\sW  :=\bigcup_{j\ge N_0} \sW_j,\quad 
\sS :=\bigcup_{j\ge N_0} \sS_j,\quad 
\mathbf S^+  :=\bigcap_{j\ge N_0}\Big( \bigcup_{L\in \sS_j} L\Big)=\basehalfball{+}{0}{3/2}\setminus \bigcup_{H\in \sW} H\, .
\end{align*}

Note, in particular, that the refinement of boundary cubes can {\em never} be stopped because of the conditions (1) and (2), as we state in the following.

\begin{lemma}\label{cor:no_stop_b}
$\mathscr{C}^\flat_j \cap \sW_i = \emptyset$ for every $i, j\geq N_0$ and in particular $\gamma \cap \basehalfball{+}{0}{3/2} \subset \mathbf S^+$. Thus boundary cubes always belong to $\sS$.
\end{lemma}

\begin{proof}
Assume there is a boundary cube in $\sW_i$, then let $L$ be a boundary cube in $\sW_i$ with largest side length. The latter must then belong to $\sW^{\mathbf{n}}_i$ because Lemma \ref{lem:raffinamento} excludes the possibility of $L$ to belong to either $\sW^{\mathbf{\be}}_j$ or $\sW^{\mathbf{\bh}}_j$. However, by definition of the family, this would imply the existence of a neighbor $L' \in \sW_i$ with $\ell (L') = 2\ell (L)$. By Lemma \ref{l:child_parent},
$L'$ would be a boundary cube in $\sW$, contradicting the maximality of $L$. 
\end{proof}

Clearly, descendants of boundary cubes might become non-boundary cubes and so their refining cubes can be stopped. We finally set $\sW_j := \sW_j^{\be}\cup\sW_j^{\bh}\cup\sW_j^{\mathbf{n}}$. From now on we specify a set of assumptions on the various choices of the constants involved in the construction.
\begin{assumption}\label{ass:cm}
$T$ and $\Gamma$ are as in Assumptions \ref{assumptions} and we also assume that

\begin{enumerate}[\upshape (a)]
\item $\alpha_\bh$ is smaller than $\frac{1}{2m}$ and $\alpha_\be$ is positive but small, depending only on $\alpha_{\bh}$;

\item $M_0$ is larger than a suitable constant, depending only upon $\alpha_\be$;

\item $2^{N_0} \geq C (m,n, M_0)$, in particular it satisfies the condition of Lemma \ref{lem:raffinamento};

\item $C_\be$ is sufficiently large depending upon $\alpha_\be$, $\alpha_\bh$, $M_0$ and $N_0$;

\item $C_\bh$ is sufficiently large depending upon $\alpha_\be, \alpha_\bh, M_0, N_0$ and $C_\be$;

\item \eqref{e:eps_condition} holds with an $\varepsilon_1$ sufficiently small depending upon all
the other parameters.
\end{enumerate}

Finally, there is an exponent $\alpha_{\mathbf{L}}$, which depends only on $m,n, Q^\star$ and $Q$ and which is independent of all the other parameters, in terms of which several important estimates in Theorem \ref{thm:center-manifold-app} will be stated.
\end{assumption}

We are ensuring that there is a nonempty set of parameters satisfying all the requirements, since the parameters are chosen following a precise hierarchy. The hierarchy is consistent with that of \cite{DS4}, the reader could compare Assumption \ref{ass:cm} with \cite[Assumption 1.8 and 1.9]{DS4}.

\subsection{Tilting optimal planes and \texorpdfstring{$L$}{Lg}-interpolating functions}\label{s:interpolating} 

In this section, we will define the interpolating functions which will give rise to the function $\boldsymbol{\varphi}^+$ whose graphs will define the center manifold $\mathcal{M}^+$. In order to begin with the construction of these objects, we shall notice that an important fact is that up to now, in the construction of the decomposition, we have nice local information about the current, i.e., inside each square of the decomposition we can do a good analysis, however, we do not know how to work among those cubes, i.e., how to compare quantities as the excess and height between two different cubes of the decomposition. To that end, we enunciate the following crucial result which is the analog of \cite[Proposition 8.24]{DDHM} which is stated for $Q^\star =1$, we mention that the proof of this proposition readily works for currents with boundary multiplicity equal to $Q^\star\geq 1$.

\begin{proposition}[Tilting and optimal planes, Proposition 8.24, \cite{DDHM}]\label{pr:tilting_cm}
Under the Assumptions \ref{assump:superlinear-lip-approx} and \ref{ass:cm}, we have $\sW_{N_0} = \emptyset$.  Then the following estimates hold for any couple of neighbors $H,L\in \sS \cup \sW$ and for every $H, L \in \sS\cup \sW$ with $H$ descendant of $L$:
\begin{enumerate}
\item[\upshape (a)]  denoting by $\pi_H, \pi_L$ the optimal planes for the excess in $\bB_H$ and $\bB_L$, respectively, we have 
$$|\pi_H - \pi_L|\leq \bar C \eps_1^{\sfrac{1}{2}} \ell (L)^{1-\alpha_\be},\qquad |\pi_H-\pi_0| \leq \bar C\varepsilon_1^{\sfrac{1}{2}},$$

\item[\upshape (b)$^\natural$] $\bh (T, \bC_{48 r_H} (p_H, \pi_0)) \leq C \varepsilon_1^{\sfrac{1}{2m}} \ell (H)$ and $\spt (T)\cap \bC_{48r_H} (p_H, \pi_0) \subset \bB_H$ if $H\in \mathscr{C}^\natural$,

\item[\upshape (b)$^\flat$] $\bh (T, \bC_{2^7 48 r_H} (p^\flat_H, \pi_0)) \leq C \varepsilon_1^{\sfrac{1}{4}} \ell (H)$ and $\spt (T)\cap \bC_{2^7 48 r_H} (p^\flat_H, \pi_0) \subset \bB^\flat_H$ if $H\in \mathscr{C}^\flat$,

\item[\upshape (c)$^{\natural}$] $\bh (T, \bC_{36 r_L} (p_L, \pi_H))\leq C \varepsilon_1^{\sfrac{1}{2m}} \ell (L)^{1+\alpha_\bh}$ and $\spt (T) \cap \bC_{36 r_L} (p_L, \pi_H) \subset \bB_L$ if $H, L \in \mathscr{C}^\natural$,

\item[\upshape (c)$^{\flat}$] $\bh (T, \bC_{2^7 36r_L} (p^\flat_L, \pi_H))\leq C \varepsilon_1^{\sfrac{1}{4}} \ell (L)^{1+\alpha_\bh}$ and $\spt (T) \cap \bC_{2^7 36 r_L} (p^\flat_L, \pi_H)) \subset \bB^\flat_L$ if $L \in \mathscr{C}^\flat$,
\end{enumerate}
where $\bar C = \bar C (\alpha_{\be}, \alpha_{\bh}, M_0, N_0, C_\be)>0$ and $C = C (\alpha_{\be}, \alpha_{\bh}, M_0, N_0, C_\be, C_\bh)>0$.
\end{proposition}

We now state the following results which is the analog of \cite[Proposition 8.7]{DDHM} will allow us to locally approximate the current by $(Q- \frac{Q^\star}{2})$-Lipschitz maps in the sense of Theorem \ref{thm:lip-approx-superlinear-decay}. We also noticed that the proof given in \cite[Proposition 8.7]{DDHM} for $Q^\star =1$ readily works for currents with boundary multiplicity equal to $Q^\star\geq 1$.

\begin{proposition}[Proposition 8.7, \cite{DDHM}]\label{prop:yes_we_can_more}\label{prop:yes_we_can}
Under the Assumptions \ref{assump:superlinear-lip-approx} and \ref{ass:cm} the following holds for every
couple of neighbors $H, L \in \sS \cup \sW$ and any $H, L\in \sS \cup \sW$ with $H$ descendant of $L$:
\begin{align*}
&\spt(T)\cap\bC_{36r_L}(p_L,\pi_H) \subset\bB_L
&\mbox{when $L\in\mathscr{C}^\natural$,}\\
&\spt(T)\cap\bC_{2^7 36r_L}(p^\flat_L,\pi_H) \subset\bB_L^\flat  &\mbox{when $L\in\sC^\flat$,} 
\end{align*}
and the current $T$ satisfies the assumptions of \cite[Theorem 2.4]{DS3} in the cylinder $\bC_{36 r_L} (p_L, \pi_H)$ (resp. of Theorem \ref{thm:lip-approx-superlinear-decay} in the cylinder $\bC_{2^7 36r_L}(p^\flat_L,\pi_H)$). 
\end{proposition}

We will now construct the ``interpolating functions'' $g_L$ for each cube $L$. To begin with the construction of this interpolation, we approximate the current $T$ by $(Q-\frac{Q^\star}{2})$-Lipschitz functions (Proposition \ref{prop:yes_we_can_more}) that will determine the boundary condition of an elliptic system which comes from the linearization of the mean curvature condition for minimal surfaces. The solution of this elliptic system will be further represented by a function $g_L$ defined in the tilted ball contained in $\pi_0$ taking values in $\pi_0^\perp$, i.e., we changed to our new reference coordinate system. Since the construction will be local, i.e., in each cube, over the set $\basehalfball{+}{0}{3/2}\setminus \mathbf{S}^+$ we will patch every $g_L$ together with a partition of the unity to obtain the function $\boldsymbol{\varphi}^+$, whose graph will be the center manifold, defined in the whole ball. So we need to define $\boldsymbol{\varphi}^+$ over $\mathbf{S}^+$ as well, at the end we introduce all the machinery needed for all cubes in $\sS\cup \sW$.

\begin{definition}[$\pi_L$-approximations]\label{def:pi-approximation}
Under the Assumptions \ref{assump:superlinear-lip-approx} and \ref{ass:cm}, we set
\begin{enumerate}[\upshape (i)]
    \item If $L\in\sS_j^\flat$ for some $j$, take the Lipschitz approximation $(f^+_L, f^-_L)$ in the cylinder $\cyltilted{p_L^\flat}{2^7 9r_L}{\pi_L}$ given by Proposition \ref{prop:yes_we_can_more}, we call \textbf{$(f^+_L, f^-_L)$ a $\pi_L$-approximation of $T$ in the cylinder $\cyltilted{p_L^\flat}{2^7 9r_L}{\pi_L}$}.
    
    \item If $L\in\sS_j^\natural\cup\sW_j$ for some $j$, we take the Lipschitz approximation $f_L$ in the cylinder $\cyltilted{p_L}{9r_L}{\pi_L}$ given by Proposition \ref{prop:yes_we_can_more}, we call \textbf{$f_L$ a $\pi_L$-approximation of $T$ in the cylinder $\cyltilted{p_L}{9r_L}{\pi_L}$}.
    
\end{enumerate}
\end{definition}

\begin{definition}[$L$-tilted harmonic interpolations]\label{def:harmonic-interpolation}
Under the Assumptions \ref{assump:superlinear-lip-approx} and \ref{ass:cm}, we define
\begin{enumerate}[\upshape (i)]
    \item if $L$ is a nonboundary cube, let $h_L:\mathrm{B}_{5r_L}(p_L,\pi_L)\to\R$ to be an harmonic function with boundary condition $h_L|_{\partial \mathrm{B}_{5r_L}(p_L,\pi_L)} = \boldsymbol{\eta}\circ f_L|_{\partial \mathrm{B}_{5r_L}(p_L,\pi_L)}$,
    
    \item if $L$ is a boundary cube, let $h_L:\mathrm{B}_{2^7 5r_L}(p_L^\flat,\pi_L)\to\R$ to be an harmonic function with boundary condition $h_L|_{\partial \mathrm{B}_{2^7 5r_L}(p_L^\flat,\pi_L)} = \boldsymbol{\eta}\circ f_L^+|_{\partial \mathrm{B}_{2^7 5r_L}(p_L^\flat,\pi_L)}$. 
\end{enumerate}
We call \textbf{$h_L$ $L$-tilted harmonic interpolating function}.
\end{definition}

We now are ready to define the final function, $g_L$, on our ``reference coordinate system'', i.e., the domain of $g_L$ is contained in $\pi_0$ and $g_L$ takes values in $\pi_0^\perp$, with the property that its graph
coincides with a suitable portion of the graph of $h_L$. The function $g_L$ is furnished by \cite[Lemma B.1]{DS4} which we state below.

\begin{proposition}[$L$-interpolating functions]\label{def:interpolating-cm}
Under the assumptions of Proposition \ref{prop:yes_we_can}, we have
\begin{enumerate}[\upshape (i)]
    \item If $L$ is a boundary cube, the function $h_L$ is Lipschitz on $\mathrm{B}^\pm_{ 2^7 9r_L/2} (p^\flat_L, \pi_L)$ and we can define a function  $g_L:\mathrm{B}^+_{2^7 4r_L}(p^\flat_L,\pi_0)\to \pi_0^\perp$ such that $\bG_{g_L}=\bG_{h_L}\res \mathrm{B}^+_{2^7 4r_L} (p^\flat_L, \pi_0) \times \R^{n}$,
    
    \item If $L$ is a non-boundary cube, the function $h_L$ is Lipschitz on $\mathrm{B}_{9r_L/2} (p_L, \pi_L)$ and we can define a function $g_L:\mathrm{B}_{4r_L}(p_L,\pi_0)\to\pi_0^\perp $ such that $\bG_{g_L}=\bG_{h_L}\res \bC_{4r_L} (p_L, \pi_0)\big)$.
    \end{enumerate}
    
The functions $g_L$ is called \textbf{$L$-interpolating function}.
\end{proposition}

\subsection{Glueing \texorpdfstring{$L$}{Lg}-interpolations}\label{ss:glued_interpolations}

We now define another set of cubes, the Whitney cubes at the step $j$, which will be similar to what we have constructed until now but we are including all the ancestors with respect to step $j$ in the same family as follows
$$\mathscr P_j:=\sS_j\cup\bigcup_{i=N_0+1}^j \sW_i .$$
Note that $\mathscr P_j$ is a ``Whitney family of dyadic cubes'' in the sense that if $K, L \in \mathscr{P}_j$ and $K\cap L\neq\emptyset$, then $\frac{1}{2} \ell (L) \leq \ell (K) \leq 2 \ell (L)$. We fix a partition of unity $\vartheta$ satisfying
\begin{equation*}
\vartheta\in C^\infty_c \left(\left[-\frac{17}{16}, \frac{17}{16}\right]^m, \left[0,1\right]\right), \ \vartheta |_{[-1,1]^m} \equiv 1, \ \text{and, for each cube} \ L, \ \tilde{\vartheta}_L(y)
 := \vartheta \left(\frac{y-c (L)}{\ell (L)}\right).
\end{equation*}
We thus set a partition of unity of $\basehalfball{+}{0}{3/2}$ defined as
\begin{equation*}
\vartheta_L:\basehalfball{+}{0}{\sfrac{3}{2}}\to\R, \ \vartheta_L(y) := \frac{\tilde{\vartheta}_L(y)}{\sum_{H \in \mathscr P_j} \tilde{\vartheta}_H(y)}.
\end{equation*} 

\begin{definition}[Glued interpolation at the step $j$]\label{def:glued_interpolations}
We set $\varphi_j:=\sum_{L\in\mathscr P_j}\vartheta_L  g_L$, this maps $\varphi_j$ are called \textbf{{glued interpolation at the step $j$.}}
\end{definition}

\subsection{Existence of a \texorpdfstring{$C^{3,\kappa}$}{Lg}-center manifold}

We are now ready to state the main theorem regarding the construction of the center manifolds needed in this paper, i.e., Theorem \ref{thm:center_manifold}, which ensures that $(\varphi_j)_j$ is a sequence that converges to a $C^{3,\kappa}$ map, $\kappa>0$, whose graph will be called the center manifold. This limit map has good properties as the smallness of the $C^{3,\kappa}$ norm, which is bounded by $\varepsilon_1$. After the main theorem, we will start the construction of the normal approximation in the sense of Theorem \ref{thm:lip-approx-superlinear-decay} but now the approximations will be defined on the center manifold and will take values on the normal bundle of the center manifold. The normal approximations will enjoy some good estimates relying in the estimates of Theorem \ref{thm:lip-approx-superlinear-decay} and the estimate obtained in the construction of our Whitney decomposition. The main theorem of this section is stated below and is the version adapted to our setting of \cite[Theorem 8.13]{DDHM}. 

\begin{theorem}[Theorem 8.13 ,\cite{DDHM}]\label{thm:center_manifold}
Under Assumptions \ref{assump:superlinear-lip-approx} and \ref{ass:cm}, there is a $\kappa := \kappa(\alpha_{\be}, \alpha_{\bh})>0$, such that

\begin{enumerate}[\upshape (i)]
\item $\varphi_j\in C^{3,\kappa}$, moreover $\|\varphi_j\|_{3, \kappa, \basehalfball{+}{0}{3/2}} \leq C \varepsilon_1^{\sfrac{1}{2}}$, for some $C:=C(\alpha_\be, \alpha_\bh, M_0, C_{\be}, C_{\bh})>0$,

\item If \(i\le j\), $L\in \sW_{i-1}$ and $H$ is a cube concentric to $L$ with $\ell (H) = \frac{9}{8} \ell (L)$, then $\varphi_j = \varphi_i$ on $H$,

\item $\varphi_j$ converges in $C^3$ to a map $\mathbf{\varphi}^+: \basehalfball{+}{0}{3/2}\to \R^n$, whose graph is a $C^{3,\kappa}$ submanifold $\mathcal M^+$, which will be called the \textup{\textbf{right center manifold}};

\item $\mathbf{\varphi}^+ = \psi$ on $\gamma\cap B_{3/2}$, i.e., $\partial \mathcal M^+ \cap \cyl{0}{3/2} = \Gamma \cap \cyl{0}{3/2}$;

\item For any $q\in \partial \mathcal M^+ \cap \cyl{0}{3/2}$, we have $T_q \mathcal M^+ = \pi (q)$ where $\pi(q)$ is the support of the unique tangent cone to $T$ at $q$.
\end{enumerate}
\end{theorem}

We will omit the proof of Theorem \ref{thm:center_manifold} since it goes along the same lines of \cite[Theorem 8.13]{DDHM}. In fact a careful inspection of its proof (for $Q^{\star}=1$) will reveal that $Q^{\star}$ comes into play to assure the compatibility of the traces, hence it also holds mutatis mutandis when $Q^\star\ge1$.

\begin{remark}
The construction of $\mathcal M^+$ made in Theorem \ref{thm:center_manifold} is based on the decomposition of $\basehalfball{+}{0}{3/2}$. Under Assumption \ref{ass:cm}, the same construction can be made for $\basehalfball{-}{0}{3/2}$ and gives a $C^{3,\kappa}$ map $\mathbf{\varphi}^-: \basehalfball{-}{0}{3/2}\to\R^n$ which agrees with $\psi$ on $\gamma\cap B_{3/2}$. The graph of $\mathbf{\varphi}^-$ is a $C^{3,\kappa}$ submanifold $\mathcal M^-$, which will be called the \textbf{{left center manifold}}. Its boundary in the cylinder $\cyl{0}{3/2}$, namely $\cyl{0}{3/2}\cap\partial\mathcal M^-$, coincides, in a set-theoretical sense, with $\cyl{0}{3/2}\cap\partial\mathcal M^+$, but it has opposite orientation, and moreover its tangent plane $T_q\mathcal M^-$ coincides with $\pi(q)$ for every point $q\in\cyl{0}{3/2}\cap\partial\mathcal{M}^-$.
\end{remark}

In particular, the union $\mathcal M:= \mathcal M^+\cup\mathcal M^-$ of the two submanifolds  is a $C^{1,1}$ submanifold in $\cyl{0}{3/2}$ without boundary in $\cyl{0}{3/2}$, which will be called the \textbf{center manifold}. Moreover, we will often state properties of the center manifold related to cubes $L$ in one of the collections $\sW_j$ described above. Therefore, we will denote by $\sW^+$ the union of all $\sW_j$ and by $\sW^-$ the union of the corresponding classes of cubes which lead to the left center manifold $\mathcal{M}^-$. We emphasize again that so far we can only conclude the $C^{1,1}$ regularity of $\mathcal M$, because we do not know that the traces of the second derivatives of $\mathbf{\varphi}^+$ and $\mathbf{\varphi}^-$ coincide on $\gamma$.

\begin{definition}
Let us define the graph parametrization map of $\mathcal M^+$ as $\mathbf{\Phi}^+(x):=(x,\mathbf{\varphi}^+(x))$. We will call \textbf{{right  contact set}} the subset $\mathbf{K}^+ := \mathbf{\Phi}^+(\mathbf S^+)$. For every cube $L\in \sW^+$ we associate a \textbf{ {Whitney region $\mathcal L$ on $\mathcal M^+$}} as follows: 
\begin{itemize}
\item $\mathcal{L}:= \mathbf{ \Phi}^+(H\cap \baseball{0}{1})$ where $H$ is the cube concentric to $L$ such that $\ell(H)=\frac{17}{16}\ell(L)$.
\end{itemize}
Analogously we define the map $\mathbf \Phi^-$, the \textbf{left contact set $\mathbf{K}^-$} and the \textbf{Whitney regions on the left center manifold $\mathcal{M}^-$}. 
\end{definition}

To keep the text flow, we postpone the proof of the Theorem \ref{thm:center_manifold} to the last part of this section.

\subsection{The \texorpdfstring{$\mathcal{M}$}{Lg}-Lipschitz approximations defined on the center manifolds}

Since the two portions $\mathcal{M}^-$ and $\mathcal{M}^+$ are $C^{3, \kappa}$ and they join with $C^{1,1}$ regularity along $\Gamma$, in a sufficiently small normal neighborhood of $\mathcal{M}$ there is a well defined orthogonal projection $\mathbf{p}$ onto $\mathcal{M}$. The thickness of the tubular neighborhood is inversely proportional to the $C^{1,1}$-norm of $\mathbf{ \varphi}^\pm$ and hence, for $\varepsilon_1$ sufficiently small, we can assume that the thickness is $2$ which leads to the next assumption.

\begin{assumption}\label{assump:center-manifold-app}
Under Assumptions \ref{assump:superlinear-lip-approx} and \ref{ass:cm}. We let $\mathcal{M}:= \mathcal{M}^+\cup\mathcal{M}^-$ and $\varepsilon_1$ be sufficiently small so that, if
\begin{align*}
\mathbf{U} := \{q\in \R^{m+n} :\,\exists! q' = \mathbf{p} (q) \in \mathcal{M} \mbox{ s.t. $|q-q'|<1$} \mbox{ and $q-q'\in T_{q'}^\perp\mathcal{M}$}\}\, ,\label{e:tub_neigh}
\end{align*}
where $T_{q'}^\perp\mathcal{M}:=(T_{q'}\mathcal{M})^\perp$, then the map $\mathbf{p}$ extends to a Lipschitz map to the closure $\overline{\mathbf{U}}$  which is  $C^{2, \kappa}$ on \(\mathbf{U}\setminus \mathbf{p}^{-1} (\Gamma)\) and
\[
\mathbf{p}^{-1} (q') = q' + \overline{\mathrm{B}_1 (0, (T_{q'} \mathcal{M})^\perp)}\mbox{ for all $q'\in \overline{\mathcal{M}}$.}
\]
\end{assumption}

As highlighted before, we construct the center manifold $\mathcal{M}$ and also a function defined on $\mathcal{M}$ that approximates, with the desired superlinear exponents for the error, the current $T$ (in the sense of Theorem \ref{thm:lip-approx-superlinear-decay}). This approximation will take values on the normal bundle of $\mathcal{M}$, we precisely define this type of approximations. Firstly, we shall define the space of $Q$-tuples on a manifold analogously to \cite{DS5} which follows the definition for the Euclidean spaces in \cite{DS1}.
\begin{definition}
Let $M$ be an $m$-dimensional manifold, and, for each $P\in M$, we denote $\a{P}$ the current with support equal to $P$, i.e., the current associated with the Dirac measure concentrated in $P$. Then we define the space of unordered $Q$-tuples in $M$, for any $Q\in\N, Q\geq 1$, as follows
$$ \mathcal{A}_Q(M) := \left \{ \sum_{i=1}^Q\a{P_i}: P_i\in M \text{ for every } i\in \{1,\dots, Q\} \right \}. $$
\end{definition}

\begin{definition}
Let $\mathcal{M}$ be the center manifold as in Theorem \ref{thm:center_manifold} without loss of generality we can assume that we are under Assumption \ref{assump:center-manifold-app}, $Q^+ = Q$ and $Q^- = Q-Q^\star$. We say that $(\mathcal{K}, F^+, F^-)$ is an \textbf{$\mathcal{M}$-normal approximation of $T$}, if

\begin{enumerate}[\upshape (i)]
\item there exist Lipschitz functions $\mathcal{N}^+:\mathcal{M}^+\cap \cyl{0}{1}\to\mathcal{A}_Q(T^{\perp}\mathcal{M}^+)$, $\mathcal{N}^+(x)\in\mathcal{A}_Q(T_x^\perp\mathcal{M}^+)$ and $\mathcal{N}^-:\mathcal{M}^-\cap \cyl{0}{1}\to\mathcal{A}_{Q-Q^\star}(T^{\perp}\mathcal{M}^-)$, $\mathcal{N}^-(x)\in\mathcal{A}_{Q-Q^\star}(T_x^\perp\mathcal{M}^-)$, where $T^{\perp}\mathcal{M}^\pm:=\sqcup_{x\in\mathcal{M}^{\pm}}T_x^\perp\mathcal{M}^\pm$ denotes the normal bundle of $\mathcal{M}^\pm$ and is seen as a subset of $\mathbb{R}^{m+n}$, 
$$\begin{array}{lrll}
    \mathcal{N}^\pm: & \mathcal{M}^\pm\cap\cyl{0}{1} & \to & \mathcal{A}_{Q^\pm}(T^\perp\mathcal{M}^\pm)   \\
    & x & \mapsto &\mathcal{N}^\pm(x):=\sum_{i=1}^{Q^\pm}\a{\mathcal{N}_i^\pm(x)},
\end{array}$$ 
where $(\mathcal{N}_i^\pm:\mathcal{M}^+\cap \cyl{0}{1}\to T^\perp\mathcal{M})_{i\in\{1,\dots, Q^\pm\}}$ are measurable sections of the normal bundle, i.e., each $\mathcal{N}_i^\pm$ is a classical $1$-valued measurable function satisfying $\mathcal{N}_i^\pm(x)\in T_x^\perp\mathcal{M}^\pm$. We then define the Lipschitz function given by
$$\begin{array}{lrll}
    F^\pm: & \mathcal{M}^\pm\cap\cyl{0}{1} & \to & \mathcal{A}_{Q^\pm}(T^\perp\mathcal{M}^\pm)   \\
    & x & \mapsto &\left(\mathcal{N}^\pm \oplus \mathrm{id}\right)(x).
\end{array}$$

\item $\mathcal{K} \subset \mathcal{M}$ is closed and $\mathbf{T}_{F^\pm} \res \mathbf{p}^{-1} (\mathcal{K}\cap \mathcal{M}^\pm) = T \res \mathbf{p}^{-1} (\mathcal{K}\cap \mathcal{M}^\pm)$, where $\mathbf T_{F^\pm}:=(F^\pm)_\sharp \a{\mathcal M}$, according to {\cite[Definition 1.3]{DS2}},

\item $\mathbf{K}^+ \cup \mathbf{K}^- \subset \mathcal{K}$, $\mathcal{N}^\pm |_{\mathcal{K}}\equiv 0$, and then $F^+ (x) =  Q \a{x}$ on $\mathcal{K}^+$ and $F^- (x) = (Q-Q^\star) \a{x}$ on $\mathcal{K}^-$.
\end{enumerate}
\end{definition}

Observe that the pairs $(F^+, F^-)$ and $(\mathcal{N}^+, \mathcal{N}^-)$ are intuitively $(Q-\frac{Q^\star}{2})$-valued maps in the spirit of Definition \ref{def:function-interface}. Although this is very intuitive, these functions are defined on manifolds, so, we make the precise definition of it as follows.

\begin{definition}\label{def:function-interface-manifold} 
Firstly, we let $Z$ be an $m$-dimensional manifold and $\Upsilon$ be a $(m-1)$-submanifold of the $m$-manifold $M$ which splits $M$ into the two connected components $M^+$ and $M^-$. Let $\Phi\in W^{\frac12,2}(\Upsilon, \mathcal{A}_{Q^\star}(Z))$, $Q,Q^{\star}\in\N$, $Q\ge Q^{\star}\ge1$. A \textbf{$(Q-\frac{Q^{\star}}2)$-valued function with interface $(\Upsilon,\Phi)$}, consists of a pair $(F^+, F^-)$ satisfying the following properties 
\begin{enumerate}[\upshape (i)]
          \item $F^+\in W^{1,2}(M^+,  \mathcal{A}_{Q}(Z))$, $F^-\in W^{1,2}(M^-,   \mathcal{A}_{Q-Q^\star}(Z))$,
          \item $F^+\restrict\Upsilon=F^-\restrict\Upsilon+\Phi$.
\end{enumerate} 
We define the \textbf{Dirichlet energy of $(F^+,F^-)$} as $\D(F^+,F^-,M):=\D(F^+,M^+)+\D(F^-,M^-)$. Such a pair will be called \textbf{$\D$-minimizing in $M$}, if for all $\left(Q-\frac{Q^{\star}}{2}\right)$-valued function $(G^+,G^-)$ with interface $(\Upsilon, \Phi)$  which agrees with $\left(F^{+}, F^{-}\right)$ outside of a compact set $K \subset\subset M$ satisfies $\D(F^+,F^-,M)\le\D(G^+,G^-, M)$.
\end{definition}

\begin{definition}\label{def:collapses-interface-manifold}
Let $(F^+,F^-)$ be a $\left(Q-\frac{Q^{\star}}{2}\right)$-valued function with interface $(\Upsilon, \Phi)$ and $\Phi = Q^{\star}\a{\hat{\Phi}}$ for the single valued function $\hat{\Phi}\in W^{\frac12,2}(\Upsilon, Z)$. We say that  $(F^+,F^-)$ \textbf{collapses at the interface}, if $F^{+}\restrict\Upsilon = Q\a{\hat{\Phi}}$.
\end{definition}

The following theorem ensures the existence of an $\mathcal{M}$-normal approximation suitable for our purposes, i.e., with the desired exponents at the bound on the Lipschitz constant, the Dirichlet energy and the size of the complement set of $\mathcal{K}$. It is the same as \cite[Theorem 8.19]{DDHM} but of course adapted to our context where $Q^\star$ is taken any arbitrary positive integer. We omit its proof because goes precisely as in the proof of \cite[Theorem 8.19]{DDHM}.

\begin{theorem}[Local behaviour of the $\mathcal{M}$-normal approximation on Whitney regions. Theorem 8.19 of \cite{DDHM}]\label{thm:center-manifold-app}
Under Assumption \ref{assump:center-manifold-app} there is a constant $\alpha_\bL:=\alpha_\bL(m, n, Q^\star, Q)>0$ such that
there exists an $\mathcal{M}$-normal approximation $(\mathcal{K}, F^+, F^-)$ satisfying the following estimates on any Whitney region $\mathcal{L}\subset \mathcal{M}$ associated to a cube $L \in \sW^+\cup \sW^-$:
\begin{align}
\Lip ({\mathcal{N}}^\pm|_{\mathcal{L}}) & \leq C\varepsilon_1^{\alpha_\bL} \ell (L)^{\alpha_\bL} \label{e:cm_app1},\\
\|{\mathcal{N}}^\pm|_{\mathcal{L}}\|_0 &\leq C \varepsilon_1^{\sfrac{1}{2m}} \ell (L)^{1+\alpha_\bh}\label{e:cm_app2},\\
\mathcal{H}^m(\mathcal{L}\setminus \mathcal{K}) + \|\mathbf{T}_F -T\| (\mathbf{p}^{-1} (\mathcal{L})) & \leq C \varepsilon_1^{1+\alpha_\bL} \ell(L)^{m+2+\alpha_\bL}\label{e:cm_app3},\\
\int_{\mathcal{L}} |D{\mathcal{N}}^\pm|^2 & \leq C \varepsilon_1 \ell (L)^{m+2-2\alpha_\be},\label{e:cm_app4}
\end{align}
for a constant $C = C (\alpha_\be,\alpha_\bh,M_0,N_0,C_\be,C_\bh)>0$. Moreover, for any $a>0$ and any Borel $\mathcal{V} \subset \mathcal{L}$,
\begin{equation}\label{e:cm_app5}
\begin{aligned}
\int_{\mathcal{V}} |\eta\circ {\mathcal{N}}^\pm|\mathrm{d}\mathcal{H}^m &\leq C \varepsilon_1 \left(\ell(L)^{m+3+\alpha_\bh/3} + a \ell (L)^{2+\alpha_\bL/2} \mathcal{H}^m(\mathcal{V})\right) \\
&+ \frac{C}{a} \int_{\mathcal{V}} \mathcal{G} ({\mathcal{N}}^\pm, Q \a{\eta\circ {\mathcal{N}}^\pm})^{2+\alpha_\bL}\mathrm{d}\mathcal{H}^m.
\end{aligned}
\end{equation}
\end{theorem}

\newcommand{\sigmaexpcm}{{\alpha_\bL}}

\section{Blowup argument by the frequency function}

In this section we finish the proof of the main theorem of this work, i.e., Theorem \ref{thm:main-theorem}, which is a consequence of Theorem \ref{thm:collapsed-is-regular}, for $m=2$, as noticed in Subsection \ref{two-dim-case}. We use the frequency function in the center manifold $\mathcal{M}$, c.f. \cite[Chapter 9]{DDHM}, originally introduced by Almgren (\cite{Alm}) and, more recently, reformulated and adapted to the boundary case in \cite{DS5,DDHM,DNS}. In Theorem \ref{thm:ff-estimate-current}, we show that the $m$-dimensional area minimizing current has to satisfies at most one of two conditions, where the first one essentially says that $0$ is a two-sided regular point of $T$ and the second one is an estimate with the Almgren's frequency function. Although, there are two alternatives in Theorem \ref{thm:ff-estimate-current}, we use a blowup argument in Theorem \ref{thm:blowup-argument} to show that the second alternative never occurs, thus implying that the only possible situation is $0$ being a two-sided regular boundary point.

\subsection{Almgren's frequency function in \texorpdfstring{$\mathcal{M}$}{Lg}}

In order to define our main quantities, we start with the following lemma that shows that exists a good perturbation function of the distance function on the center manifold .

\begin{lemma}[Lemma 9.2 ,\cite{DDHM}]\label{l:good_vector_field on M}
There exist positive continuous functions $d^\pm:\mathcal M^\pm \to \R_+$ which belong to $C^2(\mathcal M^\pm\setminus\{0\})$ and satisfies the following properties
\begin{enumerate}[\upshape (a)]
\item $d^\pm(x)=\dist_{\mathcal{M^\pm}}(x,0) + O(\dist_{\mathcal{M^\pm}}(x,0)^2)= \abs{x} + O(\abs{x}^2)$,
\item $\abs{\nabla d^\pm(x)} = 1 + O(d^\pm)$, where $\nabla$ is the gradient on the manifold $\mathcal{M}$,
\item $\frac{1}{2}\nabla^2(d^\pm)^2(x)= g + O (d^\pm)$, where \(\nabla^2\) denotes the covariant Hessian on $\mathcal{M}$ (which we regard as a $(0,2)$ tensor) and $g$ is the induced metric on $\mathcal{M}$ as a submanifold of $\mathbb R^{m+n}$,
\item $\nabla d^\pm(p) \in T_p\Gamma$ for all $p \in \Gamma$, i.e.
\begin{equation}\label{e:sign_condition}
\nabla d^\pm \cdot \vec{n}^\pm = 0 \mbox{ on }\Gamma,
\end{equation}
where $\vec{n}^\pm$ denotes the outer unit normal to $\mathcal{M}^\pm$ inside $\mathcal{M}$.
\end{enumerate}
In particular this implies 
\begin{equation}\label{e:hessian d}
\nabla^2d^\pm (x)=\frac{1}{d} \Big(g-\nabla d^\pm(x) \otimes \nabla d^\pm(x)\Big)+O(1)
\end{equation}
and 
\begin{equation}\label{e:laplacian of the good vectorfield} 
	\Delta\, d^\pm = \frac{m-1}{d^\pm} + O(1)\, 
\end{equation}
where $\Delta$ denotes the Laplace-Beltrami operator on $\mathcal{M}$, namely the trace of the Hessian $\nabla^2$. Moreover:
\begin{itemize}
\item[\upshape (S)] All the constants estimating the $O (\cdot)$ error terms in the above estimates can be made smaller than any given $\eta>0$, provided the parameter $\varepsilon_1$ in Assumption \ref{ass:cm} is chosen appropriately small (depending on $\eta$). 
\end{itemize}
\end{lemma}

Let us define three functions that we be used in the definition of the frequency function as follows
\begin{equation}
\phi (t) :=
\left\{
\begin{array}{ll}
1, &\mbox{for $0\leq t \leq \frac{1}{2}$,}\\
2 (1-t), &\mbox{for $\frac{1}{2}\leq t \leq 1$,}\\
0, &\mbox{for $t\geq 1$,}
\end{array}\right.
\end{equation}
\begin{align}
\mathrm{D}_{\phi,d^\pm} (\mathcal{N}^\pm,r) & := \hphantom{-}\int_{\mathcal{M}^\pm} \phi\left(\frac{d^\pm(x)}{r}\right)|D\mathcal{N}^\pm|^2(x)\operatorname{d}\operatorname{Vol}^\pm,\label{e:richiamo_pm}
\\
\mathrm{H}_{\phi,d^\pm} (\mathcal{N}^\pm,r) &:=-\int_{\mathcal{M}^\pm} \phi'\left(\frac{d^\pm(x)}{r}\right) |\nabla d^\pm (x)|^2 \frac{|\mathcal{N}^\pm(x)|^2}{d^\pm(x)}\operatorname{d}\operatorname{Vol}^\pm ,\label{e:richiamo_pm_2}
\end{align}
where $\operatorname{Vol}^\pm$ denotes the standard volume form on $\mathcal{M}^\pm$. 

\begin{definition}
The {\textbf{Almgren's frequency function}} is defined as the ratio
$$\mathrm{I}_{\phi,d^\pm} (\mathcal{N}^\pm,r) :=\frac{r\mathrm{D}_{\phi,d^\pm} (\mathcal{N}^\pm,r)}{\mathrm{H}_{\phi,d^\pm} (\mathcal{N}\pm,r)}.$$
\end{definition}

We also set the notation
$$
\mathcal C^{\pm}:=\left \{y\in \ball{0}{1}:  \bp(y)\in \mathcal M^{\pm}\textrm{ and } |y-\bp(y)|\le \dist(y,\Gamma)^{3/2}  \right \}
$$
for the horned neighborhoods of $\mathcal{M}^\pm$ in which $T$ is supported, compare with  Corollary \ref{cor:height-tilted-planes} and item $($v$)$ of Theorem \ref{thm:center_manifold}. We now state the following theorem that will be crucial in the blowup argument, its proof is the same given in \cite{DDHM} since the authors do not need the multiplicity condition $Q^{\star}=1$.  

\begin{theorem}[Theorem 9.3, \cite{DDHM}]\label{thm:ff-estimate-current}
Let $T$ and $\Gamma$ be as in Assumption~\ref{assump:center-manifold-app}, $Q^+ = Q, Q^- = Q-Q^\star$ and consider $\phi$ and $d^\pm$ as above. Then only one of the following alternatives holds
\begin{enumerate}[\upshape (a)]
\item $T\res \mathcal{C}^\pm = Q^\pm\a{\mathcal{M}^\pm}$ in a neighborhood of $0$,
\item $\lim_{r\to 0} \mathrm{I}_{\phi,d^\pm} (\mathcal{N}^\pm, r)>0.$
\end{enumerate} 
\end{theorem}

\subsection{Blowup argument}\label{blowup}

Letting $Q^+ := Q$ and $Q^- := Q - Q^\star$, we define the multivalued maps with domain and codomain in Euclidean spaces,
$$
N^\pm (x) = \sum_{i=1}^{Q^\pm} \a{(N^i)^\pm (x)},
$$
such selections $\{N^i\}_{i=1}^{Q^\pm}$ are given by the formulas
$$
\begin{array}{crcl}
(N^i)^\pm : & \basehalfball{\pm}{0}{1}\subset\R^m & \to & \R^n  \\
 & x & \mapsto & \bp_{\{0\}\times\R^n} \big((N^i)^\pm (x, \boldsymbol{\varphi}^\pm (x))\big).
\end{array}
$$

Observe that the pair $(N^+, N^-)$ is a $\left(Q-\frac{Q^\star}{2}\right)$-valued function with interface $(\gamma, Q^\star\a{0})$. We now set the following notation for the Dirichlet energy
$$\operatorname{Dir}(r) := \frac{1}{2}\int_{\basehalfball{+}{0}{1}} |DN^+|^2 +\frac{1}{2}\int_{\basehalfball{-}{0}{1}} |DN^-|^2 := \operatorname{Dir}^+(r) +  \operatorname{Dir}^-(r), $$
and the corresponding rescaling of $N^\pm$
$$
N_r^\pm (x) := \sum_i \a{ r^{\sfrac{m}{2}-1} \operatorname{Dir}^\pm(r)^{-\sfrac{1}{2}} (N^i)^\pm (rx)} .$$

Finally, we can state the key result to give our final contradiction argument.

\begin{theorem}\label{thm:blowup-argument}
Let $T$ and $\Gamma$ be as in Assumption \ref{assump:center-manifold-app}. If it holds
\begin{equation}\label{thm:blowup-argument:freqfunc}
    \lim_{r\to0}\mathrm{I}_{\phi, d^\pm}(N^\pm,r)>0, 
\end{equation} 
for at least one of the regions $\mathcal{C}^\pm$, then there exists a sequence $\rho_k \to 0$ as $k\to +\infty$ such that the sequence of pairs $(N_{\rho_k}^+, N_{\rho_k}^-)$ would converge in $\baseball{0}{1}$ locally strongly in $L^2$ to a $\left(Q-\frac{Q^\star}{2}\right)$ $\mathrm{ Dir}$-minimizer $(N^+_0, N^-_0)$ where $N^+_0:\basehalfball{+}{0}{1}\to \mathcal{A}_{Q}(\R^n)$ and $N^-_0:\basehalfball{-}{0}{1}\to \mathcal{A}_{Q-Q^\star}(\R^n)$, it holds that 
\begin{equation}\label{thm:blowup-argument:conv-energy}
\lim_{k\to\infty} \left(\int_{\basehalfball{+}{0}{R}} |DN_{\rho_k}^+|^2 + \int_{\basehalfball{-}{0}{R}} |DN_{\rho_k}^-|^2\right) = \int_{\basehalfball{+}{0}{R}} |DN_{0}^+|^2 + \int_{\basehalfball{-}{0}{R}} |DN_{0}^-|^2, \forall R\in(0,1),
\end{equation}
$(N^+_0, N^-_0)$ collapses at the interface $(T_0 \gamma, Q^\star\a{0})$, we have the following properties
\begin{enumerate}[\upshape (i)]
\item $(N^+_0, N^-_0)$ is nontrivial and in particular $\operatorname{Dir}(N^+_0,N^-_0,\baseball{0}{1}) = 1$; 
\item $\boldsymbol{\eta}\circ N^\pm_0 \equiv 0$.
\end{enumerate} 
\end{theorem}

As it is explained in Subsection \ref{two-dim-case}, Theorem \ref{thm:main-theorem} for currents of dimension $2$ follows from Theorem \ref{thm:collapsed-is-regular} which we are now able to prove for area minimizing currents of dimension $m\geq 2$, codimension $n\geq 2$ and multiplicity $Q^\star\geq 1$.

\begin{result}[Theorem \ref{thm:collapsed-is-regular}]
Let $T$ and $\Gamma$ be as in Assumption \ref{assumptions} with $C_0=0$. Then any two-sided collapsed point of $T$ is a two-sided regular point of $T$.
\end{result}
\begin{proof} Now, since we are under Assumption \ref{assump:center-manifold-app}, we can apply Theorem \ref{thm:ff-estimate-current} and we show that (b) of Theorem \ref{thm:ff-estimate-current} never occurs, then we are always in the case (a) of Theorem \ref{thm:ff-estimate-current} which ensures that $0$ is a boundary two-sided regular point of $T$. With this aim in mind observe that by the harmonic regularity of $\left(Q-\frac{Q^\star}{2}\right)$-$\mathrm{ Dir}$ minimizers which collapse at the interface, i.e., Theorem \ref{thm:harmonic-regularity-(Q-Q^*/2)-minimizing}, we have that $N_0^\pm = Q^\pm\a{h}$ for some classical $1$-valued harmonic function $h:\baseball{0}{1}\to\R^n$, hence we necessarily have
$$N^+_0 = Q \a{\boldsymbol{\eta} \circ N^+_0} \quad \mbox{and}\quad N^-_0 = (Q-Q^\star) \a{\boldsymbol{\eta} \circ N^-_0}.$$
By (ii) of Theorem \ref{thm:blowup-argument}, we have that $h\equiv 0$, but this is contradiction with $\operatorname{Dir}(N_0^+,N_0^-,\baseball{0}{1})=1$. Thus (b) of Theorem \ref{thm:ff-estimate-current} never occurs. This last fact surely completes the proof of the theorem.
\end{proof}

\section*{Acknowledgements}\label{ack}
\addcontentsline{toc}{section}{\nameref*{ack}}
We would like to thank Camillo De Lellis for proposing this problem, for his valuable comments, suggestions, and support throughout this work. This work was written during the Ph.D. of the second named author in the University of S\~{a}o Paulo under the supervision of prof. Camillo De Lellis and the first named author.\\
\indent S. Nardulli would like to thanks FAPESP for the grant ''Aux\'ilio \`a Pesquisa - Apoio a Jovens Pesquisadores" with financial code Fapesp: 2021/05256-0, and he also thanks CNPq for the grant "Bolsa de Produtividade em Pesquisa 1D" with financial code 312327/2021-8. \\
R. Resende would like to thanks the Coordena\c{c}\~{a}o de Aperfei\c{c}oamento de Pessoal de N\'{i}vel Superior - Brasil (CAPES) for the Ph.D. scholarship (Finance Code 88882.377954/2019-01). 

            \bibliographystyle{siam}
            \addcontentsline{toc}{section}{References} 
            \bibliography{biblio}
\end{document}